\documentclass[12pt]{amsart}

\usepackage{graphicx}
\usepackage{amsmath, comment}
\usepackage{amscd}
\usepackage{amsfonts}
\usepackage{amssymb}
\usepackage{color}
\usepackage{fullpage}
\usepackage{mathtools}
\usepackage[hypertexnames=false, colorlinks, citecolor=red,linkcolor=blue, urlcolor=black]{hyperref}
\usepackage{tikz-cd}
\usepackage[capitalise]{cleveref}
\usepackage{enumitem}

\DeclareMathOperator{\Ima}{Im}
\DeclareMathOperator{\Ker}{Ker}

\numberwithin{equation}{section}

\newtheorem*{theorem*}{Theorem}
\newtheorem{thm}{Theorem}[section]

\theoremstyle{plain}
\newtheorem{lem}[thm]{Lemma}

\newtheorem{defi}[thm]{Definition}
\newtheorem{cor}[thm]{Corollary}
\newtheorem{rema}[thm]{Remark}
\newtheorem{diag}[thm]{diagram}

\begin{document}
\author{Amin Soofiani}   
\address{Department of Mathematics, University of British Columbia, Vancouver, B.C., V6T
1Z2 Canada} 
\email{soofiani@math.ubc.ca}

\title{Hensel's lemma for the norm principle for spinor groups}

\begin{abstract}
   Let $K$ be a complete discretely valued field with residue field $k$
with $char \ k \ne 2$. Assuming that the norm principle holds for spinor groups $Spin(\mathfrak{h})$ for every regular skew-hermitian form $\mathfrak{h}$ over every quaternion algebra $\mathfrak{D}$ (with respect to the canonical involution on $\mathfrak{D}$) defined over any finite extension of $k$, we show that the norm principle holds for spinor groups $Spin(h)$ for every regular skew-hermitian form $h$ over every quaternion algebra $D$ (with respect to the canonical involution on $D$) defined over $K$.
\end{abstract} 

\maketitle

\setcounter{tocdepth}{1}
\tableofcontents

\section{Introduction} \label{103}

Let $L/K$ be a finite separable field extension, and $T$ be a commutative linear algebraic group defined over $K$. One can define the norm map $N_{L/K}: T(L) \longrightarrow T(K)$, which sends $t \mapsto \Pi_{\gamma} \gamma(t)$, where $\gamma$ runs over cosets of $Gal(K^{sep}/L)$ in $Gal(K^{sep}/K)$. For $T=\mathbb{G}_m$, this is the usual norm $N_{L/K}: L^* \longrightarrow K^*$.

Now we recall the definition of the norm principle for reductive linear algebraic groups, due to Merkurjev (see \cite{1}), which we call the $H^0$-variant of the norm principle. 

\begin{defi} \textbf{($H^0$-variant of the Norm Principle)} Let $L/K$ be a finite separable field extension. Suppose $G$ is a reductive linear algebraic group, and $T$ is a commutative linear algebraic group, both defined over $K$. Let $f:G \longrightarrow T$ be an algebraic group homomorphism defined over $K$. Consider the following diagram:

\[\ \begin{tikzcd}
G(L) \arrow[r, "f_L"] & T(L) \arrow[d, "N_{L/K}"] \\
G(K) \arrow[r, "f_K"] & T(K)                     
\end{tikzcd}
\]\

We say that the $H^0$-variant of the norm principle holds for $f:G \longrightarrow T$  over $L/K$, if $\Ima (N_{L/K} \circ f_L) \subseteq \Ima (f_K) $.

Also, we say that the $H^0$-variant of the norm principle holds for $f:G \longrightarrow T$, if for every finite separable field extension $L/K$ we have $\Ima (N_{L/K} \circ f_L) \subseteq \Ima (f_K) $.
\end{defi}

Note that if $G$ is commutative, then the $H^0$-variant of the norm principle holds for any map $G\to T$, because the norm homomorphism for $G$ makes the above diagram commutative. 

In \cite[Lemma~2.1]{1}, it is shown that the $H^0$-variant of the norm principle  for the map $f:G\to T$ described above, can be reduced to the $H^0$-variant of the norm principle  for the canonical map $G  \twoheadrightarrow G^{ab}$, where $G^{ab}$ is the abelianization of $G$, i.e. the quotient of $G$ by its commutator subgroup $G^{ab}=G/[G,G]$ (which is a torus).

\begin{defi}\label{254}
We say that the $H^0$-variant of the norm principle holds for a reductive group $G$
over $L/K$, if it holds for the map $G\longrightarrow G/[G,G]$ over $L/K$. Furthermore, we say that the $H^0$-variant of the norm principle holds for $G$
if the $H^0$-variant of the norm principle holds for $G$ over every finite separable field extension $L/K$.  
\end{defi}

The validity of the $H^0$-variant of the norm principle for reductive groups is an open question. In \cite[Theorem~1.1]{1}, it was proved that the $H^0$-variant of the norm principle holds in general for
all reductive groups of classical type without $D_n$ components.

For any central simple algebra with orthogonal involution $(A,\sigma)$, we denote the extended Clifford group of $(A,\sigma)$ by $\Omega(A,\sigma)$ (see \cite[Section~23.B]{13} for the definition of the extended Clifford group). By \cite[Proposition~5.1,~5.2, Example~4.4]{1}, in order to prove the $H^0$-variant of the norm principle for type $D_n$ groups, it suffices to prove the $H^0$-variant of the norm principle for the extended Clifford groups $\Omega(A,\sigma)$ for all central simple $K-$algebras $A$ of degree $2n$ with orthogonal involution $\sigma$.

A cohomological version of the norm principle for semisimple linear algebraic groups was introduced by Gille in \cite{29}, which we call the $H^1$-variant of the norm principle.

\begin{defi}
    \textbf{($H^1$-variant of the Norm Principle)} Let $L/K$ be a finite separable field extension, $S$ a semisimple linear algebraic group over $K$, and $Z \subseteq S$ a central subgroup. Consider the following diagram

\[
\begin{tikzcd}
{H^1(L,Z)} \arrow[r, "\alpha_L"] \arrow[d, "cor_{L/K}"'] & {H^1(L,S)} \\
{H^1(K,Z)} \arrow[r, "\alpha_K"]                         & {H^1(K,S)}
\end{tikzcd}
\]

In the above diagram, $cor_{L/K}$ denotes the corestriction map in Galois cohomology, and the map $\alpha$ is induced by the inclusion $Z \hookrightarrow S$. We say that the $H^1$-variant of the norm principle holds for the pair $(Z,S)$ over $L/K$ if for every $u \in \ker \alpha_L$, we have $cor_{L/K} (u)\in \ker  \alpha_{K}$.\
Also, we say that the $H^1$-variant of the norm principle holds for the pair $(Z,S)$, if for every finite separable field extension $L/K$, and for every $u\in\ker\alpha_L$, we have $cor_{L/K} (u)\in \ker  \alpha_{K}$.
\end{defi}
Note that if $S$ is commutative, then the $H^1$-variant of the norm principle holds for any pair $(Z,S)$, because the corestriction map for $H^1(-,S)$ makes the above diagram commutative. 

Let $C(S)$ denote the center of $S$. In Section \ref{115}, we will show that the $H^1$-variant of the norm principle for the pair $(C(S),S)$ implies the $H^1$-variant of the norm principle for the pair $(Z,S)$, where $Z$ is any central subgroup of $S$. Therefore, by "$H^1$-variant of the norm principle for $S$", we shall mean "$H^1$-variant of the norm principle for the pair  $(C(S),S)$". 

In Section \ref{115}, we will prove that the $H^0$-variant of the norm principle for reductive groups is equivalent to the  $H^1$-variant of the norm principle for semisimple groups (Theorem \ref{113}). As mentioned before, to settle the $H^0$-variant of norm principle for type $D_n$ groups, it suffices to show the $H^0$-variant of the norm principle for groups $\Omega(A, \sigma)$ for all central simple algebras of degree $2n$ with orthogonal involution $(A,\sigma)$. Under the equivalence between the $H^0$ and $H^1$ variants, the $H^0$-variant of the the norm principle for $\Omega(A, \sigma)$ is equivalent to the $H^1$-variant of the norm principle for $Spin(A, \sigma)$, where $ Spin(A, \sigma)$ is the spinor group of $(A,\sigma)$.

The $H^1$-variant of the norm principle was proved in full generality for all
semisimple groups defined over number fields by Gille (this will be discussed in Remark \ref{143}). The norm principle for semisimple groups defined over a field $K=k(C)$, function field of a curve $C$ defined over a number field $k$, is open in general. In particular it is open for the spinor groups $Spin(A,\sigma)$ defined over $K$, as mentioned before. In \cite{28}, Bhaskhar gave a scalar obstruction defined up to spinor norms, whose
vanishing implies the $H^1$-variant of the norm principle for $Spin(A,\sigma)$. In Subsection \ref{45} we will reformulate this obstruction using the $H^1$-variant of the norm principle. Let $V$ be the set of closed points of $C$, and for every $c\in V$, let $K_c$ be the completion of $K$ with respect to the discrete valuation associated to the point $c\in V$. Assume that the $H^1$-variant of the norm principle holds for $Spin(A,\sigma)$ over all completions $K_c$. It is an open question whether the kernel of the global to local map

$$\rho_{Spin(A,\sigma),V}: H^1(K, Spin(A,\sigma))\longrightarrow \prod_{c \in V} H^1(K_c, Spin(A,\sigma))$$ is trivial or not (In case it is trivial, we say that the local-global principle holds for $Spin(A,\sigma)$). If kernel of $\rho_{Spin(A,\sigma),V}$ is trivial, and if we further assume that the norm principle holds for $Spin(A,\sigma)$ over all completions $K_c$, then the norm principle holds for $Spin(A,\sigma)$ over $K$ (see \cite[Theorem~1.3.2]{270}). Equivalently, the failure of the norm principle  implies the failure of the local-global principle. Therefore, it is a natural question to ask whether the norm principle holds for $Spin(A,\sigma)$ over the completions $K_c$. In the above argument, one can  replace the set of discrete valuations associated to closed points by other discrete valuation sets, for instance all discrete valuations on $K$.

Note that each completion $K_c$ in the above argument is a complete discretely valued field whose residue field is a finite extension of $k$. Since we know that the norm principle holds for $Spin(A,\sigma)$ over number fields, a natural question to ask is the following:\\ 

\underline{Question:} Let $K$ be a complete discretely valued field with residue field $k$ with $char \ k \ne 2$. Assume that the $H^1$-variant of the norm principle holds for the spinor groups $Spin(A,\sigma)$ of all central simple $k-$algebras with orthogonal involution defined over all finite extensions of $k$. Then does the $H^1$-variant of the norm principle hold for the spinor groups $Spin(A,\sigma)$ of all central simple $K-$algebras with orthogonal involution defined over $K$?\\

An affirmative answer to this question would imply the norm principle for $Spin(A,\sigma)$ for the wide class of complete discretely valued fields $K_c$ arising from the local-global principle discussed above. When the answer to the above question is affirmative, we say that \underline{ Hensel's lemma holds for the norm principle}. In \cite{4}, an affirmative answer to this question was given in the case that the algebra $A$ splits, in which case $\sigma$ is the orthogonal involution adjoint to an even dimensional quadratic form over $K$.

 \begin{thm}\label{108}
  \cite[Theorem~5.1]{4} Let $K$ be a complete discretely valued field with residue field $k$ with
$char(k) \ne 2$. Assume that the $H^1$-variant of the norm principle holds for $Spin(q)$ for every regular quadratic form $q$ of even dimension defined over any finite extension of $k$. Then the $H^1$-variant of the norm principle holds for $Spin(Q)$ for every regular quadratic form $Q$ of even dimension over $K$.
 \end{thm}

 We remark that Theorem \ref{108} was stated in terms of the $H^0$-variant of the norm principle in \cite{4}. It is translated to the $H^1$-variant above. 
 
 In this paper, we generalize Theorem \ref{108} for spinor groups of skew-hermitian forms over quaternion algebras. We prove:
\begin{thm}\label{102} (Theorem \ref{40})
    Let $K$ be a complete discretely valued field with residue field $k$ with $char(k)\ne 2$.
    Assume that the $H^1$-variant of the norm principle holds for spinor groups $Spin(\mathfrak{h})$ for every regular skew-hermitian form $\mathfrak{h}$ over every quaternion algebra $\mathfrak{D}$ (with respect to the canonical involution on $\mathfrak{D}$) defined over any finite extension of $k$. Then, the $H^1$-variant of  norm principle holds for spinor groups $Spin(h)$ for every regular skew-hermitian form $h$ over every quaternion algebra $D$ (with respect to the canonical involution on $D$) defined over $K$.
\end{thm}

In particular, when the quaternion algebra $D$
is split, i.e. $D \cong M_2
(K)$, then $h$ corresponds to a quadratic form over $K$ by Morita equivalence. In this case, Hensel's lemma for the norm principle for $Spin(h)$ holds by Theorem \ref{108}. Thus our main result generalizes Theorem \ref{108}.

In Section \ref{115}, we will show that the $H^1$-variant of the norm principle for a semisimple group can be reduced to the $H^1$-variant of the norm principle for its simply connected covering (Lemma \ref{31}). Note that when $n$ is odd, then any absolutely simple, simply connected linear algebraic group of type $D_n$ is of the form $Spin(h)$ for a skew-hermitian form $h$ over a quaternion algebra (with respect to the canonical involution). So Theorem \ref{102} settles the norm principle for type $D_n$ groups over complete discretely valued fields $K$ (assuming that the norm principle is known over finite extensions of $k$) when $n$ is odd. This includes the case where $k$ is a number field (see Remark \ref{143}), or a $p-$adic field (in this case $H^1(k, G)$ is trivial by \cite{26}), or a global field of positive characteristic (in this case $H^1(k, G)$ is trivial by \cite{25}).

This paper contains 11 sections, including this introductory section. In Section \ref{115}, we prove the equivalence of the $H^0$ and $H^1$ variants of the norm principle. Furthermore, we reduce the $H^1$-variant of the norm principle to simply connected groups, and to the full center. We will then mainly work with the $H^1$-variant for the rest of the paper. All arguments in Section \ref{115} work in general for the $H^1$-variant of the norm principle for arbitrary semisimple linear algebraic groups (respectively, $H^0$-variant of the norm principle for arbitrary reductive linear algebraic groups). So we are not restricting ourselves to any specific type of groups in Section \ref{115}.  

From Section \ref{153} to the end of the paper, our focus will be on type $D_n$ groups. In Section \ref{153}, we recall two norm principles due to Gille and Knebusch, which will be of essential use in the rest of the paper. Then we give a new proof for a generalization of Scharlau's norm principle, which will be used in the proof of the main theorem. In Section \ref{135}, we reduce the norm principle for the spinor group $Spin(A,\sigma)$ from arbitrary finite separable field extensions to the case of separable quadratic field extensions. This reduction will simplify the proof of Hensel's lemma (our main result) significantly. Section \ref{135} also includes study of the obstruction for the norm principle for type $D_n$ groups from Galois cohomology point of view.

In Section \ref{109}, we reduce the norm principle for the spinor group of skew-hermitian forms over quaternion algebras to the case that the skew-hermitian form is anisotropic. This is a necessary hypothesis which we need in the proof of the main theorem.

All arguments in the first 5 sections of this paper work over arbitrary fields (of characteristic different from $2$), not necessarily over complete discretely valued fields. So the results can be used to study norm principle in general (in particular, groups of type $D_n$) over arbitrary fields (of characteristic different from $2$). From Section \ref{110}, we focus on complete discretely valued fields.  In Section \ref{110}, we recall theorems due to Springer and Larmour, which will be the main tools in the proof of Hensel's lemma for the norm principle. We set up the notations which are needed for the proof of the main theorem in Section \ref{111}. In Sections \ref{156} and \ref{157}, we prove some reductions of the main problem under certain hypotheses. Section \ref{128} contains a technical lemma which will be applied in several cases of the proof of our main result in the last section. Finally, in Section \ref{184}  we prove Hensel's lemma for the norm principle for the spinor group of skew-hermitian forms over quaternion algebras. \\

\textbf{Acknowledgements.} I am deeply grateful to Vladimir Chernousov and Sujatha Ramdorai for their insightful guidance and very helpful discussions.

\section{Equivalence of the two variants and reduction to simply connected groups}\label{115}

\subsection{Equivalence of $H^0$ and $H^1$ variants of the norm principle} In this subsection, we show that the two variants of the norm principle, the $H^0$-variant for reductive groups and the $H^1$-variant for semisimple groups, are equivalent. This will allow us to work with the $H^1$-variant for the rest of the paper, where we can use tools from Galois cohomology of semisimple linear algebraic groups to study the norm principle for type $D_n$ groups.

\begin{thm}\label{113}
Let $L/K$ be a finite separable field extension. Then
the following statements are equivalent:

\item \textbf{1.} $H^0$-Variant of the norm principle holds over $L/K$ for all reductive groups $G$ defined over $K$.

\item \textbf{2.} $H^1$-Variant of the norm principle holds over $L/K$ for all pairs $(J,H)$, where $H$ is a semisimple group define over $K$ and $J$ a central subgroup in $H$.\\

Note that one can replace "over $L/K$" in the statement of the theorem by "over $K$"; i.e. if (1) is true for all finite separable field extensions $L/K$, then (2) is true for all finite separable field extensions $L/K$ and vice versa. 

\begin{proof}
(1 $\Rightarrow$ 2):\\

Let $H$ be a semisimple group over $K$ and $J\subseteq H$ be a central subgroup. We have the following diagram and we need to show if $u\in \ \Ker (\alpha_L)$, then $cor_{L/K} (u) \in \ \Ker(\alpha_K)$.

\[\
\begin{tikzcd}
{H^1(L,J)} \arrow[r, "\alpha_L"] \arrow[d, "cor_{L/K}"'] & {H^1(L,H)} \\
{H^1(K,J)} \arrow[r, "\alpha_K"]                         & {H^1(K,H)}
\end{tikzcd}
\]\

Recall that there is an anti-equivalence between the groups of multiplicative type (see \cite[Chapter~7, Section~2]{24}) and the abelian groups on which $\Gamma = Gal(K^{sep}/K)$ acts continuously (see \cite[Chapter~7, Section~3]{24}). Under this anti-equivalence, we can find an embedding $J \hookrightarrow T'$ where $T'$ is a quasi-trivial torus over $K$. This can be done in the following way: Let $A$ be the character group of $J$, and $a_1, a_2, \dots, a_n$ be  generators of $A$ as a $\Gamma$ - module. For $1\leq i \leq n$, let $\Delta_i$ be the subgroup of finite index in $\Gamma$ which acts trivially on $a_i$. Then we have a surjection $\bigoplus\limits_{i=1}^n \mathbb{Z} [\Gamma / \Delta_i] \longrightarrow A$, which corresponds to an embedding $J\hookrightarrow T'$ that we were looking for.

Consider the pushout of the morphisms $J\longrightarrow H$ and $J\longrightarrow T'$:

\[\
\begin{tikzcd}
             &  & 1 \arrow[dd]                  &  & 1 \arrow[dd]            &  &                                &  &   \\
             &  &                               &  &                         &  &                                &  &   \\
1 \arrow[rr] &  & J \arrow[dd] \arrow[rr]       &  & H \arrow[rr] \arrow[dd] &  & G' \arrow[rr] \arrow[dd, "id"] &  & 1 \\
             &  &                               &  &                         &  &                                &  &   \\
1 \arrow[rr] &  & T' \arrow[dd] \arrow[rr]      &  & G \arrow[rr] \arrow[dd] &  & G' \arrow[rr]                  &  & 1 \\
             &  &                               &  &                         &  &                                &  &   \\
             &  & T \arrow[dd] \arrow[rr, "id"] &  & T \arrow[dd]            &  &                                &  &   \\
             &  &                               &  &                         &  &                                &  &   \\
             &  & 1                             &  & 1                       &  &                                &  &  
\end{tikzcd}
\]\

Note that $G$ is reductive under our construction.
Since $T^{\prime}$ is a quasi-trivial torus, $H^1(L,T^{\prime})=H^1(K,T^{\prime})=1$. This gives rise to the following commutative diagram where  $s$ and $s^{\prime}$ are the connecting maps of the first and second columns of the above diagram, respectively, and $\alpha$ is induced by the map $J \hookrightarrow H $.

%\begin{diag}\label{3}

\[ \begin{tikzcd}[row sep=3.5em]
{}  &&
  G(L) \arrow[dd,swap,"f_L" near start]  \\
& {}  &&
  G(K) \arrow[dd,"f_K"] \\
T(L) \arrow[dd, "s_L"] \arrow[rr,"id" near end] \arrow[dr,swap,"N_{L/K}"] && T(L) \arrow[dd, "s^{\prime}_L" near end] \arrow[dr,swap,"N_{L/K}"] \\
& T(K) \arrow[dd, "s_K" near end] \arrow[rr,"id", near start] && T(K) \arrow[dd, "s^{\prime}_K"] \\ 
H^1(L,J) \arrow[rr,"\alpha_L"near start] \arrow[dr,swap,"cor_{L/K}"] \arrow[dd,swap] &&
  H^1(L,H)  \\
& H^1(K,J) \arrow[dd,swap] \arrow[rr,crossing over,"\alpha_K" ] &&
  H^1(K,H)  \\
H^1(L,T')   && {} \\
& H^1(K,T') && {}
\end{tikzcd}\]
%\end{diag}

Let $v=cor_{L/K} (u)$. Since $H^1(L,T^{\prime})=1$, $u\in \Ima (s_L)$. Let $l$ be an element in $T(L)$ such that $s_L(l)=u$. Let $w=N_{L/K}(l)$. By the commutativity of the above diagram, $l \in \ \Ker (s^{\prime}_L)=\Ima (f_L)$. By assumption, the $H^0$-variant of the norm principle holds for $G$ over $L/K$; as mentioned in the introduction, this implies that the $H^0$-variant of the norm principle holds for the map $f: G \to T$ over $L/K$, hence we have $w \in \Ima (f_K)=\ \Ker (s^{\prime}_K)$, therefore $\alpha_K(v)= s^{\prime}_K (w)=0 $ and $v \in \ \Ker (\alpha_K)$.\\ 

(2 $\Rightarrow$ 1):

Let $H=[G,G]$, $T=G^{ab}$, $T^{\prime}=R(G)$ (the radical of $G$), $J=H\cap T^{\prime}$, $G^{\prime}=\frac{G}{R(G)}$. Then $H$ is a semisimple group and $J \subseteq H$ is a central subgroup; with this set-up, the two diagrams above will again be commutative with exact rows and columns. 

Let $s \in G(L)$, $l=f_L(s)$, $w=N_{L/K}(l)$. We would like to show $w \in \Ima (f_K)$. Let $u=s_L(l)$, $v=cor_{L/K} (u)$. So we have $s^{\prime}_L (l)=0$ hence $u \in \ \Ker (\alpha_L)$. By assumption, the $H^1$-variant of the norm principle holds for the pair $(J,H)$ over $L/K$, which implies that $\alpha_K (v)=0$, therefore $w\in \ \Ker (s^{\prime}_K)= \Ima (f_K)$.
\end{proof}
\end{thm}

\subsection{Reduction to the center, and to simply connected covers} The next lemma states that in order to prove the $H^1$-variant of the norm principle for a pair $(Z_1,G)$ it is enough to show the $H^1$-variant of the norm principle for $(Z_2,G)$, where $Z_2$ is a central subgroup of $G$ containing $Z_1$.

\begin{lem}\label{37}
    Let $L/K$ be a finite separable field extension and $G$ be a semisimple linear algebraic group defined over $K$. Let   $Z_1 \subseteq Z_2 \subseteq G$ be central subgroups. Then the $H^1$-variant of the norm principle for the pair $(Z_2, G)$ over $L/K$ implies the $H^1$-variant of the norm principle for the pair $(Z_1,G)$ over $L/K$.

    \begin{proof}

Let $G_1:= G/Z_1$ and $G_2:= G/Z_2$. We have the following diagram where $i$ is the identity map:

\[\
\begin{tikzcd}
G_1(-) \arrow[d, "f"] \arrow[r] & {H^1(-, Z_1)} \arrow[d, "g"] \arrow[r] & {H^1(-, G)} \arrow[d, "i"] \\
G_2(-) \arrow[r]                & {H^1(-,Z_2)} \arrow[r]                 & {H^1(-,G)}                
\end{tikzcd}
\]\
        
It gives rise to the following commutative diagram whose rows are exact:

\[\  \begin{tikzcd}[sep=2em]
G_1(L) \arrow[rr, "\alpha_L"] \arrow[dr,"f_L"]  && H^1(L, Z_1) \arrow[dr,"g_L"]  \arrow[dd, "cor_{L/K}" near end]  \arrow[rr,"\beta_L" near end]    && H^1(L, G) \arrow[dr, "i_L"] \\
& G_2(L)  \arrow[rr,"\gamma_L" near end]  && H^1(L, Z_2 )  \arrow[dd, "cor_{L/K}" near end]  \arrow[rr,"\theta_L", near start] && H^1(L, G)  \\ 
G_1(K)  \arrow[dr, "f_K"]  \arrow[rr, "\alpha_K" near end] && H^1(K, Z_1)  \arrow[dr,"g_K"]  \arrow[rr, "\beta_K", near end] &&
  H^1(K, G) \arrow[dr, "i_K"] \\
& G_2(K) \arrow[rr, "\gamma_K"] && H^1(K, Z_2)   \arrow[rr,crossing over,"\theta_K" ] &&
  H^1(K, G)  
  \end{tikzcd}
  \]\ 
%\end{diag}

Let $u\in \Ker (\beta_L)= \Ima  (\alpha_L)$. Choose  $v\in {\alpha_L}^{-1}(u)$, and put $w=\gamma_L(f_L(v))$. By assumption, $cor_{L/K}(w) \in \Ima (\gamma_K)=\Ker (\theta_K)$. Hence $cor_{L/K}(u) \in \Ker (\beta_K)$ and the $H^1$-variant of the norm principle holds for the pair $(Z_1, G)$ over $L/K$.
    \end{proof}
\end{lem}
From now on, by the \underline{norm principle for $G$} (where $G$ is a semisimple group), we shall mean the \underline{$H^1$-variant of the norm principle for the pair $(C(G),G)$}, where $C(G)$ is the center of $G$, unless otherwise stated. By Lemma \ref{37}, it is enough to show the $H^1$-variant of the norm principle only for the full center of any semisimple group.

Now we show that the norm principle for a simply connected covering of a semisimple group implies the norm principle for the group itself.

 \begin{lem}
     \label{31}
Let $H$ be a semisimple group defined over $K$, $L/K$ be a finite separable field extension, and $\phi: \Tilde{H}\longrightarrow H$ be the simply connected covering of $H$. Then the norm principle for $\Tilde{H}$ implies the norm principle for $H$ (over $L/K$).
     \begin{proof}
     Let $Z$ and $\Tilde{Z}$ be the centers of $H$ and $\Tilde{H}$, respectively, and $G=H/Z= \Tilde{H}/ \Tilde{Z}$.   The following sequence (which is exact in rows) 

\[\
\begin{tikzcd}
1 \arrow[r] & \tilde{Z} \arrow[r] \arrow[d] & \tilde{H} \arrow[r] \arrow[d, "\phi"] & G \arrow[r] \arrow[d, "id"] & 1 \\
1 \arrow[r] & Z \arrow[r]                   & H \arrow[r]                           & G \arrow[r]                 & 1
\end{tikzcd}
\]\

 gives rise to the following commutative diagram whose rows are exact:

%\begin{diag}
\[ \begin{tikzcd}[sep=2em]
G(L) \arrow[rr, "\alpha_L"] \arrow[dr,"id"]  && H^1(L, \Tilde{Z}) \arrow[dd, "cor_{L/K}" near end]  \arrow[rr,"\beta_L" near end]    \arrow[dr,"j_L"] && H^1(L, \Tilde{H}) \arrow[dr, "\tau_L"] \\
& G(L)  \arrow[rr,"\gamma_L" near end]  && H^1(L, Z )  \arrow[dd, "cor_{L/K}" near end]  \arrow[rr,"\theta_L", near start] && H^1(L, H)  \\ 
G(K)  \arrow[dr, "id"]  \arrow[rr, "\alpha_K" near end] && H^1(K, \Tilde{Z})  \arrow[rr, "\beta_K", near end] \arrow[dr, "j_K"]  &&
  H^1(K, \Tilde{H}) \arrow[dr, "\tau_K"] \\
& G(K) \arrow[rr, "\gamma_K"] && H^1(K, Z)   \arrow[rr,crossing over,"\theta_K" ] &&
  H^1(K, H)  
\end{tikzcd}\]
%\end{diag}

Let $v\in \Ker (\theta_L)$. Choose  $g\in {\gamma_L}^{-1}(v)$ and put $u=\alpha_L(g)$. By assumption, $cor_{L/K}(u) \in \Ker (\beta_K)=Im  (\alpha_K)$, therefore there exists $w\in G(K)$ such that $\alpha_K(w)=cor_{L/K}(u)$. Hence $cor_{L/K}(v)=\gamma_K(w) \in \Ima (\gamma_K)=\Ker (\theta_K)$ and the norm principle holds for $(Z,H)$.
     \end{proof}
 \end{lem}

\begin{rema}\label{114}
    By Lemmas \ref{37} and \ref{31}, in order to prove the $H^1$-variant of the norm principle for a pair $(Z,G)$, it is enough to prove it for $(C(\Tilde{G}), \Tilde{G})$, where $\Tilde{G}$ is the simply connected covering of $G$, and $C(\Tilde{G})$ is its center.
\end{rema}

\section{Gille's, Knebusch's, and Scharlau's norm principles}\label{153}

All fields in the rest of the paper are assumed to be of characteristic different from $2$.

\subsection{Rationality, Gille's norm principle, and Knebusch's norm principle}\label{456}

The study of the norm principle is closely related to the study of $R$-equivalence classes of linear algebraic groups. The notion of $R$-equivalence was first introduced by Manin. Although $R$-equivalence is defined in general for algebraic varieties, we will restrict our attention to the underlying varieties of connected linear algebraic groups, where one knows some additional structural results. We briefly recall its definition below.

Let $G$ be a connected linear algebraic group defined over a field $K$. Then $x, y \in G(K)$ are
said to be $R$-equivalent if there exists a $K$-rational map 
${\mathbb{A}^1}  -\!\! \rightarrow G$ defined at $0$ and $1$, which sends $0$ to $x$ and $1$ to $y$. This defines an
equivalence relation on $G(K)$ (see \cite[Definition~6.2~and~Lemma~6.3]{14}). The set of all points in $G(K)$ which are $R-$equivalent to identity is denoted by $RG(K)$. This set is in fact a normal subgroup of $G(K)$, and there is a bijection of sets between $G(K)/R$ and $G(K)/RG(K)$ \cite[Section~6.2, Lemma~6.3]{14}. This bijection induces a group structure on the set
of $R$-equivalence classes of $G(K)$. We  consider $G(K)/R$ along with this induced group structure. The group $G$ is said to be \textbf{R-trivial} if for all field extensions $F/K$, the group $G(F)$ has only one $R$-equivalence class, i.e., $RG(F)=G(F)$. 

We also have the notion of rationality: a linear algebraic group defined over $K$ is called  $K-$rational if its function field is purely transcendental over $K$. If $G$ is $K-$rational, then it is $R-$trivial (see \cite[Section~6.1,~Page 71,~and~Lemma~6.5]{14}).

Let $\beta: \Tilde{G} \longrightarrow G$ be a central $K-$isogeny of reductive algebraic groups defined over a field $K$ with smooth kernel $Z$. We have the exact sequence of the Galois cohomology pointed sets

$$1 \longrightarrow Z(K) \longrightarrow  \Tilde{G}(K) \overset{\beta_K}\longrightarrow G(K) \overset{\phi_K}\longrightarrow H^1(K,Z) \overset{e_K}\longrightarrow H^1(K, \Tilde{G}). $$

\begin{thm}\label{7}
    (\textbf{Gille}, \cite[Theorem~A]{8}) For any finite separable field extension $L/K$, we have  $cor_{L/K}(\phi_L (RG(L))) \subseteq \phi_K(RG(K)).$
\end{thm}

%Assume that the variety $G$ is $K-$rational. Then a theorem due to Gille immediately implies that the norm principle holds for the pair $(Z,\Tilde{G})$ over $K$.

If $G$ is $K-$rational, then it is in  particular $R-$ trivial. Therefore Gille's theorem implies that for any finite separable field extension $L/K$, we have $cor_{L/K}(Im \   \phi_L) \subseteq Im \ \phi_K$, or equivalently, we have  $cor_{L/K}(\ker e_L) \subseteq \ker e_K$, i.e., the norm principle holds for $(Z,\Tilde{G})$ over $L/K$. Since this is the case for every finite separable field extension $L/K$, then the norm principle holds for $(Z,\Tilde{G})$ over $K$.

Now we recall Knebusch's norm principle from the algebraic theory of quadratic forms, and we then revisit it using Gille's theorem. Let $q$ be an even dimensional quadratic form defined over $K$. We say that an element $u\in K^*$ is a \emph{spinor norm} of $q$ over $K$, if $u$ is product of an even number of values in $K^*$ represented by $q$.

\begin{thm}\label{147}
 (\textbf{Knebusch's Norm Principle})\  
 
Let $L/K$ be a finite separable field extension, and $q$ be an even dimensional regular quadratic form defined over $K$.   Let $x\in L^*$ be a spinor norm of $q$ over $L$. Then $N_{L/K}(x)$ is a spinor norm of $q$ over $K$.
 \begin{proof}
  See \cite[Chapter~7, Section~5, Page~206]{14}.
 \end{proof}
\end{thm}

Theorem \ref{147} can be formulated as an $H^0$-variant of the norm principle. Let $\Gamma^+(q)$ be the \emph{special Clifford group} of $q$ defined over $K$ (see \cite[Section~3.4.2]{3}). Then Theorem \ref{147} can be formulated as follows: the $H^0$-variant of the norm principle holds for the \emph{spinor norm map} $\Gamma^+(q)\to\mathbb{G}_m$ (see \cite[Section~8.2~and~Section~3.4.4]{3}). 

Now consider the following sequence of algebraic groups over $K$, which is exact over the separable closure of $K$:
$$1 \longrightarrow \mu_2 \longrightarrow Spin (q) \longrightarrow O^+(q) \longrightarrow 1,$$  where $Spin(q)$ is the \emph{spinor group} of $q$ and $O^+(q)$ is the \emph{special orthogonal group} of $q$. The above sequence induces the following map in Galois cohomology $$H^1(K,\mu_2) \longrightarrow H^1(K, Spin(q)).$$

The kernel of this map exactly consists of (square classes of) spinor norms (see \cite[Page~187]{13}). In the proof of Theorem \ref{113}, if we consider $G=\Gamma^+(q)$, $H=Spin(q)$, and $J=\mu_2=ker \ (Spin(q) \to O^+(q))$, then we conclude that the $H^0$-variant of the norm principle for the spinor norm map $\Gamma^+(q)\to \Gamma^+(q)/ Spin(q)=\mathbb{G}_m$ is equivalent to the $H^1$-variant of the norm principle for the pair $(\mu_2, Spin(q))$. 

$$\textit{Knebusch's Norm Principle} \iff$$ $$ H^0-\textit{variant of the norm principle for the spinor norm map} \ \Gamma^+(q)\to \mathbb{G}_m \iff $$  $$H^1-\textit{variant of the norm principle for the pair}  \ (\mu_2, Spin(q)).$$

In fact, Knebusch's norm principle can also be derived from Theorem \ref{7}: since the variety $O^+(q)$ is rational (see \cite[Proposition~2.4]{6}), then Theorem \ref{7} combined with the exact sequence 

$$1 \longrightarrow \mu_2 \longrightarrow Spin (q) \longrightarrow O^+(q) \longrightarrow 1$$ implies Knebusch's norm principle. 

Recall that for any even-dimensional central simple \( K \)-algebra \( A \)
equipped with an orthogonal involution \( \sigma \), one can define the
\emph{spinor group} of \( (A,\sigma) \), denoted \( Spin(A,\sigma) \)
(see \cite[Definition~13.30]{13}), as well as the
\emph{special orthogonal group} of \( (A,\sigma) \), denoted
\( O^{+}(A,\sigma) \) (see \cite[Definition~12.24]{13}). If $A$ is split over $K$, then by Morita equivalence the group  \( Spin(A,\sigma) \) (resp. \( O^{+}(A,\sigma) \)) coincides with the usual spinor group (resp. special orthogonal group) of the quadratic form on $K$ adjoint to $\sigma$.

Considering the exact sequence$$1 \longrightarrow \mu_2 \longrightarrow Spin (A,\sigma) \longrightarrow O^+(A,\sigma) \longrightarrow 1$$in conjunction with the rationality of the  variety $O^+(A,\sigma)$ (see \cite[Proposition~2.4]{6}) and Theorem \ref{7}, one obtains a generalization of Knebusch's norm principle:

\begin{thm}\label{8} (Generalization of Knebusch's norm principle, due to Gille) For any even dimensional central simple $K-$algebra with orthogonal involution $(A,\sigma)$, the $H^1$-variant of the norm principle holds for the pair $(\mu_2, Spin(A,\sigma))$ over $K$.
\end{thm}

We will use the term "spinor norm" not only for quadratic forms, but in general for elements in the kernel of the map $H^1(K,\mu_2)\to H^1(K,Spin(A,\sigma))$ for any central simple $K-$algebra with orthogonal involution $\sigma$.

\begin{rema}\label{143}
    It is known that for any semisimple adjoint linear algebraic group of classical type $G$ over a number field $K$, the group $G(K)/R$,  the group of $R$-equivalence classes of the
$K$-points of $G$, is trivial (\cite[Corollaire III.4.2]{8} and \cite[page~1]{1000}). Therefore, by Theorem \ref{7}, the norm principle holds for all classical semisimple linear algebraic groups over number fields.
\end{rema} 

\subsection{Scharlau's norm principle} In this subsection, we recall Scharlau's norm principle from the algebraic theory of quadratic forms, and then we prove a generalization of it for central simple algebras with orthogonal involution.

Let $q$ be a regular $n$-dimensional quadratic form defined over $K$. We say that an element $\lambda\in K^*$ is a \emph{multiplier} for $q$ over $K$, if we have a $K-$isometry $\lambda q \cong_K q$.

\begin{thm}\label{1470}
 (\textbf{Scharlau's Norm Principle})\  
Let $L/K$ be a finite separable field extension, and $q$ be a regular quadratic form defined over $K$.  Suppose $x \in L^{*}$ be a multiplier for $q$ over $L$. Then $N_{L/K}(x)\in K^*$ is a multiplier for $q$ over $K$.
 \begin{proof}
  See     \cite[Chapter~7, Section~4, Page~205]{14}.
 \end{proof}
\end{thm}

Note that even though the above theorem holds for quadratic forms of arbitrary dimension, but in the rest of this subsection we assume that the dimension of $q$ is even.  

Let $GO(q)$ be the group of \emph{similitudes} of $q$, which is an algebraic group defined over $K$ (see \cite[Section~12.A]{13} for a definition). Theorem \ref{1470} is equivalent to the $H^0$-variant of the norm principle for the \emph{multiplier map $GO(q)\to \mathbb{G}_m$}  (See \cite[Section~3.2]{3} for the definition of the multiplier map, and \cite[Section~8.2]{3} for the equivalent formulation of Theorem \ref{1470} in terms of the $H^0$-variant of the norm principle). 

 The subgroup of $GO(q)$ consisting of proper similitudes is denoted by $GO^+(q)$ (see \cite[Definition~12.12]{13}). Note that by \cite[Section~10.4]{3} (considering the case that the central simple algebra in the argument of \cite[Section~10.4]{3} splits), the $H^0$-variant of the norm principle also holds for the restriction of the multiplier map to proper similitudes  $GO^+(q) \to \mathbb{G}_m$.

The semisimple part of $GO^+(q)$ is the special orthogonal group of $q$, i.e. $O^+(q)$. In the proof of  Theorem \ref{113}, let $G=GO^+(q)$, $H=O^+(q)$, and $J=\mu_2$ (center of $O^+(q)$). Then, we conclude that the $H^0$-variant of the norm principle for the multiplier map $GO^+(q)\to GO^+(q)/ O^+(q)=\mathbb{G}_m$ is equivalent to the $H^1$-variant of the norm principle for the pair $(\mu_2, O^+(q))$. 

$$\textit{Scharlau's Norm Principle (for multipliers of proper similitudes)} \iff$$ $$ H^0-\textit{variant of the norm principle for the multiplier map} \ GO^+(q)\to \mathbb{G}_m \iff $$  $$H^1-\textit{variant of the norm principle for the pair}  \ (\mu_2, O^+(q)).$$\\

 In the rest of this subsection, we generalize Scharlau's norm principle. Consider a central division algebra with symplectic involution $(D,\sigma)$ over a field $K$, and an $n-$dimensional skew-hermitian form $h$ over $(D,\sigma)$. Recall that the skew-hermitian form $h$ induces the adjoint orthogonal involution $\sigma_h$ on $M_n(D)$ \cite[Page~43, Theorem~4.2]{13}. Then the special orthogonal group of $h$ is defined as $O^+(h):= O^+(M_n(D),\sigma_h)$, and the spinor group of $h$ is defined as $Spin(h):= Spin(M_n(D),\sigma_h)$. 
We also have the projective   group of proper similitudes of $h$: $PGO^+(h):=PGO^+(M_n(D),\sigma_h)$ (see \cite[Definition~12.24]{13}).  We have the following commutative exact diagram:

\[\
\begin{tikzcd}
            & 1 \arrow[d]                                  & 1 \arrow[d]                        &                                    &   \\
1 \arrow[r] & \mu_2 \arrow[d] \arrow[r]                    & O^+(h) \arrow[d] \arrow[r]         & PGO^+(h) \arrow[r] \arrow[d, "id"] & 1 \\
1 \arrow[r] & \mathbb{G}_m \arrow[d, "\times 2"] \arrow[r] & GO^+(h) \arrow[d, "\mu"] \arrow[r] & PGO^+(h) \arrow[r]                 & 1 \\
            & \mathbb{G}_m \arrow[d] \arrow[r, "id"]       & \mathbb{G}_m \arrow[d]             &                                    &   \\
            & 1                                            & 1                                  &                                    &  
\end{tikzcd}
\]\

Recall that the center of $O^+(
h)$ is $\mu_2$, and $\mu$ is the multiplier map (see  \cite[Section~10.4]{3}). The proof of Theorem \ref{113} shows that the $H^0$-variant of the norm principle for the multiplier map $\mu$ is equivalent to the $H^1$-variant of the norm principle for the pair $(\mu_2, O^+(h))$. Let $\beta_K$ be the induced Galois cohomology map $\beta_K: H^1(K,\mu_2)\longrightarrow H^1(K,O^+(h))$. Recall that the elements of $H^1(K,O^+(h))$ are in one-to-one correspondence with isometry classes of skew-hermitian forms with the same dimension and the same discriminant as $h$, and the distinguished element of $H^1(K,O^+(h))$ is the $K-$isometry class of $h$ (see \cite[Section 29.D]{13}).
Let $\eta$ be an element in $H^1(K,\mu_2)={K^*} / {{K^{*}}^2}$ and choose a representative for $\eta$, say $\eta=[a]=a{K^{*}}^2 \in {K^*} / {{K^{*}}^2}$ for some $a\in K^*$. The map $\beta_K$ sends the square class $\eta=[a]$ to the $K-$isometry class of the skew-hermitian form $ah$. Assume that $\eta \in \Ker \beta_K$. Then $ah$ is isometric to the skew-hermitian form $h$ over $K$, hence the element $a$ will be a multiplier for $h$ over $K$. Conversely, if $a$ is a multiplier for $h$ over $K$, i.e. if $ah\cong h$ over $K$, then the element $\eta$ will be in the kernel of the map $\beta_K$. We will show the norm principle for the pair $(\mu_2, O^+(h))$ over any finite separable field extension $L/K$, i.e., $cor_{L/K}(ker \ \beta_L) \subseteq ker \ \beta_K$.

The next theorem is related to the hyperbolicity of skew-hermitian forms over central division algebras and over the function fields of their Severi-Brauer varieties. Theorem \ref{4} will be used to generalize Scharlau's norm principle (see Theorem  \ref{5}). 

\begin{thm}\label{4}
Let $(D, \sigma)$ be a central division algebra with symplectic involution and $h$ a skew-hermitian form over $(D, \sigma)$. Let $E$ be the function field of $SB(D)$ (the Severi-Brauer variety over $D$). If $h$ is hyperbolic over $D_E$, then it is hyperbolic over $D$. 

\begin{proof}
See \cite[Theorem~1.1]{11}.
\end{proof}
\end{thm}

\begin{thm}\label{5}\textbf{(Generalization of Scharlau's norm principle)}
Let $(D, \sigma)$ be a central division algebra with symplectic involution $\sigma$ over $K$, and $h$ be a skew-hermitian form over $(D, \sigma)$. Then the norm principle holds for  $O^+(h)$ over any finite separable field extension $L/K$.

\begin{proof}

Let $E$ be the function field of $SB(D)$, the Severi-Brauer variety of $D$. Consider the following commutative diagram, whose  arrows are inclusions and $LE:=L\underset{K}{\otimes}E$:

\[\begin{tikzcd}
L \arrow[r] 
& LE \\
K \arrow[r] \arrow[u] & E \arrow[u]
\end{tikzcd}
\]

Suppose the dimension of $h$ over $D$ is $d$.
Let $V:=D^d$, $D_L:=D\underset{K}{\otimes}L$ , $D_E:=D\underset{K}{\otimes}E$, $D_{LE}:=D\underset{K}{\otimes}LE$, $V_L:=V\underset{K}{\otimes}L$,
$V_E:=V\underset{K}{\otimes}E$, and
$V_{LE}:=V\underset{K}{\otimes}{LE}$. Also let $h_{L}$, $h_{E}$, and $h_{LE}$ be the corresponding skew-hermitian forms over $V_L$, $V_E$, and $V_{LE}$, respectively.

Let $\beta_L$ be the induced Galois cohomology map $\beta_L: H^1(L,\mu_2)\longrightarrow H^1(L,O^+(h))$. Let $\eta$ be an element in $H^1(L,\mu_2)={L^*} / {{L^{*}}^2}$ and choose a representative for $\eta$, say $a {L^{*}}^2 \in {L^*} / {{L^{*}}^2}$ for some $a\in L^*$. Assume that $\eta\in \Ker \beta_L$, so $ah_L\cong h_L$. Hence $h_{L} \perp -ah_{L} $ is hyperbolic, and so is $h_{LE} \perp -ah_{LE}$. Since $SB(D)$ is a Severi-Brauer variety, $D_E$ splits, i.e. $D_E\cong M_n(E)$ for some $n \in \mathbb{N}$. We also have  $D_{LE}\cong M_n(LE)$. By Morita equivalence (\cite[Chapter 1, Section 9]{12}), skew-hermitian forms over $D_E$ (respectively $D_{LE}$) correspond to quadratic forms over $E$ (respectively $LE$). So we get a quadratic form $f_E$ (respectively $f_{LE}$). Under this correspondence, direct sums and hyperbolicity is preserved. Therefore, $h_{LE} \perp -ah_{LE}$ corresponds to a hyperbolic quadratic form $f_{LE} \perp -a f_{LE}$, which means  $a$ is a multiplier of $f_{LE}$. By Scharlau's norm principle, $f_E \perp -N_{LE/E}(a) f_{E}$ is hyperbolic, and by Theorem \ref{4}, $h \perp -N_{L/K}(a)h$ is hyperbolic over $K$. So $N_{L/K} (a)$ is a multiplier of $h$ over $K$, and $cor_{L/K} (\eta) \in \ \Ker (\beta_K)$. 
\end{proof}

\end{thm}

Note that if $D$ splits in Theorem \ref{5}, then one obtains Scharlau's norm principle for quadratic forms, by using Morita equivalence.

\section{Obstruction to the norm principle for the spinor group, and reduction to quadratic field extensions}\label{135}

Assume that $A$ is a central simple $K$-algebra with an orthogonal involution $\sigma$. In this section, we study the obstruction to the norm principle for the spinor group $Spin(A,\sigma)$. We also reduce the norm principle for the group $Spin(A, \sigma)$ from arbitrary finite separable field extensions to quadratic field extensions. 

The spinor group $Spin (A,\sigma)$ is a semisimple group of type $D_n$, whose center will be denoted by $Z$. 
Let $F:= K(\sqrt{disc   \ \sigma})$ (see \cite[Section~7, Definition~7.2]{13} for the definition of $disc \ \sigma$). Then $Z=\mathbf{R}_{F/K}(\mu_2)$ if $n$ is even, and $Z=\Ker \ (\mathbf{R}_{F/K} (\mu_4)\overset{Norm}\longrightarrow \mu_4)$ if $n$ is odd
 (see \cite[Page~332]{19}).

\subsection{Obstruction to the norm principle for $Spin(A,\sigma)$} \label{45} We describe the obstruction to the norm principle for $Spin(A,\sigma)$ in terms of spinor norms. The obstruction was studied by Bhaskhar in \cite{28}, but we study the obstruction from Galois cohomology point of view.

Consider the following commutative diagram whose vertical arrows are inclusions:

\[\
\begin{tikzcd}
1 \arrow[r] & \mu_2 \arrow[r] \arrow[d] & Z \arrow[r] \arrow[d]      & \mu_2 \arrow[d] \arrow[r] & 1 \\
1 \arrow[r] & \mu_2 \arrow[r]           & {Spin(A,\sigma)} \arrow[r] & {O^+(A,\sigma)} \arrow[r] & 1
\end{tikzcd}
\]\

One of the groups $\mu_2$ is center of $O^+(A,\sigma)$, and the other $\mu_2$ is the kernel of the map $Spin(A,\sigma) \to O^+(A,\sigma)$. 

Let $L/K$ be a finite separable field extension. The above diagram gives rise to:

%\begin{diag}\label{10}
\[ \begin{tikzcd}[sep=0.7em]
 H^1(L, \mu_2) \arrow[rr, "f_L"] \arrow[dr,"id"]  \arrow[dd, "cor_{L/K}"]&& H^1(L, Z) \arrow[dd, "cor_{L/K}" near end] \arrow[rr,"g_L" near end]    \arrow[dr,"\alpha_L"] && H^1(L, \mu_2) \arrow[dr, "\beta_L"] \arrow[dd, "cor_{L/K}" near end] \\
& H^1(L, \mu_2) \arrow[dd, "cor_{L/K}" near end] \arrow[rr,"e_L" near end]  && H^1(L,Spin(A,\sigma) )  \arrow[rr,"s_L", near start] && H^1(L,O^+(A,\sigma) )  \\ 
H^1(K, \mu_2)  \arrow[dr, "id"]  \arrow[rr, "f_K" near end] && H^1(K, Z)  \arrow[rr, "g_K"] \arrow[dr, "\alpha_K"]  &&
  H^1(K,\mu_2) \arrow[dr, "\beta_K"] \\
& H^1(K, \mu_2 ) \arrow[rr, "e_K"] && H^1(K, Spin(A,\sigma))   \arrow[rr,crossing over,"s_K" ] &&
  H^1(K,O^+(A,\sigma))
\end{tikzcd}\]
%\end{diag}

It is commutative and exact in all rows. Let $u \in \Ker (\alpha_L)$. We will associate an element $b\in K^*$ (not unique) to $u$, and will show that $cor_{L/K} (u) \in \ \Ker (\alpha_K)$ if and only if $[b]\in H^1(K,\mu_2)$ is a spinor norm.

Put $v=g_L(u)$, $w=cor_{L/K} (u)$. By commutativity of the above diagram, $v\in \ \Ker (\beta_L)$, so by Theorem \ref{5}, $cor_{L/K}(v) \in \ \Ker (\beta_K)$. We also have

$$s_K(\alpha_K(w))= \beta_K(g_K(w))=\beta_K(cor_{L/K}(v))=1,$$hence $\alpha_K (w) \in \Ker (s_K)= \Ima (e_K)$. Let ${bK^{*}}^2 \in H^1(K, \mu_2)$ be a preimage of $\alpha_K(w)$ under the map $e_K$. We have

$${bK^{*}}^2 \in \Ker (e_K) \iff w \in \Ker (\alpha_K)$$

\textbf{
Therefore, the \textit{obstruction to the norm principle}    for   $Spin(A,\sigma)$ over $L/K$ lives in $H^1(K,\mu_2)$. If one shows that this obstruction is a spinor norm, then the norm principle holds.}

\subsection{Reduction to quadratic field extensions} In this subsection, we reduce the norm principle for the group $Spin(A, \sigma)$ from arbitrary finite separable field extensions to quadratic field extensions.

\begin{thm}
\label{9}
Let $L/K$ be a finite separable field extension of an odd degree $n$, and $a\in K$. Then $a{K^{*}}^2 \in H^1(K,\mu_2)$ is a spinor norm iff  $a{L^{*}}^2 \in H^1(L,\mu_2)$ is a spinor norm.
\begin{proof}

\textbf{From K to L}: trivial.

\textbf{From L to K}:  

Since $n$ is odd we have 

$$cor_{L/K} (a{{L}^*}^2)= cor_{L/K} \circ  res_{L/K}(a{{K}^*}^2)=a^n {{K}^*}^2 = a {{K}^*}^2,$$ and by Theorem \ref{8}, 

$$cor_{L/K} (a{{L}^*}^2) \in \Ker (h_K: H^1(K,\mu_2) \longrightarrow H^1(K, Spin(A,\sigma))),$$i.e. $a{{K}^*}^2$ is a spinor norm.
\end{proof}
\end{thm}

\begin{thm}
\label{152}
Let $E/K$ be an infinite algebraic extension such that any finite subextension $L/K$ is of odd degree. Then for any $a \in K$ , $a {{K}^*}^2 \in H^1(K, \mu_2)$ is a spinor norm iff $a{{E}^*}^2 \in H^1(E, \mu_2) $ is a spinor norm.
\begin{proof}
\textbf{From K to E}: trivial.

\textbf{From E to K}: Consider the following exact sequence

$$O^+(A,\sigma)(E)  \overset{\delta_E}\longrightarrow H^1(E,\mu_2) \overset{h_E}\longrightarrow H^1(E,Spin (A,\sigma)).$$

Suppsoe $a{{E}^*}^2 \in H^1(E, \mu_2) $ is a spinor norm, i.e.
 $a {{E}^*}^2 \in \ \Ker h_E$. Then there is an element $M \in O^+(A,\sigma)(E)$ such that $\delta_E(M)=a {{E}^*}^2$.

Let $L$ be a subfield of $E$ which is a finite extension of $K$, such that $M \in O^+(A,\sigma)(L)$. Then the element $a {{L}^*}^2$   becomes a spinor norm over $L$, and since $L/K$ is of odd degree, we conclude the result by Theorem \ref{9}.
\end{proof}
\end{thm}

Now we prove the main result of this subsection:

\begin{thm}\label{11}

Let $(A,\sigma)$ be a central simple algebra with orthogonal involution defined over a field $K$. Assume that the norm principle holds for  $Spin(A,\sigma)$ over any separable quadratic field extension $E/F$, where $F$ contains $K$. Then the norm principle holds for  $Spin(A,\sigma)$ over any finite separable field extension $L/K$.
 \begin{proof}

Consider the following commutative diagram whose vertical arrows are inclusions:

\[\
\begin{tikzcd}
1 \arrow[r] & \mu_2 \arrow[r] \arrow[d] & Z \arrow[r] \arrow[d]      & \mu_2 \arrow[d] \arrow[r] & 1 \\
1 \arrow[r] & \mu_2 \arrow[r]           & {Spin(A,\sigma)} \arrow[r] & {O^+(A,\sigma)} \arrow[r] & 1
\end{tikzcd}
\]\

It gives rise to:

%\begin{diag}\label{10}
\[ \begin{tikzcd}[sep=0.7em]
 H^1(L, \mu_2) \arrow[rr, "f_L"] \arrow[dr,"id"]  \arrow[dd, "cor_{L/K}"]&& H^1(L, Z) \arrow[dd, "cor_{L/K}" near end] \arrow[rr,"g_L" near end]    \arrow[dr,"\alpha_L"] && H^1(L, \mu_2) \arrow[dr, "\beta_L"] \arrow[dd, "cor_{L/K}" near end] \\
& H^1(L, \mu_2) \arrow[dd, "cor_{L/K}" near end] \arrow[rr,"e_L" near end]  && H^1(L,Spin(A,\sigma) )  \arrow[rr,"s_L", near start] && H^1(L,O^+(A,\sigma) )  \\ 
H^1(K, \mu_2)  \arrow[dr, "id"]  \arrow[rr, "f_K" near end] && H^1(K, Z)  \arrow[rr, "g_K"] \arrow[dr, "\alpha_K"]  &&
  H^1(K,\mu_2) \arrow[dr, "\beta_K"] \\
& H^1(K, \mu_2 ) \arrow[rr, "e_K"] && H^1(K, Spin(A,\sigma))   \arrow[rr,crossing over,"s_K" ] &&
  H^1(K,O^+(A,\sigma))
\end{tikzcd}\]
%\end{diag}

It is commutative and exact in all rows. Let $u \in \Ker (\alpha_L)$. We want to show that $cor_{L/K} (u) \in \ \Ker (\alpha_K)$.

Let $b\in K^*$ be the element associated to $u$ which was introduced in Subsection \ref{45}. We need to show that 
$b$ is a spinor norm over $K$.

Let $\theta \in L$ be a generator for the extension $L/K$, i.e. $L=K(\theta)$. Assume that $f$ is the minimal polynomial of $\theta$ over $K$. Take any  algebraic extension $E/K$. We have
 
 $$L \underset{K}{\otimes} E= \frac{K[x]}{(f(x))} \underset{K}{\otimes} E = \frac{E[x]}{(f(x))}$$
 
 Let $f= f_1 f_2 \dots f_m$ be the factorization of $f$ in $E[x]$ into irreducible polynomials (note that $f_i$'s are distinct because of separability of $L/K$), and for every $1 \leq i \leq m$, let $L_i$ be $\frac{E[x]}{(f_i)}$. Then 
 
 $$L \underset{K}{\otimes} E= \frac{E[x]}{(f_1)} \times \frac{E[x]}{(f_2)} \times \dots \times  \frac{E[x]}{(f_m)}= L_1 \times L_2  \times \cdots \times L_m.$$
 
 So $L \underset{K}{\otimes} E= L_1 \times L_2  \times \cdots \times L_m$ is an \'{e}tale algebra where each $L_i$ is a finite separable extension of $E$ of degree equal to $deg f_i$.
 
 Let $\Gamma$ be the Galois group of $K^{sep}/K$ and $\Delta \subseteq \Gamma$ be a 2-Sylow subgroup. This corresponds to a tower $K \subseteq E \subseteq K^{sep}$ of fields where E is the elementwise fixed field of $\Delta$, so that $Gal(K^{sep}/E)=\Delta$. The degree of every finite subextension $L_i/E$ is going to be  a power of 2, and the degree of every  subextension of finite degree of $E/K$ is odd. By Theorem \ref{152}, in order to show that $b \frac{K^*}{{K^{*}}^2}$  is a spinor norm over $K$, it suffices to show it is a spinor norm over $E$.

 Consider the restriction map, which sends $u$   to $(u_1 , \dots,    u_m )$:
 
 $$res_{L \underset{K}{\otimes} E/L}: H^1(L, Z) \longrightarrow H^1(L \underset{K}{\otimes} E, Z) = H^1(L_1 \times L_2  \times \cdots \times L_m, Z)= H^1(L_1, Z) \times \dots \times  H^1(L_m, Z)$$

 $$   u \mapsto (u_1, \cdots, u_m) \  \  \  \ 
 \  \  \  \  \  \  \  \     \  \  \  \  \  \  \  \   \  \   \  \   \  \  \   \   \ \  \    \   \  \  \ \   \ \      \    \   \  \ \   \  \  \  \ \   $$
 
We have the following diagram (all arrows are inclusion maps)

\[\
\begin{tikzcd}
E \arrow[rr]            &  & L \underset{K}{\otimes} E \\
                        &  &                           \\
K \arrow[uu] \arrow[rr] &  & L \arrow[uu]             
\end{tikzcd}
\]\

 which gives rise to the following (the group $Z$ is abelian, so the corestriction maps exist):

\[\
\begin{tikzcd}
{H^1(E,Z)}                         &  & {H^1(L \underset{K}{\otimes} E,Z)= H^1(L_1, Z) \times \dots \times  H^1(L_m, Z)} \arrow[ll, "cor_{L \underset{K}{\otimes} E/E}"'] \\
                                   &  &                                                                                                                                   \\
{H^1(K,Z)} \arrow[uu, "res_{E/K}"] &  & {H^1(L,Z)} \arrow[uu, "res_{L \underset{K}{\otimes} E/L}"'] \arrow[ll, "cor_{L/K}"]                                              
\end{tikzcd}
\]\

By chasing $u$ in the diagram above, we get
 
\begin{equation}
    \label{12}
    res_{E/K} (cor_{L/K}(u))=  \prod_{i=1}^{m} cor_{{L_i} / E}(u_i).
\end{equation}

 Since $\Delta$ is a 2-Sylow subgroup, for every $1\leq i \leq m $ there exists a tower $E=L_{i,0}\subseteq L_{i,1} \subseteq L_{i,2} \subseteq \dots \subseteq L_{i,deg(f_i)}=L_i$ where each $L_{i,j}$ is a quadratic extension of $L_{i,j-1}$. We have
 
$$ cor_{L_i/E}(u_i)= cor_{L_{i,1}/L_{i,0}} \circ cor_{L_{i,2}/L_{i,1}} \circ \dots \circ cor_{L_{i, deg f_i}/L_{i,deg f_i -1}}(u_i).$$

By assumption, the norm principle holds for any separable quardratic field extension $F_1/ F_2$. Hence each $cor_{{L_i}/E} (u_i)$ becomes a spinor norm, and so does their product \eqref{12}.
 \end{proof}
\end{thm}

\section{Reduction to anisotropic skew-hermitian forms}\label{109}
In this section, we reduce the norm principle for spinor groups of skew-hermitian forms to the case where the skew-hermitian form is anisotropic.

Recall that the main result of this paper, namely Theorem \ref{40}, was proved in the special case where the underlying algebra 
$D$
 splits over 

$K$ in \cite{4}. By Morita theory, this split case corresponds to spinor groups of quadratic forms. A key step in the proof of the split case in \cite{4} is the reduction of the norm principle for spinor groups of quadratic forms to the case where the underlying quadratic form is anisotropic (see \cite[Corollary 2.6]{4}).

In the present section, we generalize this reduction to the setting of spinor groups of skew-hermitian forms, showing that it suffices to treat the anisotropic case. The main result of this section is stated in Remark \ref{4000}, and it will be used later in the proof of Theorem \ref{40}.

%For spinor group of quadratic forms, the reduction to anisotropic case can be done by directly showing the norm principle for the spinor groups of isotropic quadratic forms. For an isotropic quadratic form $q$ over $K$, we know that the spinor norm map of $q$ is surjective over $K$, hence every scalar is a spinor norm. In Subsection \ref{45}, we showed that the obstructions to norm principle for $Spin(q)$ over $K$ are elements in $H^1(K, \mu_2)$, and the norm principle holds for $Spin(q)$ over $K$ if the obstructions are spinor norms of $q$ over $K$. Therefore, the norm principle holds for $Spin(q)$ over $K$ if $q$ is $K-$isotropic. What we will show for skew-hermitian forms is slightly different: we will prove that if the norm principle holds for the spinor groups of all anisotropic skew-hermitian forms, then it does so for the spinor groups of all skew-hermitian forms.

We first set up the notations which will be used in this section. Let $K$ be a field of characteristic $\ne 2$, $(D,\tau)$ a central division algebra with symplectic involution defined over $K$, and $h$ a skew-hermitian form defined over $(D,\tau)$. 

Throughout the current section, $F$ denotes any field containing $K$. Let $h_F=h_{1,F} \perp h_{2,F}$ be the Witt decomposition of $h$ over $F$: $h_1$ and $h_2$ are skew-hermitian forms over $(D_F,\tau_F)$ which are totally hyperbolic and anisotropic over $F$, respectively. Note that $h_1$ or $h_2$ can be trivial, i.e., we could have $h=h_1$ or $h=h_2$. 

If $h$ is isotropic over $F$, i.e. $h_{1,F}$ is nontrivial, then we can define the algebraic groups $H$, $T$, and $\widetilde{T}$ as follows: the special orthogonal group $O^+(h_1)$ is an $F$-split group, and $O^+(h_2)$ is an $F$-anisotropic group.  Let $T$ be a maximal $F$-split torus in $O^+(h_1)$. 
 The direct product $O^+(h_1) \times O^+(h_2)$ is a subgroup of $O^+(h)$. Let $H$ be the preimage of  $O^+(h_1) \times O^+(h_2)$  under the map $Spin(h) \to O^+(h)$. The group $H$ is in fact an almost direct product of $Spin(h_1)$ and $Spin(h_2)$, and the direct product $Spin(h_1) \times Spin(h_2)$ is a simply connected cover of $H$. 
We also define $\widetilde{T}$ as the connected component of the preimage of $T$ in $H$ under the map $H \to O^+(h_1) \times O^+(h_2)$. The group $\widetilde{T}$ is a maximal $F-$split torus in $H$.

        We denote the center of $Spin(h)$ by $Z=C(Spin(h))$. 

\begin{lem}\label{257}
If $h$ is isotropic over $F$, then  the group  $Z$ can be embedded in the centralizer $C_H(\widetilde{T})$ (as algebraic groups over $F$).
\begin{proof}
Recall that the center of the group $O^+(h)$ is $\mu_2= \{\pm 1 \}$.
If $h$ is totally hyperbolic over $F$, then $H=Spin(h)$ and the result follows immediately. So assume that $h$ is not totally hyperbolic, hence the forms $h_1$ and $h_2$ are both nontrivial. The group $O^+(h_1) \times O^+(h_2)$ contains the element $(-1, -1)$, which is the non-identity element of the center of $O^+(h)$. Hence the preimage of  $O^+(h_1) \times O^+(h_2)$ contains the center of $Spin(h)$.
\end{proof}
\end{lem}
        
%By fixing a matrix representation of the group $O^+(h)(F)$, we can view the non-identity element of the center of $O^+(h)(F)$, which is $\mu_2(F)$, as the block matrix
 %        $
%-I=\begin{bmatrix}
%- I_1 & 0 \\
%0 & -I_2
%\end{bmatrix}$
 Consider the centralizer $C_H(\widetilde{T})$.  We clearly have the embeddings:
        $$Z \hookrightarrow  C_H(\widetilde{T})\hookrightarrow C_{Spin(h)}(\widetilde{T}).$$

      \begin{lem}\label{752}

          If $h$ is isotropic over $F$, then 
          $C_H(\widetilde{T})=C_{Spin(h)}(\widetilde{T})$ (as algebraic groups over $F$).
\begin{proof}
If $h$ is totally hyperbolic over $F$, then $H=Spin(h)$ and we are done. So assume that $h$ is not totally hyperbolic, i.e., $h_1$ and $h_2$ are both nontrivial.

    We prove that $C_{O^+(h)}(T)$ embeds in $C_{O^+(h_1) \times O^+(h_2)}(T)$ over $F$, and then the claim follows immediately, by considering the preimage of the centralizers under the map $Spin(h) \to O^+(h)$. Note that the restriction of the map $H \to  O^+(h_1) \times O^+(h_2)$ to $\widetilde{T}$ is surjective.

    Let us fix a matrix representation of the group $O^+(h)(F)$. View any element $g\in C_{O^+(h)}(T)(F)$ as a block matrix $
\begin{bmatrix}
A & B \\
C & D
\end{bmatrix}$, where the matrices $A$, $B$, $C$, and $D$ have coefficients in $F$ (the blocks $A$ and $D$ are square matrices). We need to show that $g\in O^+(h_1)(F) \times O^+(h_2)(F)$. By assumption, for every element $t\in T(F)$ we have $gt=tg$. The matrix representing $t$ has the form  $\begin{bmatrix}
M & 0 \\
0 & I
\end{bmatrix}$, where $M$ is a  diagonal matrix (because $T$ is a maximal $F-$split torus in $O^+(h_1)$), and $I$ is the identity matrix. A straightforward computation shows that $MB=B$ and $CM=C$. Therefore both matrices $B$ and $C$ must be zero, hence $g\in O^+(h_1)(F) \times O^+(h_2)(F)$. 
\end{proof}
      \end{lem}

The following theorem due to Merkurjev will be used in this section:

\begin{thm}\label{33}
  
$$PGO^+ (h)(F)/R \cong G_F(h)/ {F^*}^2 Hyp_F(h)$$

where:

$PGO^+ (h)(F)/R$ is the group of $R-$equivalence classes of $PGO^+ (h)$ over $F$,

 $G_F(h)=\{ s\in F^* \ | \ sh \cong h\}$ is the group of multipliers of $h$, and
 
$Hyp_F(h)$ is the subgroup of $G_F(h)$ generated by $\{N_{E/F} (E^*) \  |  \ h  \ \textit{is hyperbolic over} \  E\}$ where $E$ 
runs over finite extensions of $F$.
\begin{proof}
See \cite[Theorem~1]{16}.
    \end{proof}
\end{thm}

Before stating the next lemma, we recall that in this section, the algebra $D$ is defined over the field $K$, the skew-hermitian form $h$ is defined  over $(D_K,\tau_K)$, and $F$ is a field containing $K$ (see the notations set up before Lemma \ref{257}). The objects $H$, $T$, and $\widetilde{T}$ are defined over $F$, only when $h$ is $F$-isotropic.
In Lemmas \ref{257} and \ref{752}, the form $h$ is assumed to be $F-$isotropic.

\begin{lem}\label{215}
Assume that the extension $F/K$ is  finite and separable, and 
 $h$ is totally hyperbolic over $F$. Then the norm principle holds for $Spin(h)$ over $F/K$.
\begin{proof}
Since $h$ is hyperbolic over $F$, then $Hyp_F(h)=F^*$. This implies $G_F(h)=F^*$, and by Theorem \ref{33}, $PGO^+(h)(F)$ is $R-$trivial. 

Consider the diagram

\[\ 
\begin{tikzcd}
PGO^+(h)(F) \arrow [r, "\epsilon_F"] & H^1(F,Z) \arrow[r, "\delta_F"] \arrow[d, "cor_{F/K}"'] & H^1(F,\widetilde{T}) \arrow[r, "\gamma_F"] & H^1(F, C_{Spin(h)}(\widetilde{T})) \arrow[r, "\beta_F"] & H^1(F,Spin(h)) \\
PGO^+(h)(K) \arrow[r, "\epsilon_K"] & H^1(K,Z) \arrow[rrr, "\alpha_K"]                        &   &   & H^1(K,Spin(h))
\end{tikzcd}
\]\

Since the group $Spin(h)$ is split over $F$, we have $H=Spin(h)$, and $C_H(\widetilde{T})=\widetilde{T}$, so the map $\gamma_F$ is the identity map (see Lemma \ref{752}). The torus $\widetilde{T}$ is split over $F$, therefore by Hilbert 90  we have $H^1(F, \widetilde{T})=\{1\}$. So the map $\beta_F \circ \gamma_F \circ \delta_F$ is the zero map, hence the map $\epsilon_F$ is surjective. This fact, in conjunction with $R-$triviality of $PGO^+(h)(F)$ and Theorem \ref{7} implies that $cor_{F/K}(ker \ \beta_F \circ \gamma_F \circ \delta_F) \subseteq ker \alpha_K$, as desired. 
\end{proof}
    
\end{lem}

 A straightforward implication of Lemma \ref{215} is the following:
 
\begin{cor}\label{258}
   Assume that the extension $F/K$ is  finite and separable, and 
 $h$ is totally hyperbolic over $K$. Then the norm principle holds for $Spin(h)$ over $F/K$.
\end{cor}

\begin{thm}\label{2000}
Assume that $L/K$ is a finite separable field extension, and 
 $h$ is isotropic over $K$. Assume that for every skew-hermitian form $\mathfrak{h}$ defined over $(D_K,\tau_K)$ with dimension $< dim(h)$ the norm principle holds for $Spin(\mathfrak{h})$ over $L/K$. Then the norm principle holds for $Spin(h)$ over $L/K$.
    \begin{proof}
  We retain the notation established throughout this section, with $F=K$. So the forms $h_1$ and $h_2$ are defined over $K$ and they are totally hyperbolic and anisotropic over $K$, respectively.  We  shall make use of the groups $H$ and $\widetilde{T}$, which are defined over $K$.

  Consider the following embeddings over $K$:

  $$Z \hookrightarrow  C_H(\widetilde{T})= C_{Spin(h)}(\widetilde{T}) \hookrightarrow H \hookrightarrow Spin(h),$$

  which give rise to the following diagram

\[\ 
\begin{tikzcd}
H^1(L,Z) \arrow[d, "cor_{L/K}"] \arrow[r, "\delta_L"] & H^1(L, C_H(\widetilde{T})) \arrow[r, "\gamma_L"]  & H^1(L,H) \arrow[r, "\beta_L"]  & H^1(L,Spin(h))  \\
H^1(K,Z) \arrow[r, "\delta_K"]                 & H^1(K, C_H(\widetilde{T})) \arrow[r, "\gamma_K"]     & H^1(K,H) \arrow[r, "\beta_K"]                             & H^1(K,Spin(h))
\end{tikzcd}
\]\

The map $\beta_L \circ \gamma_L$ has trivial kernel \cite[Corollary~6.4]{235}. 
Let $u\in ker \ \beta_L \circ \gamma_L \circ \delta_L$. Then $u\in ker \ \delta_L$, and therefore $u\in ker \ \gamma_L \circ \delta_L$. So it is sufficient to prove the norm principle for the pair $(Z,H)$ over $L/K$. 

Assume that $h_2$ is trivial (equivalently, $h$ is totally hyperbolic over $K$), then we are done by Corollary \ref{258}.

Now assume that $h$ is not totally hyperbolic over $K$. Therefore, $dim(h_1)$ and $dim(h_2)$ are both strictly less than $dim(h)$. By our construction, the direct product $Spin(h_1) \times Spin(h_2)$ is a simply connected cover of $H$. So it is sufficient to prove the norm principle for $Spin(h_1) \times Spin(h_2)$ over $L/K$ (see Lemma \ref{31}). By assumption, the norm principle holds for each group  $Spin(h_1)$ and $Spin(h_2)$ over $L/K$. The norm principle for a direct product trivially reduces to the norm principle for each component. Hence the norm principle holds for $Spin(h_1) \times Spin(h_2)$ over $L/K$. 
 \end{proof}
\end{thm}

Using induction on $dim(h)$, we obtain the following conclusion:
\begin{cor}\label{365}
     Let $L/K$ be a finite separable field extension, $(D,\tau)$ a central $K$-division algebra with symplectic involution defined over $K$. Assume that for every $K-$anisotropic skew-hermitian form $\mathfrak{h}$ defined over $(D,\tau)$ the norm principle holds for $Spin(\mathfrak{h})$ over $L/K$. Then the norm principle holds for $Spin(h)$ over $L/K$, for  every skew-hermitian form $h$ defined over $(D,\tau)$.
\end{cor}

According to Theorem \ref{11}, in the proof of the main theorem (in Section \ref{184}) we will have the hypothesis that the degree of the field extension in the question is $2$. So in this section, from now on, we assume that $L/K$ is a  quadratic extension (hence separable, recall that $char \ K \ne 2$).
 
 Also, in the rest of this section, we assume that $D$ is a quaternion algebra defined over $K$ (which we denote by $D_K$ as well). We will further reduce the norm principle for $Spin(h)$ over $L/K$ to the case that $h$ is $L-$anisotropic as well. In order to prove this reduction, we need the following fact.

\begin{thm}\label{100}
    Let $(D,\tau)$ be a quaternion division algebra with canonical involution over $K$, and $L/K$ a quadratic separable field extension. Let $h$ be a skew-hermitian form over $(D,\tau)$. Assume that $D$ splits over $L$, and $h$ is $L-$ isotropic. Then there is a diagonalization  $h=<d_1, \dots, d_n>$ over $K$ such that $K(d_1)\cong L$. 
\begin{proof}
See \cite[Theorem~A.1]{18}.
\end{proof}
\end{thm}

\begin{lem}\label{366}

Let $L/K$ be a quadratic field extension,  $(D,\tau)$ be a  quaternion division algebra with canonical involution over $K$, and $h$ be a skew-hermitian form over  $(D_K,\tau_K)$ which is $L-$isotropic. Assume that for every skew-hermitian form $\mathfrak{h}$ over  $(D,\tau)$ with dimension less than dimension of $h$, the norm principle holds for $Spin(\mathfrak{h})$ over $L/K$. Then the norm principle holds for $Spin(h)$ over $L/K$.
    \begin{proof}

        If $h$ is $K-$isotropic, then we are done by Theorem \ref{2000}. So assume that $h$ is $K-$anisotropic. 

We introduce skew-hermitian forms $h_1$ and $h_2$ over $(D_K,\tau_K)$, such that $h_K=h_{1,K} \perp h_{2,K}$, $h_{1,L}$ is totally hyperbolic, and $h_{2,L}$ is anisotropic.

\begin{itemize}
   
        \item Case 1. Assume that $D$ splits over $L$. Then we take $h_1$ to be the form $<d_1>$ introduced in Theorem \ref{100}, and $h_2=:<d_2, \dots, d_n>$. Since $L=K(d_1)$, then $O^+(h_1)$ is the norm one torus $R_{L/K}^{1}(\mathbb{G}_m)$ defined over $K$. Therefore, $O^+(h_1)$ is split over $L$, and $h_{1,L}$ is isotropic. Being one-dimensional, it follows that in fact $h_{1,L}$ is hyperbolic.  If $h_{2,L}$ is anisotropic, then we are done. Otherwise, we apply Theorem~\ref{100} to $h_2$
and repeat the above procedure. Iterating this process finitely many times, we eventually obtain
a decomposition $h_K \simeq h_{1,K} \perp h_{2,K}$
such that $h_{1,L}$ is totally hyperbolic and $h_{2,L}$ is anisotropic.

      \item 
   Case 2: Assume that $D$ does not split over $L$.

   let $\theta\in L^*$ be a generator of the  extension $L/K$ such that $\theta^2 \in K$ (such element exists because $char \ K \ne 2$). Let $V=V_K$ be the left vector space over $D_K$ on which $h$ is defined. Take an isotropic vector  $v \in V_L$. Write $v=v_1 + \theta v_2$ where $v_1 , v_2 \in V$. Let $W:=<v_1,v_2>$ be the subspace in $V$ generated by $v_1$, $v_2$ and 
    $h_1$ be the restriction of $h$ to $W$. The subspace $W$ has to be two dimensional, because otherwise there is a scalar $\lambda \in D$ such that $v_1= \lambda v_2$ and we get $v=v_1 + \theta v_2 = (\lambda + \theta) v_2$, hence $h(v_2)=0$, and $h$ becomes isotropic over $K$ as well. The skew-hermitian form $h_{1,L}$ is hyperbolic, because it is two-dimensional and isotropic. By taking an orthogonal complement of $W$ in $V$, we can write $h_K=h_{1,K} \perp h_{2,K}$. If $h_{2,L}$ is anisotropic, we are done. Otherwise, similar to case 1, we repeat this process until we get the desired decomposition $h_K=h_{1,K}\perp h_{2,K}$.
\end{itemize}

    Note that in both cases above, the form $h_2$ may be trivial. This would happen if $h$ is totally hyperbolic over $L$. In that case we have $h=h_1$. 

 In this proof, the decomposition $h_K=h_{1,K} \perp h_{2,K}$ is obtained over $K$, but after extending the scalars to $L$, it coincides with the Witt decomposition $h_L=h_{1,L} \perp h_{2,L}$ (over $L$) introduced earlier in this section (taking $F=L$), hence we retain the notations and the objects defined at the beginning of this section: the groups $H$ and $\widetilde{T}$ are now defined over $L$. By Lemmas \ref{257} and \ref{752}, we have the following $L-$embeddings: 

  $$Z \hookrightarrow  C_H(\widetilde{T})= C_{Spin(h)}(\widetilde{T}) \hookrightarrow H \hookrightarrow Spin(h).$$

Note that even though the forms $h_1$ and $h_2$ are defined over $K$, but the groups $H$ and $\widetilde{T}$ are only defined over $L$.
 
Consider the diagram

\[\ 
\begin{tikzcd}
H^1(L,Z) \arrow[r, "\delta_L"] \arrow[d, "cor_{L/K}"'] & H^1(L, C_H(\widetilde{T})) \arrow[r, "\gamma_L"] & H^1(L,H) \arrow[r, "\beta_L"] & H^1(L,Spin(h)) \\
H^1(K,Z) \arrow[rr, "\alpha_K"]                        &                         & H^1(K,H)  \arrow[r, "\beta_K"] & H^1(K,Spin(h))
\end{tikzcd}
\]\

The map $\beta_L \circ \gamma_L$ has trivial kernel \cite[Corollary~6.4]{235}. 
Let $u\in ker \ \beta_L \circ \gamma_L \circ \delta_L$. Then $u\in ker \ \delta_L$, and therefore $u\in ker \ \gamma_L \circ \delta_L$. So it is sufficient to prove the norm principle for the pair $(Z,H)$ over $L/K$. 

Assume that $h_2$ is trivial (equivalently, $h$ is totally hyperbolic over $L$). Then by Lemma \ref{215} we are done.

Now assume that $h$ is not totally hyperbolic over $K$. Therefore, $dim(h_1)$ and $dim(h_2)$ are both nonzero and strictly less than $dim(h)$. By our construction, the direct product $Spin(h_1) \times Spin(h_2)$ is a simply connected cover of $H$. So it is sufficient to prove the norm principle for $Spin(h_1) \times Spin(h_2)$ over $L/K$ (see Lemma \ref{31}). By assumption, the norm principle holds for each group  $Spin(h_1)$ and $Spin(h_2)$ over $L/K$. The norm principle for a direct product trivially reduces to the norm principle for each component. Hence the norm principle holds for $Spin(h_1) \times Spin(h_2)$ over $L/K$. 
\end{proof}
\end{lem}

\begin{rema}\label{4000}
    According to Theorem \ref{11}, Corollary \ref{365}, and using induction on the result of Lemma \ref{366}, we may assume that the extension $L/K$ on which we want to show the norm principle for $Spin(h)$ is quadratic, and $h$ is $L$-anisotropic.
\end{rema}

\section{Larmour's theorem}\label{110}\

 Let $K$  be a complete discretely valued field with uniformizer $\pi$ and residue field $k$ with $char \ k\ne 2$. Let $\mathcal{O}$ be the ring of integers of $K$ and $m$ be its maximal ideal.

 If a quadratic form over $K$ has a diagonalization with unit elements in $\mathcal{O}$, then we call it a unit quadratic form (or unramified quadratic form). To any unit quadratic form $q=<u_1, u_2, \dots, u_n>$ over $K$, one can associate the quadratic form $<\bar{u_1}, \bar{u_2}, \dots, \bar{u_n}>$ over $k$, which we denote by $\Bar{q}$.

 Any quadratic form $q$ over $K$ can be diagonalized by elements in $\mathcal{O}$ and each diagonal entry can assumed to be of the form $u$ or $u \pi$, where $u$ is a unit (because any element of $K$ is in one of those forms up to squares). Therefore $q$ can be written as $q_0 \perp \pi q_1$, where $q_0$ and $q_1$ are unit quadratic forms. The quadratic forms $q_0$ and $q_1$ are determined uniquely up to isometry if $q$ is $K$- anisotropic (see \cite[Chapter~6, Section~1]{14}). The quadratic forms $\Bar{q_0}$ and $\Bar{q_1}$ are called the first and second residue forms of $q$, respectively. In what follows, $W(\star)$ denotes the Witt ring of a field $\star$.

\begin{thm}\label{53}
    \textbf{(Springer's Theorem)} Let $K$  be a complete discretely valued field with uniformizer $\pi$ and residue field $k$ with $char \ k\ne 2$. Then we  have a group isomorphism

$$\partial= (\partial_0, \partial_1): W(K) \longrightarrow  W(k) \oplus  W(k)$$
$$\   \  \ \ \   \  \    \   \    \   \   \  \  \   \     \  [q]   \     \ \mapsto  \  \   \   {([\bar{q_0}], [\bar{q_1}]})$$where the quadratic forms $\bar{q_0}$ and $\bar{q_1}$ are the first and second residue forms of $q$, respectively (defined above).  
  
\end{thm}

\begin{rema}
   In Springer's Theorem, $\partial_1$ depends on the choice of $\pi$, but $\partial_0$ does not.
\end{rema}

 Let us denote the valuation on $K$ by $\nu_K$:

$$\nu_K: K \backslash \{0\} \longrightarrow \mathbb{Z}.$$

Let $D$ be a quaternion division algebra over $K$ with canonical (symplectic) involution $\tau$. By \cite{23}, the valuation on $K$ extends uniquely to the following valuation on $D$: 

$$\nu_{D}: D \backslash \{0\}   \longrightarrow \mathbb{Q}$$

$$ \   \  \   \   \  \   \   \   \   \  \  \   \  \  \   \   \  \   \   \   \   \  \  \   \  \  \  \  \  \     x \mapsto \frac{1}{2} \  \nu_K (Nrd (x)).$$

 Let $\mathcal{O}_D:= \{d \in D | \nu_{D}(d) \geq 0\}$ be the ring of integers of $D$, $m_D:=\{d \in D | \nu_{D}(d) > 0\}$ its maximal ideal, $\bar{D}:={\mathcal{O}_D}/{m_D}$ the residue division ring, and $\Gamma_D$ the image of $\nu_D$. We will denote by $\bar{*}$ the image of $* \in \mathcal{O}_D$ in $\bar{D}$.

 The inclusion $\mathcal{O}\hookrightarrow \mathcal{O}_D$  induces an embedding $k={\mathcal{O}}/m \hookrightarrow \bar{D}= {\mathcal{O}_D}/{m_D}$, and elements in the image of $k$ under this embedding commute with all elements of $\bar{D}$, because $K$ is the center of $D$. Therefore, $\bar{D}$ is a $k-$algebra. 
 
 Since the reduced norm map $Nrd: D \longrightarrow K$ is invariant under the involution $\tau$ (see \cite[Page~14, Corollary~2.2]{13}), the map

$$\bar{\tau}: \bar{D} \longrightarrow \bar{D}$$
$$\  \   \  \   \    \ x + m_D \mapsto \tau(x) + m_D$$
is an involution induced on $\bar{D}$.

Let $x, y$ be two elements in $D$ such that  $x^2=a$, $y^2=b$, $D=(\frac{a,b}{K})$. Also put $z:=xy=-yx$. The set $\{1,x,y,z\}$ is a $K-$basis of $D$.

We have two cases of ramification for $D$: it is either unramified (i.e. $a$ and $b$ are units), or ramified (i.e. $a$ is a unit and $b$ is a non-unit). The case where $D$ is generated by two non-units $a$ and $b$ can be reduced to the ramified case, since $(\frac{a,b}{K})\cong (\frac{a, -ab}{K})$ and if $a$ and $b$ are non-units, then $ab$ will be unit up to squares; the norm forms associated to $(\frac{a,b}{K})$ and $\cong (\frac{a, -ab}{K})$ are both isometric to the form $<1,-a,-b,ab>$, and the isometry of the norm forms implies the isomorphism of the quaternion algebras.

Any skew-hermitian form over $(D, \tau)$ has a diagonalization by skew-symmetric elements. 
 If a skew-hermitian form over $(D, \tau)$ has a diagonalization with unit elements, then we call it a \underline{unit skew-hermitian form} or \underline{unramified skew-hermitian form}, in a similar manner to unit quadratic forms. To any unit skew-hermitian form $h=<u_1, u_2, \dots, u_n>$ over $(D, \tau)$, one can associate the skew-hermitian form $<\bar{u_1}, \bar{u_2}, \dots, \bar{u_n}>$ over $(\Bar{D}, \Bar{\tau})$, which we denote by $\Bar{h}$ (this argument works for hermitian forms as well).

Let $h$ be a skew-hermitian form over $(D, \tau)$. Then $h$ can be decomposed as

$$h = h_0 \perp h_1,$$
where $h_0$ is a unit skew-hermitian form, and any element $s$ in a diagonalization of $h_1$ satisfies $\nu_D(s)= \min \{\nu_D(t) \ | \ t \in D\} \cap \mathbb{Q}^+$. This value depends on the ramification of the quaternion algebra. Larmour's theorem (\cite[Theorem~3.6]{15}) states that if $h$ is $K-$anisotropic, then this decomposition is unique up to isometry classes of $h_0$ and $h_1$.

Now we discuss the two cases of ramification of $D$:

\textbf{Case 1}: $D$ is unramified: $D=(\frac{a,b}{K})$ where both $a$ and  $b$ are units.

The assumption implies that $\pi$ is still a uniformizer for $\nu_D$ (unless otherwise stated, when we refer to the uniformizer of $\nu_D$, we mean $\pi$). We have $\Gamma_D=\mathbb{Z}$ and $\nu_D(x)=\nu_D(y)=\nu_D(z)=0$. Then $h$ can obviously be written as  $h= h_0 \perp h_1$, where $h_0$ is a unit skew-hermitian form and in the diagonalization of $h_1$ each element has value $1$, because multiplying any element in the diagonalization of $h$ by $\pi^2$ does not change the isometry class of $h$.

By \cite[Page~21]{22}, $\bar{D}$ is isomorphic to $(\frac{\bar{a}, \bar{b}}{k})$, which is a division quaternion algebra (the norm form of $D$ is anisotropic and by Hensel's lemma, the norm form of $\bar{D}$ is anisotropic too). Then the involution $\bar{\tau}$ is the following:

\begin{align*}
  \bar{\tau} : (\frac{\bar{a}, \bar{b}}{k}) & \rightarrow (\frac{\bar{a}, \bar{b}}{k}) \\
         p+q\bar{x} + r \bar{y} + s \bar{z} & \mapsto p-q\bar{x} - r \bar{y} - s \bar{z},
\end{align*}
for $p,q,r,s \in k$.\\

\textbf{Case 2}: $D$ is ramified: $D=( \frac{a, \pi}{K})$ where $a$ is a unit. In fact if $D=(\frac{a,b}{K})$ with $\nu_K(a)=0$ and $\nu_K(b)=1$, then we can choose $b$ to be the uniformizer of $K$, and hence without loss of generality assume that $b=\pi$.

In this case $y=\sqrt{\pi}$  is a uniformizer for $\nu_D$ (similar to case 1, unless otherwise stated, when we refer to the uniformizer of $\nu_D$, we mean $y=\sqrt{\pi}$). We have $\Gamma_D= \{\frac{l}{2}| l \in \mathbb{Z}\}$, $\nu_D(x)=0$, $\nu_D(y)=\nu_D(z)=\frac{1}{2}$.

In case 2, by \cite[Page~22]{22}), $\bar{D}$ is isomorphic to $k(\bar{x})$, which is a quadratic field extension of $k$. Note that in general, $k(\bar{x})$ is a quadratic \'{e}tale extension of $k$, but in our case we know that $\bar{x}$ is not a square in $k$ since otherwise by Hensel's lemma $x$ will be a square in $K$ which contradicts the fact that $D$ is a division quaternion algebra. Therefore, $k(\bar{x})$ is a quadratic field extension of $k$. 

The involution $\bar{\tau}$ will just be the nontrivial $k-$automorphism of $k(\bar{x})$, which is an involution of the second kind:

\begin{align*}
  \bar{\tau} :  k(\bar{x}) & \rightarrow k(\bar{x})  \\
         p+q\bar{x} & \mapsto p-q\bar{x},
\end{align*} for $p,q \in k$.

By \cite[Lemma~2.1]{27}, $h$ can be written as  $h= h_0 \perp h_1$, where $h_0$ is a unit skew-hermitian form and in the digonalization of $h_1$ each element has value $\frac{1}{2}$ (up to squares). Indeed, replacing any skew-symmetric element $\alpha$ in the diagonalization of $h$ by $- \pi^{-2} y\alpha y$ does not change the isometry class of $h$.

So in both cases, we get a similar decomposition as Springer's decomposition for quadratic forms: 

\begin{equation}\label{52}
    h= h_0 \perp h_1.
\end{equation}

Whenever we write $h= h_0 \perp h_1$, this is what we mean, and we will call this the \textbf{Larmour  decomposition of $h$}. So the elements in the diagonalization of $h_0$ are units, and in $h_1$ they are non-units. This notation will be used in both cases where $D$ is unramified or ramified.

We seek a more concrete understanding of how the diagonal entries in $h_0$ and $h_1$ look like.

\begin{lem}
    \label{63}

Let $h= h_0 \perp h_1$ be a skew-hermitian form over $(D, \tau)$ as before, over a ramified $D$ (case 2).

(1) The skew-hermitian form $h_0$ is isometric over $D$ to a unit skew-hermitian form whose diagonal entries are of the form $\alpha x$ where $\alpha  \in {\mathcal{O}_K}^{*}$. 

(2) The skew-hermitian form $h_1$ is isometric over $D$ to a skew-hermitian form whose diagonal entries are of the form $\beta  y + \gamma z$ where $\beta , \gamma  \in \mathcal{O}_K$, such that at least one of $\beta $ and $\gamma$ is unit and  $\nu_D(\beta  y + \gamma z)=\frac{1}{2}$. 

    \begin{proof}

     \cite[Remark~4.4 and Remark~4.6]{27}.
\end{proof}
\end{lem}

\begin{rema}\label{66}
    Let $h$ be a skew-hermitian form over  $(D, \tau)$, where $D$ is ramified, and $h=h_0 \perp h_1$ be Larmour decomposition of $h$. The proof of Lemma \ref{63} shows that if $u= \alpha x + \beta y + \gamma z$ is an element in the diagonalization of $h_0$, then replacing $u$ with $\alpha x$ will not change the isometry class of $h_0$ (and therefore $h$). Similarly,  if $u= \alpha x + \beta y + \gamma z$ is an element in the diagonalization of $h_1$, then replacing $u$ with $\beta y + \gamma z$ will not change the isometry class of $h_1$ (and therefore $h$). 
\end{rema}

The next theorem states that  Larmour decomposition (\ref{52}) is unique up to isometry for anisotropic skew-hermitian forms.

\begin{thm}\label{54}
    \textbf{(Larmour's theorem)} The skew-hermitian form $h=h_0 \perp h_1$ is $K-$anisotropic, if and only if $h_0$ and $h_1$ are $K-$ anisotropic. Furthermore, if $h$ is $K-$anisotropic, then the decomposition $h=h_0 \perp h_1$ is unique up to isometry classes of $h_0$ and $h_1$.
    \begin{proof}
See \cite[Theorem~3.3 and Theorem~3.6]{15}. Note that this fact is true regardless of the ramification of $D$.
    \end{proof}
\end{thm}

\section{Notations and some explicit formulas}\label{111}

\subsection{Notations}\label{112} Now we set up the notations which will be used in the proof of the main theorem, stated below:

\begin{thm}
    \label{40}
     Let $K$ be a complete discretely valued field with residue field $k$ with $char(k)\ne 2$.
    Assume that the norm principle holds for spinor groups $Spin(\mathfrak{h})$ for every regular skew-hermitian form $\mathfrak{h}$ over every quaternion algebra $\mathfrak{D}$ (with respect to the canonical involution on $\mathfrak{D}$) defined over any finite extension of $k$. Then, the norm principle holds for spinor groups $Spin(h)$ for every regular skew-hermitian form $h$ over every quaternion algebra $D$ (with respect to the canonical involution on $D$) defined over $K$.
\end{thm}

Now, we set up the notations which will be used in the proof of Theorem \ref{40}:

\begin{itemize}
    \item $K$: a complete discretely valued field with $char \ K \ne 2$.
    \item $\nu$: the valuation on $K$. 
          \item $\pi$: a uniformizer for $K$.
    \item $\mathcal{O}$: the ring of integers of $K$.
    \item $m$: the maximal ideal of $\mathcal{O}$.
    \item $k$: the residue field of $K$, with $char \ k \ne 2$.
    \item Whenever we say that an object is \underline{ramified}, we mean it is \underline{totally ramified} (for field extensions, quaternion algebras, elements in the fields and algebras, and skew-hermitian forms).
    \item $L/K$: a quadratic field extension. The reduction of the main theorem to quadratic field extensions was proved in Section \ref{135}. 
    \item $\sqrt{t}$: a generator of $L$ over $K$; since $char \ K \ne 2$, we may assume that $L=K(\sqrt{t})$, where $t \in K^* \backslash {K^*}^2$. If $L/K$ is an unramified extension, then $t$ is a unit of $K$. If $L/K$ is ramified, then  $t=c\pi$ for a unit $c\in \mathcal{O}^*$.
    \item $\overline{*}$ for $* \in \mathcal{O}$: the image of $*$ in $k$.  
    \item $l$: residue field of $L$. If $L/K$ is unramified, then $l=k(\sqrt{\overline{t}})$. If $L/K$ is ramified then $l=k$.

\item $\pi'$: a uniformizer of $L$. If $L/K$ is unramified, then $\pi'=\pi$. If $L/K$ is ramified then  $\pi'=\sqrt{c\pi}$ for a unit $c\in \mathcal{O}^*$ (recall that we assumed that $t=\sqrt{c\pi}$ is the generator of $L/K$).
\item $D$: division quaternion algebra over $K$. So $D=(\frac{a,b}{K})$ where $a$ and $b$ are units in $K$, or  $D=(\frac{a,\pi}{K})$ where $a$ is a unit in $K$. We call $D$ unramified in the first case, and ramified in the second one. Note that a ramified quaternion division algebra is of the form $D=(\frac{a,a' \pi}{K})$ where $a, a'$ are units in $K$. Without loss of generality, by replacing $a'\pi$ with $\pi$ if necessary, we may assume that $a'=1$. 
\item $\nu_L$, $\nu_D$, $\nu_{D_L}$ respectively denote the extensions of $\nu$ to $L$, $D$, and $D_L$. The extension to $D_L$ exists only if $D_L$ is an $L$-division algebra.
\item $1, x, y, z$: the basis for $D$ over $K$ such that if $D$ is unramified then $x^2=a$, $y^2=b$, $z=xy=-yx$, and if $D$ is ramified, then $x^2=a$, $y^2=\pi$, $z=xy=-yx$.

\item $\mathcal{O}_L$, $\mathcal{O}_D$, $\mathcal{O}_{D_L}$ respectively denote the ring of integers of $L$, $D$, and $D_L$. Note that $\mathcal{O}_{D_L}$ is only defined if $D_L$ is an $L$-division algebra. 
\item $m_L$, $m_D$, $m_{D_L}$ respectively denote the maximal ideals of $L$, $D$, and $D_L$. Note that $m_{D_L}$ is only defined if $D_L$ is an $L$-division algebra. 
\item $\overline{D}:=D/{m_D}$: the residue $k$-algebra. 
\item $\tau:$ the canonical involution on $D$.
\item $\overline{\tau}$: the induced canonical involution on $\overline{D}$.
\item $\overline{D}_l:= {D_L}/{m_{D_L}}$: the residue $l$-algebra. Note that  
$\overline{D}_l$ is only defined if $D_L$ is an $L$-division algebra. If $D_L$ is unramified over $L$, then $\overline{D}_l=(\frac{\overline{a}, \overline{b}}{l})$ will be a division quaternion $l$-algebra. 
If $D_L$ is ramified over $L$, then $\overline{D}_l=l(\sqrt{\bar{a}})$ is a quadratic field extension of $l$ (see Section \ref{110}).

 \item $\overline{*}$ for $* \in \mathcal{O}_L$, $\mathcal{O}_D$, or  $\mathcal{O}_{D_L}$: the image of $*$ in $l$, $\overline{D}$, or $\overline{D}_l$, respectively. Note that the image in $\overline{D}_l$ is only defined if $D_L$ is an $L$-division algebra. 

\item The following diagram contains all the fields, maximal ideals, rings of integers, and the algebras. All the maps are inclusions:

\begin{tikzcd}
                                     & m_{D_L} \arrow[rrrr]            &  &  &                                      & \mathcal{O}_{D_L} \arrow[rrrr]            &  &  &                         & D_L            \\
m_D \arrow[rrrr] \arrow[ru]            &                           &  &  & \mathcal{O}_D \arrow[rrrr] \arrow[ru]            &                           &  &  & D \arrow[ru]            &              \\
                                     & m_L \arrow[rrrr] \arrow[uu] &  &  &                                      & \mathcal{O}_L \arrow[rrrr] \arrow[uu] &  &  &                         & L \arrow[uu] \\
m \arrow[rrrr] \arrow[ru] \arrow[uu] &                           &  &  & \mathcal{O} \arrow[rrrr] \arrow[ru] \arrow[uu] &                           &  &  & K \arrow[ru] \arrow[uu] &             
\end{tikzcd}\\ \\ \\

\item $h$: a skew-hermitian form over $(D,\tau)$, where $\tau$ is the canonical involution on $D$. We assume that $h$ is anisotropic over $L$ (hence anisotropic over $K$ too). The reduction of the main question to anisotropic skew-hermitian forms was proved in Section \ref{109}.
\item $n$: dimension of  $h$ over $D$.
\item $\overline{h}$: the residue form of an unramified form $h$ defined over the residue algebra with canonical involution $(\overline{D}_l , \overline{\tau})$ (see Section \ref{110}).
\item $\sigma_h$: the involution adjoint to $h$ on the matrix algebra $M_n (D)$. 
\item $O^+(h):= O^+(M_n(D),\sigma_h)$, the special orthogonal group of $h$.
\item $Spin(h):= Spin(M_n(D),\sigma_h)$, the Spinor group of $h$.
\item $PGO^+(h):= PGO^+(M_n(D),\sigma_h)$, the projective   group of proper similitudes of $h$.
\item $C(G)$: center of any group $G$.
\item $Z=C(Spin(h))$.

\item The following diagram relates the groups of type $D_n$ which will be used in our argument. The rows are exact.

\[\
\begin{tikzcd}
1 \arrow[r] & \mu_2 \arrow[r] \arrow[d, "id"] & Z \arrow[r] \arrow[d]      & \mu_2 \arrow[d] \arrow[r] & 1 \\
1 \arrow[r] & \mu_2 \arrow[r]           & {Spin(h)} \arrow[r] & {O^+(h)} \arrow[r] & 1
\end{tikzcd}
\]\

\item The above diagram gives rise to the following commutative diagram:

\begin{diag}\label{139}
\[ \begin{tikzcd}[sep=0.7em]
 H^1(L, \mu_2) \arrow[rr, "f_L"] \arrow[dr,"id"]  \arrow[dd, "cor_{L/K}"]&& H^1(L, Z) \arrow[dd, "cor_{L/K}" near end] \arrow[rr,"g_L" near end]    \arrow[dr,"\alpha_L"] && H^1(L, \mu_2) \arrow[dr, "\beta_L"] \arrow[dd, "cor_{L/K}" near end] \\
& H^1(L, \mu_2) \arrow[dd, "cor_{L/K}" near end] \arrow[rr,"e_L" near end]  && H^1(L,Spin(h) )  \arrow[rr,"s_L", near start] && H^1(L,O^+(h) )  \\ 
H^1(K, \mu_2)  \arrow[dr, "id"]  \arrow[rr, "f_K" near end] && H^1(K, Z)  \arrow[rr, "g_K"] \arrow[dr, "\alpha_K"]  &&
  H^1(K,\mu_2) \arrow[dr, "\beta_K"] \\
& H^1(K, \mu_2 ) \arrow[rr, "e_K"] && H^1(K, Spin(h))   \arrow[rr,crossing over,"s_K" ] &&
  H^1(K,O^+(h))  
\end{tikzcd}\]
\end{diag}

\item $f$, $g$, $e$, $s$, $\alpha$, $\beta$ are the Galois cohomology maps in Diagram \ref{139}.

\item $u$: an element in $\Ker (\alpha_L)$ (see Diagram \ref{139}).

\item $\lambda\in L^*$ is a representative for $g_L(u)$ (see Diagram \ref{139}). Without loss of generality, we will assume that $\nu_L(\lambda)\in \{0, \nu_L(\pi')\}$. Note that if $\nu_L(\lambda)\ne 0$, then one of the following holds: $\nu_L(\lambda)= 1$ if $L/K$ is unramified, or $\nu_L(\lambda)= \frac{1}{2}$ if $L/K$ is ramified.
\item $h=h_0 \perp h_1$: Larmour decomposition of $h$ over $K$. Note that one of $h_0$ or $h_1$ could be trivial, i.e. with dimension $0$; in that case we will have $h=h_0$ (unramified) or $h=h_1$ (ramified). The base change of the skew-hermitian forms $h_0$ and $h_1$ from $K$ to $L$ are denoted by $h_{0,L}$ and $h_{1,L}$, respectively. 

\item $F$: the discriminant extension of $h$ over $K$, i.e. 
$$F:= K[t]/(t^2 - (-1)^{n}Nrd_{D/K} \ (det \ (h))).$$ 
This is a quadratic \'{e}tale extension of $K$. So $F$ is either a quadratic separable field extension of $K$, or it is isomorphic to $K \times K$. 

\item $\psi$: the nontrivial $K-$ automorphism of $F$.

\item $M:= L \underset{K}{\otimes} F$ the discriminant extension of $h$ over $L$. Therefore, $M$ is an  \'{e}tale quadratic extension of $L$: it is either a quadratic field extension of $L$, or isomorphic to $L \times L$.  
\item $\psi_L$: the nontrivial $L-$ automorphism of $M$.
\item The following diagram shows the discriminant extensions of $h$ over $K$ and $L$ (all arrows are inclusions):

\[\begin{tikzcd}
L \arrow[r] 
& M:= L \underset{K}{\otimes} F \\
K \arrow[r] \arrow[u] & F \arrow[u]
\end{tikzcd}
\]

\item $\nu_F$: the extension of the valuation $\nu$ from $K$ to $F$, in case $F$  is a field.

\item $\nu_M$: the extension of the valuation $\nu_L$ from $L$ to $M$, in case $M$ is a field.

Note: If $h$ is an unramified skew-hermitian form over $K$ and $M$ is a field, then the extension $M/L$ is unramified. Therefore, $\pi '$ is a uniformizer for $M$. 
\end{itemize}

\begin{lem}\label{94}
Recall that $\sqrt{t}$ is a generator of the field extension $L/K$. Assume that $D_L$ splits. Then 
$D_K\cong (\frac{t,w}{K})$ for some $w \in K$.
\begin{proof}
See \cite[Page~67, Theorem~4.1]{14}.
\end{proof}
\end{lem}

\textbf{Recall that to prove Theorem \ref{40}, we need to show that $cor_{L/K}(u)\in Ker   (\alpha_K)$  (see Diagram \ref{139})}.\\

In Sections \ref{111}, \ref{156}, \ref{157}, and \ref{128}, we set up the necessary steps towards the proof of Theorem \ref{40}. 

\subsection{Explicit formulas for the maps $f$ and $g$}\label{154} We recall explicit formulas for the maps $f$ and $g$ in Diagram \ref{139}.

Recall that when $n$ is even, we have $Z=R_{F/K}(\mu_2)$, and $H^1(K, Z) = {F^*}/{{F^*}^2}$ (see \cite[Page~332]{19}). Also in Diagram \ref{139}, the map $f_K :  {K^*}/{{K^*}^2}   \longrightarrow {F^*}/{{F^*}^2}   $        is the natural map
and $g_K$ is given by the following (see \cite[Proposition~13.33]{13}): $$g_K :  {F^*}/{{F^*}^2}       \longrightarrow  {K^*}/{{K^*}^2}$$

$$  \   \   \   \  \   \   \  \  \  \   \   \ \   [p] \mapsto [N_{F/K}(p)]$$for every $p \in F^*$. 

When $n$ is odd, then $Z=\Ker (Norm: R_{F/K}(\mu_4) \longrightarrow \mu_4)$ (see \cite[Page~332]{19}), and $H^1(K,Z)={U(K)}/{U_0(K)}$, where $U \subseteq \mathbb{G}_m \times  R_{F/K} \mathbb{G}_m$ is the subgroup defined by $$U(K):=\{(q,p) \in K^* \times F^* \ | \ q^4=Norm_{F/K}(p)\}$$ and $U_0 \subseteq U$ is the subgroup defined by $$U_0(K):=  \{(N_{F/K}(p),p^4) \ | \ p \in F^*\}.$$

Then the maps $f$ and $g$ in Diagram \ref{139} will be the following (by \cite[Proposition~13.36]{13}):

$$f_K: {K^*}/{{K^*}^2}   \longrightarrow {U(K)}/{U_0(K)} $$

$$[q] \mapsto [q,q^2]$$and 

$$g_K:  {U(K)}/{U_0(K)}  \longrightarrow   {K^*}/{{K^*}^2}    $$

$$[q,p] \mapsto [N_{F/K}(p_0)],$$where $p_0 \in F^*$ is such that $p_0 \psi {(p_0)}^{-1}=q^{-2}p$ (recall that $\psi$ is the nontrivial $K-$automorphism of $F$).\\

\subsection{The discriminant extension}\label{155} We will need the following facts about the discriminant extension of $h$ in the proof of the main theorem. Recall that $F$ denotes the discriminant extension of $h$ over $K$.

\begin{lem}\label{67}
    
    Assume that $D$ is ramified over $K$, so $D=( \frac{a, \pi}{K})$ where $a$ is a unit and $\pi$ is a uniformizer of $K$, and $h$ is unramified with dimension $n$ over $D$. Then 

 $$F \cong K \times K   \  \  \  \iff   \ n   \ \  \    \textit{is even.}$$

Note that the above fact is equivalent to the following: 

$$F \ \textit{is a field}    \iff   \ n  \  \  \  \textit{is odd.}$$ 

In the second case, $F=K(\sqrt{a})$.
\begin{proof}

By Remark \ref{66} we can assume that the skew-hermitian form $h$ has a diagonalization $<u_1 x, u_2 x, \dots, u_n x>$ where $u_i \in {\mathcal{O}_K}^*$.

$F\cong K \times K \  \   \iff  \   \ (-1)^{n} Nrd(u_1 u_2 \dots u_n x^n)\in {K^*}^2    \  \  \  \  \iff  a^n \in {K^*}^2. $

If $n$ is odd, then $$a^n \in {K^*}^2  \  \  \iff     \  \ a \in {K^*}^2  \ \  \Rightarrow   \ \  D   \  \textit{is split},$$
therefore when $F\cong K \times K$, $n$ must be even since we know that $D$ is a division algebra. Conversely, if $n$ is even, then $a^n \in {K^*}^2$ and hence $F\cong K \times K$.

This implies that 

$$F \ \textit{is a field} \     \iff   \ n  \  \  \  \textit{is odd.}$$ 

So when $F$ is a field (and $n$ is automatically odd), we have $F=K(\sqrt{disc  \  h})=K(\sqrt{a})$.
\end{proof}
\end{lem}

\begin{rema}\label{161}
   Assume that $D$ is unramified over $K$. Then any element in the diagonalization of the ramified part of $h$, i.e. $h_1$, has the form $\pi d$, where $d$ is a unit skew element in $D_K$. Since $Nrd(\pi d)=- \pi^2 d^2$ is a unit in $K$ up to squares, the discriminant of $h$ has to be a unit of $K$. Therefore, $F$ (the discriminant extension of $h$ over $K$) is either an unramified quadratic field extension of $K$, or it is isomorphic to $K\times K$.
\end{rema}

\section{Reduction to unramified and ramified skew-hermitian forms under a specific assumption}\label{156}

First, we need to prove a lemma in Galois cohomology. 

%\begin{lem}
 %   \label{18}
%Let $K$ be a field and $G$ an algebraic group defined over $K$ with smooth center $Z$, and $\bar{G}=G/Z$. Then the coboundary map  $$\gamma_K: \bar{G}(K) \longrightarrow H^1(K,Z)$$ is a homomorphism of abstract groups.    
%\begin{proof}
%Since $Z$ is smooth, we have the following exact sequence:

%$$1 \to Z(K^{sep}) \to G(K^{sep})\to \bar{G}(K^{sep}) \to 1.$$

 %   Let $h_1$ be an element in $\bar{G} (K)$. Take a preimage $g_1$ of $h_1$  in $G(K^{sep})$. Then $\gamma_K(h_1)$ is the class of the following cocyle $$a_1: Gal(\bar{K}/K) \longrightarrow Z$$
  %  $$\  \  \ \   \  \ \  \  \ \ \  \   \  \ \  \  \ \ \   \  \ \ s \mapsto {g_1}^{-1} {g_1}^s.$$

   % We can take another arbitrary element $h_2$ in $\bar{G}(K)$, a preimage $g_2 \in G(K)$, and a cocycle $a_2$ in a similar manner. Now $\gamma_K(h_1 h_2)$ will be the class of the  cocycle $s \mapsto (g_1 g_2)^{-1} (g_1 g_2)^s$. The element ${g_1}^{-1} {g_1}^s$ is in $Z$, so it commutes with ${g_2}^{-1}$. Then we have:
%$$(g_1 g_2)^{-1} (g_1 g_2)^s= {g_2}^{-1} {g_1}^{-1} {g_1}^s {g_2}^s=({g_1}^{-1} {g_1}^s)  ({g_2}^{-1} {g_2}^s)$$
 %   which shows that the map $\gamma_K$ is a homomorphism of abstract groups. 
%\end{proof}
%\end{lem}

\begin{lem}\label{85}
Let $K$ be a field and $G$ an algebraic group defined over $K$ with a central subgroup $Z$. Let $\eta_K$  be the induced map $\eta_K: H^1(K,Z) \longrightarrow H^1(K,G)$. Assume that $a$ and $b$ are two cocycles with coefficients in $Z$ such that $\eta_K([a])=\eta_K([b])$. Then  $\eta_K([ab^{-1}])=1$. 
\begin{proof}
    Let $i$ be the inclusion map $i: Z\longrightarrow G$. For any cocycle $c$ with with coefficients in $Z$, the induced map $\eta_K$ sends $[c]$ to $[i\circ c]$. Since $\eta_K([a])=\eta_K([b])$, there exist an element $g\in G$ such that for every $s\in Gal(\bar{K}/K)$ we have $(i\circ a) (s)= g^{-1} ((i\circ b) (s)) g^s $. The subgroup $Z$ is central in $G$, hence for any element $s \in \ Gal(\bar{K}/K)$ the elements $((i\circ a)(s))$ and $((i\circ b)(s))$ commute with every element of $G$. Therefore
$$(i \circ ab^{-1}) (s)=g^{-1} g^s,$$
 which means $i \circ ab^{-1}$ is cohomologous to $1$.
\end{proof}
\end{lem}

Assume that for the skew-hermitian form $h$ with Larmour decomposition  $h_K \cong h_{0,K}\perp h_{1,K}$ there exists an element $\lambda_0 \in L^*$ such that:

\begin{equation}\label{56}
    \lambda_0 h_{0,L}\cong h_{0,L}  \  \    \   \    \  \   \  \   \  \  \   \  \    \lambda_0 h_{1,L}\cong h_{1,L}. 
\end{equation}

We allow one of $h_0$ or $h_1$ to be trivial.

In this section, we will prove a reduction of the norm principle under the assumption that the isometries (\ref{56}) hold. When we apply the main result of this section in the proof of the main theorem, we will take $\lambda_0$ to be the element $\lambda$ introduced in Subsection \ref{112} as a representative for $g_L(u)$ (see Diagram \ref{139}). The following argument is taken from \cite{4}. 

Let $R=O^+(h_0) \times O^+ (h_1)$ be a subgroup of $O^+(h)$, and $\widetilde{R}$ be the preimage of $R$ under the map $Spin(h) \longrightarrow O^+(h)$. 

\[\
\begin{tikzcd}
1 \arrow[r] & \mu_2 \arrow[r] \arrow[d, "id"] & Z \arrow[r] \arrow[d]         & \mu_2 \arrow[r] \arrow[d] & 1 \\
1 \arrow[r] & \mu_2 \arrow[r] \arrow[d, "id"] & \widetilde{R} \arrow[r] \arrow[d] & R \arrow[r] \arrow[d]     & 1 \\
1 \arrow[r] & \mu_2 \arrow[r]                 & Spin(h) \arrow[r]             & O^+(h) \arrow[r]          & 1
\end{tikzcd}
\]\

The above diagram gives rise to the following commutative diagram with exact rows: 
\[\
\begin{tikzcd}
{H^1(L,\mu_2)} \arrow[dd, "id"] \arrow[rr, "f_L"]  &  & {H^1(L,Z)} \arrow[rr, "g_L"] \arrow[dd, "{\alpha_{1,L}}"]              &  & {H^1(L,\mu_2)} \arrow[dd, "{\beta_{1,L}}"] \\
                                                   &  &                                                                        &  &                                            \\
{H^1(L,\mu_2)} \arrow[rr, "e'_L"] \arrow[dd, "id"] &  & {H^1(L,\widetilde{R})} \arrow[rr, "s'_L"] \arrow[dd, "{\alpha_{2,L}}"] &  & {H^1(L,R)} \arrow[dd, "{\beta_{2,L}}"]     \\
                                                   &  &                                                                        &  &                                            \\
{H^1(L,\mu_2)} \arrow[rr, "e_L"]                   &  & {H^1(L,Spin(h))} \arrow[rr, "s_L"]                                     &  & {H^1(L,O^+(h))}                           
\end{tikzcd}
\]\

So we have $\alpha_L = \alpha_{2,L} \circ \alpha_{1,L}$, and $\beta_L = \beta_{2,L} \circ \beta_{1,L}$ (see also Diagram \ref{139}). \\

\begin{lem}\label{51}

Assume that:

- the $H^1$-variant of the norm principle holds for the pair $(Z, \widetilde{R})$, and

- isometries (\ref{56}) hold.

Then $cor_{L/K}(u) \in \Ker (\alpha_K)$. 
\begin{proof}
   The proof involves chasing elements in the diagram above. Let $v=\alpha_{1,L}(u)$. We have $\lambda_0 h_{0,L} \cong h_{0,L}$ and $\lambda_0 h_{1,L} \cong h_{1,L}$, therefore, $\beta_{1,L}([\lambda_0])=1$. This shows that $v \in \Ima (e'_L)$, so there exists $a \in L$ such that $e'_L([a])=v$.

The element $b=f_L([a])^{-1} u \in H^1(L,Z)$ is in the kernel of the map $\alpha_{1,L}$ (by Lemma \ref{85}), and therefore it is in the kernel of the map $\alpha_L$. But $u$ is also in the kernel of $\alpha_L$, hence $f_L([a])^{-1} \in \Ker (\alpha_L)$, i.e. $a$ is a spinor norm for $h$. By Theorem \ref{8} we have $cor_{L/K}(f_L([a])^{-1}) \in \Ker (\alpha_K)$. 

By our assumption, the norm principle holds for the pair $(Z, \widetilde{R})$, therefore $cor_{L/K}(b) \in \Ker (\alpha_{1,K}) \subseteq \Ker (\alpha_K)$. Again, by Lemma \ref{85}, $cor_{L/K}(u) \in \Ker (\alpha_K)$.\end{proof}
\end{lem}

\begin{rema}\label{55}
    Lemma \ref{51} shows that to prove Theorem \ref{40}, it will be enough  to show the norm principle for the spinor groups of $h_0$ and $ h_1$ under the following additional hypotheses:\\ \\
- isometries (\ref{56}) hold,\\ 
- the skew-hermitian form $h$ has a nontrivial decomposition $h= h_0 \perp h_1$ (i.e. $dim \ h_0 \ne 0, \ dim \ h_1 \ne 0$).

Therefore, we have reduced  the norm principle to the following two cases: 

\begin{itemize}
    \item unramified skew-hermitian forms (i.e. $h_K \cong h_{0,K}$), and 
    \item ramified (i.e. totally ramified) skew-hermitian forms (i.e. $h_K \cong h_{1,K}$).
\end{itemize} 
\end{rema}
\begin{rema}\label{62}
  Note that if $D$ is unramified over $K$, then since $$Spin(h_{1,K}) = Spin(\frac{1}{\pi} h_{1,K}),$$ by Remark \ref{55} we can assume that $h$ is unramified over $K$ (the skew-hermitian form $\frac{1}{\pi} h_{1,K}$ is unramified).
\end{rema}

\section{A reduction when $h$, $D$, and $\lambda$ are unramified over $L$}\label{157}

Throughout this section we assume that $h$, $D$, and $\lambda$ are unramified over $L$  (by $\lambda$ being unramified we mean it is a unit over $L$ up to squares) and will prove a reduction of the main theorem under this assumption. In this case, we may assume that the group $Spin(h)$ and its center $Z$ are actually defined over $\mathcal{O}_L$ and by base change we obtain $Spin(h)$ and $Z$ over $L$. Consider the following commutative exact diagram, induced by the base change $\mathcal{O}_L \longrightarrow L$:

\includegraphics[width=\textwidth]{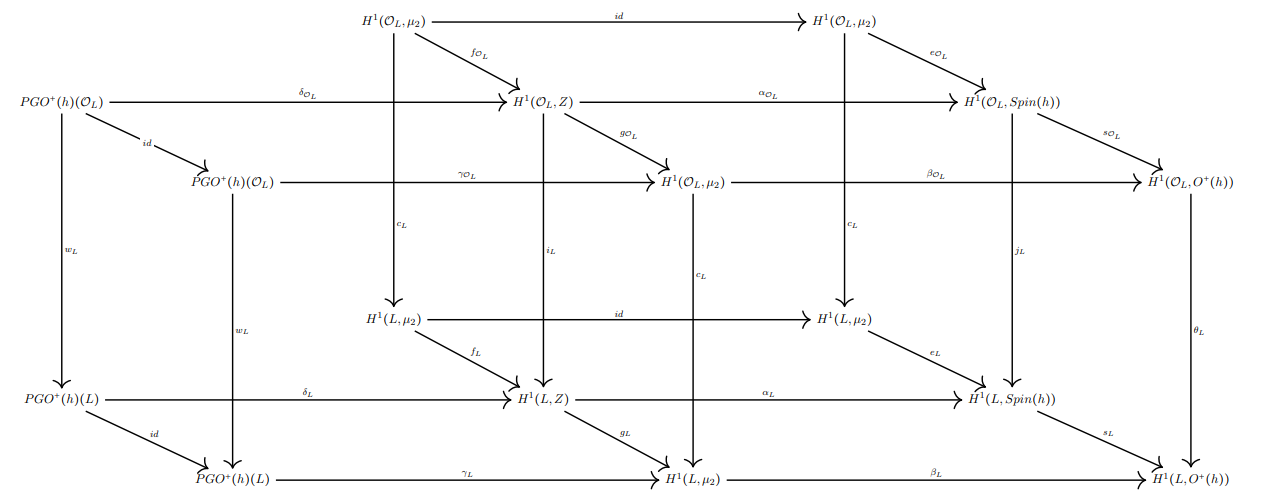}

In the above diagram, the maps $c_L$, $i_L$, and $w_L$ are injective group homomorphisms, and the maps $j_L$ and $\theta_L$ have trivial kernels by \cite{17}.

 We need the following two lemmas  in case $h$, $D$, and $\lambda$ are unramified over $L$:\\

\begin{lem}\label{68} 
 If $h$, $D$, and $\lambda$ are unramified over $L$, we can write the element $u$ as a product $u=i_L(u_1) f_L(u_2)$, for some $u_1 \in H^1(\mathcal{O}_L, Z)$,  $u_2 \in H^1(L, \mu_2)$, such that

    $g_L(i_L(u_1))=[\lambda]$, 
    
    $i_L(u_1)\in \Ker (\alpha_L)$, 
    
    $e_L(u_2)=1$, and 
    
    $\alpha_{\mathcal{O}_L}(u_1)=1$.\\

    %(c) If $e(L/K)=2$ (hence $\pi '=\sqrt{c\pi}$, i.e. case E), then $cor_{L/K}(i_L(u_1))=1$.\\

\begin{proof}

Consider the element $\lambda {{\mathcal{O}_L}^*}^2 \in H^1(\mathcal{O}_L, \mu_2)$. Then

 \begin{equation*} 
\begin{split}
\theta_L(\beta_{\mathcal{O}_L}(\lambda {{\mathcal{O}_L}^*}^2)) & =  \beta_L(c_L(\lambda {{\mathcal{O}_L}^*}^2))  \\
 & = \beta_L(\lambda {L^*}^2) \\
 & = \beta_L(g_L(u))  \\
 &= s_L(\alpha_L(u)) \\
 &= 1.
\end{split}
\end{equation*}

By \cite{17}, $\Ker (\theta_L)$ is trivial. So we have $\beta_{\mathcal{O}_L}(\lambda {{\mathcal{O}_L}^*}^2)=1$, hence $\lambda {{\mathcal{O}_L}^*}^2 \in \Ima (\gamma_{\mathcal{O}_L})$, and there exists an $\mathcal{O}_L$-point $s$ in $PGO^{+} (h'')(\mathcal{O}_L)$ such that $\gamma_{\mathcal{O}_L}(s)=\lambda {{\mathcal{O}_L}^*}^2$. Let $u_1:=\delta_{\mathcal{O}_L}(s) \in H^1(\mathcal{O}_L, Z)$. Then

 \begin{equation*} 
\begin{split}
g_L(i_L(u_1)) & =  c_L(g_{\mathcal{O}_L}(u_1)) \\
 & = c_L(g_{\mathcal{O}_L}(\delta_{\mathcal{O}_L}(s))) \\
 & = c_L(\gamma_{\mathcal{O}_L}(s)) \\
 &= c_L(\lambda {{\mathcal{O}_L}^*}^2)  \\
 &= \lambda {L^*}^2 .
\end{split}
\end{equation*}

Since $g_L(u)=[\lambda]$, so by Lemma \ref{85} we have $u i_L(u_1)^{-1}\in \Ker (g_L)=\Ima (f_L)$, hence there exists $u_2\in H^1(L, \mu_2)$ such that $f_L(u_2)=u i_L(u_1)^{-1}$, therefore $u=i_L(u_1) f_L(u_2)$. 

We have $u_1 \in \Ima (\delta_{\mathcal{O}_L})=\Ker (\alpha_{\mathcal{O}_L})$, so $\alpha_{\mathcal{O}_L}(u_1)=1$. Hence 

$$\alpha_L(i_L(u_1))=j_L(\alpha_{\mathcal{O}_L}(u_1))=j_L(1)=1,$$so $i_L(u_1) \in \Ker (\alpha_L)$. 

Since $u\in \Ker (\alpha_L)$, by Lemma \ref{85}, $f_L(u_2)=u i_L(u_1)^{-1}\in \Ker (\alpha_L)$. Hence $e_L(u_2)=\alpha_L(f_L(u_2))=1$. 

Finally, since $j_L(\alpha_{\mathcal{O}_L}(u_1))=\alpha_L(i_L(u_1))=1$ and the kernel of the map $j_L$ is trivial by \cite{17}, we have $\alpha_{\mathcal{O}_L}(u_1)=1$. 
\end{proof}
\end{lem}

\begin{lem}\label{69}
Let $h$, $D$, and $\lambda$ be unramified over $L$. Recall the elements $u_1$ and $u_2$ from Lemma \ref{68}. Then
$$cor_{L/K}(i_L(u_1)) \in \Ker (\alpha_K) \Rightarrow  cor_{L/K}(u) \in \Ker (\alpha_K).$$
\begin{proof}
    Since $e_L(u_2)=1$, by Theorem \ref{8} we have 
$cor_{L/K}(u_2)\in \Ker (e_K)$, so  $cor_{L/K}(f_L(u_2))\in \Ker (\alpha_K)$. Then
$$cor_{L/K}(u)=cor_{L/K}(i_L(u_1) f_L(u_2))=cor_{L/K}(i_L(u_1)) cor_{L/K}(f_L(u_2)),$$and by Lemma \ref{85} we have $cor_{L/K}(u) \in \Ker (\alpha_K)$. 
\end{proof}
\end{lem}

\section{A technical lemma}\label{128}

We will apply the following lemma in the proof of the main theorem in Section \ref{159}.

\begin{lem}\label{88}
    Consider the linear algebraic groups $G_1$, $G_2$, $G_3$, $G_4$, $G_5$, $G_6$ and $G_7$ all defined over the field $K$, which fit into the following commutative diagram with exact rows.

\[\
\begin{tikzcd}
1 \arrow[r] & G_1 \arrow[r] \arrow[d, "id"] & G_2 \arrow[r] \arrow[d] & G_3 \arrow[r] \arrow[d] & 1 \\
1 \arrow[r] & G_1 \arrow[r] \arrow[d, "id"] & G_4 \arrow[r] \arrow[d] & G_5 \arrow[r] \arrow[d] & 1 \\
1 \arrow[r] & G_1 \arrow[r]                 & G_6 \arrow[r]           & G_7 \arrow[r]           & 1
\end{tikzcd}
\]\

Assume that all the vertical arrows are inclusions, $G_2$ is central in $G_6$ (and hence in 
$G_4$), and $G_3$ is central in $G_7$ (and hence in 
$G_5$). 

Let $L/K$ be a finite separable field extension. The above diagram gives rise to the following commutative diagram with exact horizontal rows:

\includegraphics[width=\textwidth]{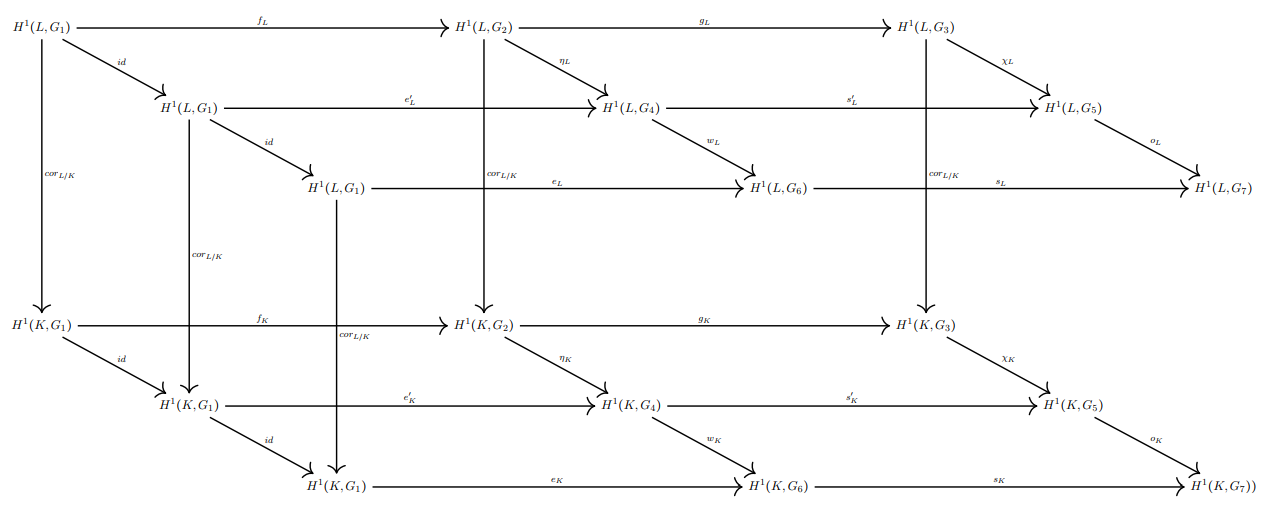}

Let $\alpha:= w \circ \eta$, and $u\in \Ker (\alpha_L)$. 

Assume that:

$\bullet$ (a) $u \in \Ker (\chi_L \circ g_L)$.

$\bullet$ (b) $cor_{L/K} (\Ker ( \chi_L))\subset \Ker (\chi_K)$.

$\bullet$ (c) $cor_{L/K} (\Ker ( e_L))\subset \Ker (e_K)$.

$\bullet$ (d) $cor_{L/K} (\Ker (\eta_L))\subset \Ker (\eta_K)$.

Then $cor_{L/K}(u)\in \Ker ( \alpha_K)$. 
\begin{proof}

 Put $v:=\eta_L(u)$ and $\lambda:=g_L(u)$. By assumption (a), $u \in \Ker (\chi_L  \circ g_L)= \Ker ({s_L}' \circ \eta_L)$, so $s'_L(v)=1$. Let $u_2$ be a preimage of $v$ under ${e_L}'$. By (b) we have $\chi_K(cor_{L/K}(\lambda))=1$. Let $u_1:=u{f_L(u_2)}^{-1}$. We have 

$$e_L(u_2)=w_L(e'_L(u_2))=w_L(v)=\alpha_L(u)=1.$$

So $u_2 \in \Ker (e_L)$ and assumption (c) implies that $cor_{L/K}(u_2)\in \Ker (e_K)$ and therefore $cor_{L/K}(f_L(u_2))\in \Ker (\alpha_K)$.

Lemma \ref{85} shows that $u_1 \in \Ker (\eta_L)$, because $\eta_L(u)=\eta_L(f_L(u_2))=v$. By assumption (d) we have $\eta_K(cor_{L/K}(u_1))=1$, and $\alpha_K(cor_{L/K}(u_1))=w_K \circ \eta_K (cor_{L/K}(u_1))=1$. 

Therefore $cor_{L/K}(u_1)$ and $cor_{L/K}(f_L(u_2))$ are both in $\Ker (\alpha_K)$, and by Lemma \ref{85}, we have $cor_{L/K}(u)=cor_{L/K}(u_1 f_L(u_2)) \in \Ker (\alpha_K)$. 
\end{proof}
\end{lem}

\section{Proof of the main theorem}\label{184}

In this section, we will give a proof of the main theorem, i.e. Theorem \ref{40}, in the following 4 cases separately. All the notations introduced in Section \ref{111} (and then used in Sections
 \ref{156} and \ref{157}) will be used in the proof. 

\begin{itemize}
    \item \textbf{Case 1}: $L/K$ is an unramified extension, and $D_K$ is unramified over $K$.
    \item \textbf{Case 2}: $L/K$ is an ramified extension, and $D_K$ is unramified over $K$.

    \item \textbf{Case 3}: $L/K$ is an ramified extension, and $D_K$ is ramified over $K$.
    \item \textbf{Case 4}: $L/K$ is an unramified extension, and $D_K$ is ramified over $K$.
\end{itemize}

In each case we have several subcases, depending on the discriminant extension $F$, the element $\lambda$, ramification of the skew-hermitian form $h$, splitting of $D$ over $L$, and the parity of $n$. Recall that by Remark \ref{4000}, we assume that \textbf{$h$ is anisotropic over $L$ in all cases}. We list all the subcases below.

\includegraphics[width=\textwidth]{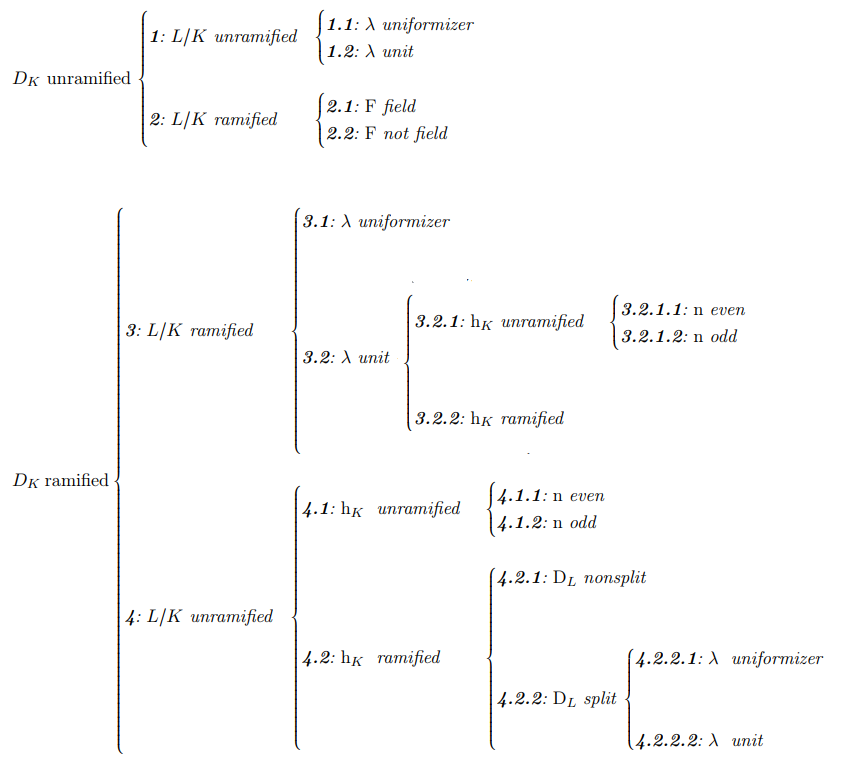}

\subsection{Proof of the main theorem: part 1}\label{141} In this subsection we will prove the main theorem in case 1: $L/K$ is unramified  and $D_K$ is unramified. 

The elements $a$, $b$, and $t$ are units; $x^2=a, \ y^2=b,$ and  $z=xy=-yx$ for $x, y, z \in D_K$. If $D$ splits over $L$, then by Lemma \ref{94}, we assume that $t=a$.

 Let us consider the following subcases: 
 \begin{itemize}
     \item \textbf{Subcase 1.1}: $\lambda$ is a uniformizer of $L$.

      \item \textbf{Subcase 1.2}: $\lambda$ is a unit of $L$.
 \end{itemize}

 \textbf{Subcase 1.1: $\lambda$ is a uniformizer of $L$}.

We prove a fact about the R-triviality of the adjoint group $PGO^+(h)$ under a specific assumption.

\begin{lem}\label{118}
    Let $E$ be a complete discretely valued field and $char \ E\ne 2$.  Assume that $Q$ is an unramified quaternion division algebra over $E$ and $\mathfrak{h}$ is an anisotropic skew-hermitian form over $Q$ with respect to the canonical involution on $Q$. Let $\mu$ be a uniformizer of $E$ such that we have the $E-$isometry $\mu \mathfrak{h} \cong \mathfrak{h}$. Then $PGO^+(\mathfrak{h})(E)$ is $R-$trivial.
    \begin{proof}
        If $\mathfrak{h}$ is unramified or ramified, then we will have $\mu \mathfrak{h}_0 \cong \mathfrak{h}_0$ or  $\mu \mathfrak{h}_1 \cong \mathfrak{h}_1$ which is impossible by Theorem \ref{54}. Hence, we assume that $\mathfrak{h}$ is neither unramified, nor ramified over $E$. Let  $\mathfrak{h}=\mathfrak{h}_0 \perp \mathfrak{h}_1$ be  Larmour decomposition of $\mathfrak{h}$ over $E$ (so $dim \ \mathfrak{h}_0$ and $dim \ \mathfrak{h}_1$ are nonzero).

        Since $\mu$ is a multiplier of $\mathfrak{h}$ over $E$, we have the following isometries (over $E$):

$$ \mathfrak{h}_{0,E} \perp \mathfrak{h}_{1,E}  \cong \mathfrak{h}_E \cong \mu \mathfrak{h}_E \cong \mu \mathfrak{h}_{0,E} \perp \mu \mathfrak{h}_{1,E}.$$

The element $\mu$ is a uniformizer in $E$, hence by Theorem \ref{54}, we have the $E-$ isometry $\mathfrak{h}_{1,E} \cong \mu \mathfrak{h}_{0,E}$. Furthermore, $h_E \cong \mathfrak{h}_{0,E} \perp \mu \mathfrak{h}_{0,E}$.

Let $\theta \in G(\mathfrak{h})(E)$ (i.e. $\theta$ be a multiplier for $\mathfrak{h}$ over $E$). Up to squares in $E^*$, either $\theta$ is a unit in $E^*$, or $\theta = \theta' \mu$ for a unit $\theta'$ in $E^*$. Since $PGO^+ (\mathfrak{h})(E)/R \cong G(\mathfrak{h})/ {E^*}^2 Hyp(\mathfrak{h})$ (Theorem \ref{33}), we need to show
$\theta \in {E^*}^2 Hyp(h)$.

The totally ramified extension $E(\sqrt{-\mu})$ splits $\mathfrak{h}$ and $\mu = N_{E(\sqrt{-\mu})/E}(\sqrt{-\mu})$. Thus $\mu \in Hyp(\mathfrak{h})$. 

Case 1: Let $\theta$ be a unit in $E^*$. The $E-$isometry $\theta \mathfrak{h}_E \cong \mathfrak{h}_E$ implies that  

$$\theta \mathfrak{h}_{0,E} \perp \theta \mu \mathfrak{h}_{0,E} \cong 
\mathfrak{h}_{0,E} \perp \mu  \mathfrak{h}_{0,E},$$
hence by Theorem \ref{54}, we have $\mathfrak{h}_{0,E}\cong \theta \mathfrak{h}_{0,E}$.

Consider the field $E_1=E(\sqrt{-\mu 
 \theta})$, which is a complete discretely valued field with the same residue field as $E$. Over $E_1$, we have the isometries 

$$\mathfrak{h}_{E_1} \cong \theta \mathfrak{h}_{E_1} \cong  \theta \mathfrak{h}_{0,E_1} \perp \theta \mu \mathfrak{h}_{0,E_1} \cong \mathfrak{h}_{0,E_1} \perp \theta \mu \mathfrak{h}_{0,E_1} \cong \mathfrak{h}_{0,E_1} \perp -\mathfrak{h}_{0,E_1}.$$

So $\mathfrak{h}$ is hyperbolic over $E_1$, and $\theta \mu = N_{E_1/E}(\sqrt{- \theta \mu}) \in Hyp(\mathfrak{h})$. Since $\mu \in Hyp(\mathfrak{h})$, we have $\theta \in Hyp(\mathfrak{h})$, which shows that $PGO^+ (\mathfrak{h})(E)$ is $R-$trivial by Theorem \ref{33}.

Case 2: Let $\theta=\theta' \mu$  for a unit $\theta'$ in $E^*$. Then $\theta' = \theta {\mu}^{-1}\in G(\mathfrak{h})$, and repeating the same argument as in case 1 for $\theta'$ shows that $PGO^+ (\mathfrak{h})(E)$ is $R-$trivial.
    \end{proof}
\end{lem}

The following lemma is going to be used in the current section, and also later in the proof of case 3.

\begin{lem}\label{39}
Assume that $\lambda$ is a uniformizer in $L$, and  furthermore one of the following assumptions holds:

\begin{itemize}
    \item $L/K$ is unramified and $D_K$ is unramified, or

    \item $L/K$ is ramified and $D_K$ is ramified.
\end{itemize}

 Then $PGO^+(h)(L)$ is $R-$trivial.
   \begin{proof}
   Assume that $D_L$ is an $L$- division algebra. Then $D_L$ is  unramified over $L$. By Lemma \ref{118}, $PGO^+(h)(L)$ is $R-$trivial (recall that $h$ is anisotropic over $L$, by Remark \ref{4000}).

Now assume that $D_L$ is split. The skew-hermitian form $h$ over $D_L$ corresponds to a quadratic form $q$ over $L$ with $\lambda \in G(q)(L)$ (the group of multipliers of $q$ over $L$) by Morita equivalence, once we fix an isomorphism $M_2(L)\cong D_L$ (Morita equivalence preserves multipliers). Since $\lambda$ is not a square in $L$, the Springer  decomposition of quadratic forms implies that $q\cong q' \otimes <1, \lambda>$ for an unramified quadratic form $q'$ defined over $L$. Then by \cite[Lemma~4.2]{4},
 $PGO^+ (q)(L)$ is $R-$ trivial, which also means that $PGO^+ (h)(L)$ is $R-$ trivial.
\end{proof}
\end{lem}

\textbf{Proof of the main theorem in subcase 1.1:}\\ 

\begin{proof}
    By Lemma \ref{39} $PGO^+(h)(L)$ is $R-$trivial. Then Theorem \ref{7} implies that $cor_{L/K}(u) \in \Ker (\alpha_K)$.\\ 
\end{proof}

\textbf{Subcase 1.2: $\lambda$ is a unit of $L$.}\\

\begin{lem}\label{158}
In subcase 1.2, we have the following isometries over $L$:

$$h_{0,L}\cong \lambda h_{0,L},   \  \   \ \textit{and}  \  \  \  \  \  h_{1,L}\cong \lambda h_{1,L}.$$

\begin{proof}
If one of the skew-hermitian forms $h_0$ or $h_1$ is trivial, then we already have the desired result since $h_L \cong \lambda h_L$. So assume that both of them are nontrivial.

First assume that $D_L$ is an $L$-division algebra. Then we have $\lambda h_{0,L} \perp \lambda h_{1,L} = \lambda h_L \cong h_L = h_{0,L} \perp h_{1,L}$. 
Every element in the diagonalization of $h_0$ (resp. $h_1$) has value $0$ (resp. $1$) over $L$. By Theorem \ref{54} we have $h_{0,L}\cong \lambda h_{0,L}$ and   $h_{1,L}\cong \lambda h_{1,L}$.

Now, assume that $D$ splits over $L$. Fix an isomorphism $M_2(L)\cong D_L$. By Morita equivalence, the forms $h_{0,L}$ and $h_{1,L}$ correspond to quadratic forms $q_{0,L}$ and $q_{1,L}$ over $L$, respectively. Changing the isomorphism $M_2(L)\cong D_L$ affects Morita equivalence by a scalar. We choose the isomorphism such that the ramification of the forms are preserved under Morita equivalence, i.e. $q_0$ is unramified and $q_1$ is ramified. The isometry $h_L \cong \lambda h_L$ implies that $q_L \cong \lambda q_L$, so we have $\lambda q_{0,L} \perp \lambda q_{1,L} \cong q_{0,L} \perp q_{1,L}.$ Note that the elements in the diagonalization of $q_0$ (resp. $q_1$) have value $0$ (resp. $1$) over $L$. Hence 
by Springer's decomposition theorem for quadratic forms over complete discretely valued fields we have $q_{0,L}\cong \lambda q_{0,L}$ and $  q_{1,L}\cong \lambda q_{1,L}.$
Finally, by Morita equivalence $h_{0,L}\cong \lambda h_{0,L}$ and $ h_{1,L}\cong \lambda h_{1,L}.$
\end{proof}
\end{lem}

 By Remarks \ref{55}, \ref{62} and Lemma \ref{158}, we may assume that $h=h_0$, i.e. $h$ is unramified. Let $\bar{h}$ be the residue skew-hermitian form defined over $(\overline{D_K}, \overline{\tau})$.  \\

\textbf{Proof of the main theorem in subcase 1.2:}

\begin{proof}
    Recall that $k$ denotes the residue field of $K$ and $l=k(\sqrt{\bar{t}})$ is the residue field of $L$. Note that even though we do not have a canonical valuation on $D_L$ if it splits, but we still have a skew-hermitian form $\bar{h}$ defined over $l$, obtained by base change of $\bar{h}$ from $k$ to $l$. 

Recall the element $u_1$ from Lemma \ref{68}. By specializing to the residue field we get 

$$\overline{u_1}\in \Ker ({\alpha}_l: H^1(l, \bar{Z}) \longrightarrow H^1(l, Spin(\bar{h}))).$$

Since the norm principle holds for $Spin(\bar{h})$ over $l/k$, by assumption we get $$N_{l/k}(\overline{u_1})\in \Ker ({\alpha}_k: H^1(k, \bar{Z}) \longrightarrow H^1(k, Spin(\bar{h}))).$$

Consider a principal homogeneous space $\mathcal{X}$ under $Spin(\bar{h})$ over $k$ whose isomorphism class is $\alpha_k(N_{l/k}(\overline{u_1}))$. As a variety, $\mathcal{X}$  is isomorphic to the  variety $Spin(\bar{h})$ over $k$ and has a $k-$rational point. By Hensel's lemma (see \cite[Section~2.3, Prop~5]{1001}), this rational point lifts to a $K-$rational point on a principal homogeneous space under $Spin(h)$ over $K$ whose isomorphism class is
$\alpha_K(cor_{L/K}(i_L(u_1)))$. Therefore, the image of the element $cor_{L/K}(i_L(u_1))$ under the map $H^1(K,Z) \to H^1(K, Spin(h))$ is the distinguished element of $H^1(K, Spin(h))$. By Lemma \ref{69}, $cor_{L/K}(u) \in \Ker (\alpha_K)$.\\ 

\end{proof}

\subsection{Proof of the main theorem: part 2}\label{142} In this subsection we will prove the main theorem in case 2: when $L/K$ is ramified and $D_K$ is unramified.

\begin{lem}\label{93}
    If $D_K$ is unramified over $K$, and $D_L$ splits, then $L/K$ must be unramified.
    \begin{proof}
        By Lemma \ref{94}, $L$ can be embedded over $K$ in $D_K$. Since $D_K$ is unramified over $K$, any maximal subfield, in particular $L$, has to be unramified over $K$ too. 
    \end{proof}
\end{lem}

By Lemma \ref{93}, $D_L$ is an unramified $L-$division algebra in case 2.

 The skew-hermitian form $h$ can be neither unramified nor ramified over $K$, i.e. both $h_0$ and $h_1$ could be nontrivial. However, $h_L$ is unramified over $L$. 

Recall that by Remark \ref{161}, $F$  is either an unramified quadratic field extension of $K$, or it is isomorphic to $K\times K$. Now we give the proof of the main theorem in these two different cases separately, first when $F$ is a field, and then in the case $F=K \times K$. \\ 

\textbf{Subcase 2.1:} $F$ is a field. 

In this case, $M= F \otimes_K L$ is a ramified quadratic field extension of $F$.\\ 

\textbf{Proof of the main theorem in subcase 2.1}: 

\begin{proof}
    
Assume that $n$ is even. Then $H^1(K,Z)= F^* / {F^*}^2$ and $H^1(L,Z)= M^* / {M^*}^2$. The field extension $F/K$ is unramified,  hence the field extension $M/L$ is also unramified,  and the element $\lambda$ is a unit in $L^*$ up to squares (recall that the element $[\lambda]$ is image of the element $u$ under the map $g_L: H^1(L,Z) \to H^1(L,\mu_2)$, see Diagram \ref{139} and the explicit formula for the map $g_L$ is Subsection \ref{154}). Recall the element $u_1$ from Lemma \ref{68}. The element $u_1$ has a representative $v$ in ${\mathcal{O}_M}^*$. Since the extension $M/F$ is ramified, $N_{M/F}(v)$ is a square in $F^*$, so $cor_{L/K} (i_L(u_1))=1 \in H^1(K,Z)$. Hence $cor_{L/K}(u)\in \Ker \ \alpha_K$ by Lemma \ref{69}.

Now assume that $n$ is odd. By the explicit description of the map $g_L$ given in Subsection \ref{154}, it follows that the element $\lambda$ is a unit up to squares in $L^*$.
Recall the element $u_1$ from Lemma \ref{68}.
The element $u_1$ has a representative $(q,p) \in {\mathcal{O}_L}^* \times {\mathcal{O}_M}^*$. Then $cor_{L/K}(u_1)=[(N_{L/K}(q), N_{M/F}(p))]$. Since both extensions $M/F$ and $L/K$ are ramified, both elements $N_{L/K}(q)$ and $N_{M/F}(p)$ are squares in $K^*$ and $F^*$, respectively. Hence we have $cor_{L/K} (i_L(u_1))=1 \in H^1(K,Z)$ and  $cor_{L/K}(u)\in \Ker \ \alpha_K$ by Lemma \ref{69}. \\
\end{proof}

\textbf{Subcase 2.2:} $F\cong K \times K$.  

In this case, $M\cong L \times L$.\\ 

\textbf{Proof of the main theorem in subcase 2.2}:

\begin{proof}
    
Assume that $n$ is even. Then $H^1(K,Z)= (K^* \times K^*) / ({K^*}^2 \times {K^*}^2)$ and $H^1(L,Z)= (L^* \times L^*) / ({L^*}^2 \times {L^*}^2)$. The element $u$ has a representative $(v,w)$ in ${L}^* \times {L}^*$. Since the extension $L/K$ is ramified, $cor_{L/K}(u)=[(N_{L/K}(v), N_{L/K}(w))]=[(1,1)]=1 \in H^1(K,Z)$. Hence $cor_{L/K}(u)\in \Ker \ \alpha_K$.

 Now assume that $n$ is odd. Then the element $u$ has a representative $(q,(p,p')) \in {L}^* \times ({L}^* \times {L}^*)$. Then $cor_{L/K}(u)=[(N_{L/K}(q), (N_{L/K}(p), N_{L/K}(p')))]$. Since $L/K$ is a ramified extension, all elements $N_{L/K}(q)$, $N_{L/K}(p)$ and $N_{L/K}(p')$ are squares in $K^*$. Hence we have $cor_{L/K} (u)=1 \in H^1(K,Z)$ and  $cor_{L/K}(u)\in \Ker \ \alpha_K$.\\

\end{proof}

\subsection{Proof of the main theorem: part 3}\label{144} The goal in this subsection is to prove the main theorem in case 3, when $L/K$ is ramified and $D_K$ is ramified. Note that $D_L$ is then unramified over $L$.

 Recall that $D=(\frac{a,\pi}{K})$ and $L=K(\sqrt{c\pi})$. The element $a$ is unit over $K$; $x^2=a, \ y^2=\pi,$ and   $z=xy=-yx$ for $x, y, z \in D_K$. 

 Let us consider he following subcases: 

 \begin{itemize}
     \item \textbf{Subcase 3.1}: $\lambda$ is a uniformizer of $L$.

      \item \textbf{Subcase 3.2}: $\lambda$ is a unit of $L$.
 \end{itemize}

 \textbf{Proof of the main theorem in subcase 3.1}:

\begin{proof}
If $\lambda$ is a uniformizer of $L$, then by Lemma \ref{39} $PGO^+(h)(L)$ is $R-$trivial, and Theorem \ref{7} implies that $cor_{L/K}(u) \in \Ker (\alpha_K)$.\\ 
\end{proof}

In the rest of this subsection, we prove the main theorem in subcase 3.2: we assume that $\lambda$ is a unit of $L$.

\begin{lem}\label{160}
In subcase 3.2, we have the following isometries over
L:

$$h_{0,L}\cong \lambda h_{0,L},   \  \   \ \textit{and}  \  \  \  \  \  h_{1,L}\cong \lambda h_{1,L}.$$

\begin{proof}
If one of the skew-hermitian forms $h_0$ or $h_1$ is trivial, then we already have the desired result since $h_L \cong \lambda h_L$. So assume that both of them are nontrivial.

First assume that $D_L$ is an $L$-division algebra. Then we have $\lambda h_{0,L} \perp \lambda h_{1,L} = \lambda h_L \cong h_L = h_{0,L} \perp h_{1,L}$. Every element in the diagonalization of $h_0$ (resp. $h_1$) has value $0$ (resp. $\frac{1}{2}$) over $L$. By Theorem \ref{54} we have $h_{0,L}\cong \lambda h_{0,L}$ and   $h_{1,L}\cong \lambda h_{1,L}$.

Now, assume that $D$ splits over $L$. Fix an isomorphism $M_2(L)\cong D_L$. By Morita equivalence, the forms $h_{0,L}$ and $h_{1,L}$ correspond to quadratic forms $q_{0,L}$ and $q_{1,L}$ over $L$, respectively. Changing the isomorphism $M_2(L)\cong D_L$ affects Morita equivalence by a scalar. We choose the isomorphism such that the ramification of the forms are preserved under Morita equivalence, i.e. $q_0$ is unramified and $q_1$ is ramified. The isometry $h_L \cong \lambda h_L$ implies that $q_L \cong \lambda q_L$, so we have $\lambda q_{0,L} \perp \lambda q_{1,L} \cong q_{0,L} \perp q_{1,L}.$ Note that the elements in the diagonalization of $q_0$ (resp. $q_1$) have value $0$ (resp. $\frac{1}{2}$) over $L$. Hence 
by Springer's decomposition theorem for quadratic forms over complete discretely valued fields we have $q_{0,L}\cong \lambda q_{0,L}$ and $  q_{1,L}\cong \lambda q_{1,L}.$
Finally, by Morita equivalence $h_{0,L}\cong \lambda h_{0,L}$ and $ h_{1,L}\cong \lambda h_{1,L}.$ 
\end{proof}
\end{lem}

  By Remark \ref{55} and Lemma \ref{160}, we may assume that either $h=h_0$ or $h=h_1$, i.e. $h$ is either unramified or ramified.\\

  \textbf{Subcase 3.2.1: $h_K$ is unramified}

The group $Spin(h)$ and its center $Z$ are actually defined over $\mathcal{O}_K$ and by base change we also get the original groups over $K$. Note that $D_L$ is either an unramified $L-$division algebra, or it is a split algebra.\\ 

\textbf{Proof of the main theorem in subcase 3.2.1:}

Recall that the discriminant extension of $K$ is $F$ (and the discriminant extension of $L$ is $M$), and also the element $u_1$ from Lemma \ref{68}.\\

\textbf{Subcase 3.2.1.1: when $n$ is even}\\

\begin{proof}
    
In this case be Lemma \ref{67}, we have $F=K \times K$, hence $M=L \times L$ and $H^1(L, Z)=H^1(L, \mu_2) \times H^1(L, \mu_2)=({L^*}/{L^*}^2)\times ({L^*}/{L^*}^2)$. Since the quadratic extension $L/K$ is totally ramified, any element in ${L^*}/{L^*}^2$ has a representative in ${K^*}/{K^*}^2$, so $cor_{L/K}(i_L(u_1))$ is trivial in $H^1(K, Z)=({K^*}/{K^*}^2)\times ({K^*}/{K^*}^2)$, and in particular it implies that $cor_{L/K}(u)\in \Ker \ \alpha_K$ by Lemma \ref{69}.
\end{proof}

\textbf{Subcase 3.2.1.2: when $n$ is odd}\\

\begin{proof}
    
In this case, $M=L(x)$ and $H^1(L, Z)= U(L)/U_0(L)$ (see Section \ref{154} and Lemma \ref{67}). Since the quadratic extension $L/K$ (resp.  $L(x)/K(x)$) is totally ramified, any element in $L^* / (L^*)^2$ (resp.  ${L(x)^*}/{L(x)^*}^2$) has a representative in $K^* / (K^*)^2$ (resp. ${K(x)^*}/{K(x)^*}^2$). So we can represent $u_1$ as $u_1=[(q,p)]$ where $q \in {\mathcal{O}_K}^*$ and $p \in \mathcal{O}_{{K(x)}^*}$, and then  $cor_{L/K}(i_L(u_1))=[(q^2, p^2)]=[(1,1)]=1$, hence  $cor_{L/K}(i_L(u_1))\in \Ker (\alpha_K)$ and we are done by Lemma \ref{69}. 
\end{proof}

\textbf{Subcase 3.2.2: $h_K$ is ramified}

We will construct an unramified hermitian form $h'$ in this subcase, and then will reduce the main question from $h$ to $h'$. 
By Remark \ref{66} we can assume that $$h=<(\beta_1 + \gamma_1 x)y, \dots, (\beta_n + \gamma_n x)y>,$$where $\beta , \dots, \beta_n, \gamma_1, \dots, \gamma_n \in {\mathcal{O}_K}$ and the elements $\beta_i + \gamma_i x$ are all units.

Consider the following orthogonal involution on $D$:

$$\tau_{y}: D   \longrightarrow D$$

$$ \   \  \   \  \   \   \   \  \   \  \  \   \  \  \      s \mapsto  y^{-1} \tau(s) y.$$\\

Elements in the subfield $K(x)$ are invariant under $\tau_y$. The following is a unit hermitian form over $(D, \tau_y)$ (since the elements of the maximal subfield $K(x)$ are invariant under $\tau_y$):

$$h'=<\beta_1 + \gamma_1 x, \dots, \beta_n + \gamma_n x>.$$
  
Let $H=diag(\beta_1 + \gamma_1 x, \dots, \beta_n + \gamma_n x)\in M_n(D)$ be the diagonal matrix representing $h'$. Then the matrix $Hy$ represents $h$.

\begin{cor}\label{101}
    The norm principle for  $Spin(h')$ is equivalent to the norm principle for   $Spin(h)$.
\begin{proof}

Let $\sigma_h$ and $\sigma_{h'}$ be the involutions on $M_n (D)$, adjoint to $h$ and $h'$, respectively. For all $M\in M_n(D)$ we have

 \begin{equation*} 
\begin{split}
\sigma_{h'}(M) & =  H \tau_y(M^t)H^{-1} \\
 & = Hy \tau(M^t)y^{-1}H^{-1} \\
 & = Hy \tau(M^t) (Hy)^{-1} \\
 & = \sigma_h(M). 
\end{split}
\end{equation*}

Hence the adjoint involutions $\sigma_h$ and $\sigma_{h'}$ on $M_n (D)$ are identical. Therefore the algebraic groups $O^+(h)$ and $O^+(h')$ (respectively, $Spin(h)$ and $Spin(h')$) are identical.
\end{proof}
\end{cor}

\textbf{Proof of the main theorem in subcase 3.2.2:}

\begin{proof}
    
In subcase 3.2.2, in order to prove the main result for $Spin(h)$, it is enough to show it for $Spin(h')$ by Corollary \ref{101}. Since $h'$ is an unramified form, then the same argument as in subcase 3.2.1 works for this subcase, and we are done.
\end{proof}

\subsection{Proof of the main theorem: part 4}\label{159} In this subsection, we prove the main theorem in case 4: $D_K\cong (\frac{a,\pi}{K})$ is ramified over $K$ and $ L=K(\sqrt{t})$ is an unramified extension of $K$. Recall that the elements $a$ and $t$ are units. We have  $x^2=a, \ y^2=\pi,$ and  $z=xy=-yx$ for $x, y, z \in D_K$. 

\begin{lem}\label{163}
    Assume that $dim (h_0)\ne 0$ (recall that $h_0$ is the unramified part of the Larmour decomposition of the form $h$ over $K$). Then $D_L$ is an $L$-division algebra and ramified over $L$. 
    \begin{proof}
        If $D_L$ splits over $L$, then the unitary group $U(h_0)$ contains the norm one torus ${R_{K(x)/K}}^{(1)} (\mathbb{G}_m)$ which becomes a split torus over $L$ (the field $K(x)$ is isomorphic to $L$ over $K$; see Lemma \ref{94}), contradicting the fact that $h$ is anisotropic over $L$. Therefore $D_L$ is an $L-$division algebra. Since $L/K$ is unramified, $D_L$ is ramified as an algebra over $L$.
    \end{proof} 
\end{lem}

\begin{lem}\label{162}
In case 4, we have the following isometries over
L:

$$h_{0,L}\cong \lambda h_{0,L},   \  \   \ \textit{and}  \  \  \  \  \  h_{1,L}\cong \lambda h_{1,L}.$$

\begin{proof}
If one of the skew-hermitian forms $h_0$ or $h_1$ is trivial (over $K$), then we already have the desired result since $h_L \cong \lambda h_L$. So assume that both of them are nontrivial over $K$.

We know that $D_L$ is an $L$-division algebra by Lemma \ref{163}. Then we have $\lambda h_{0,L} \perp \lambda h_{1,L} = \lambda h_L \cong h_L = h_{0,L} \perp h_{1,L}$ and $\nu_L(\lambda)\in \{0,1\}$. Every element in the diagonalization of $h_0$ (resp. $h_1$) has value $0$ (resp. $\frac{1}{2}$) over $L$. By Theorem \ref{54} we have $h_{0,L}\cong \lambda h_{0,L}$ and   $h_{1,L}\cong \lambda h_{1,L}$. 
\end{proof}
\end{lem}

By Remark \ref{55} and Lemma \ref{162}, we may assume that either $h=h_0$ or $h=h_1$, i.e. $h$ is either unramified or ramified.\\

\begin{lem}\label{956}
    Assuem that $\lambda$ is a uniformizer of $L$. Let $\lambda ':=\lambda \pi^{-1}$.

    (1) There exists $u'\in Ker \ \alpha_L$ such that $g_L(u')=\lambda'$.

    (2) Let $u'$ be an element in $Ker \ \alpha_L$ such that $g_L(u')=\lambda'$ and $cor_{L/K}(u')\in Ker \ \alpha_K$. Then $cor_{L/K}(u)\in Ker \ \alpha_K$.
    \begin{proof}
        (1) Recall that $D_K$ is ramified over $K$ and $\pi$ is a square in $D_K$. So we have the isometry $h \cong \pi h$ over $K$ and $\pi$ is a multiplier for $h$ over $K$, hence  $\pi$ is also a multiplier for $h$ over $L$. Therefore  $\lambda':=\lambda \pi^{-1}$ is a unit multiplier of $h$ over $L$.

  Consider the following commutative diagram with exact rows:

\begin{diag}\label{500}
    
\[\
\begin{tikzcd}
 PGO^+(h)(-) \arrow[r,"\delta"] \arrow[d, "id"] & H^1(-,Z) \arrow[r, "\alpha"] \arrow[d, "g"]      & H^1(-,Spin(h)) \arrow[d, "s"]  \\
 PGO^+(h)(-) \arrow[r, "\gamma"]           & {H^1(-,\mu_2)} \arrow[r, "\beta"] & {H^1(-,O^+(h)).} 
\end{tikzcd}
\]\

\end{diag}

All elements $[\lambda],[\pi]$, and   $[\lambda']$ lie in the kernel of the map $\beta_L$. Since the leftmost vertical arrow in Diagram \ref{500} is the identity map, there exists $u'\in \Ker \ \alpha_L$ such that $g_L(u')=\lambda'$. \\

(2) Let $u'':=u {u'}^{-1}$. By Lemma \ref{85}, $u''\in \Ker \ \alpha_L$. Since $cor_{L/K}(u')\in\Ker  \ \alpha_K$, to prove $cor_{L/K}(u)\in \Ker \ \alpha_K$, it is sufficient to show that $cor_{L/K}(u'')\in\Ker  \ \alpha_K$.

 Since $\pi\in K^*$ is a multiplier for $h$ over $K$, by the exactness of the second row of Diagram \ref{500}, there exists an element $\kappa_K \in PGO^+(h)(K)$ such that $\gamma_K(\kappa_K)=[\pi]$. Let $\kappa_L$ be the image of $\kappa_K$ under the map $PGO^+(h)(K) \to PGO^+(h)(L)$. Let ${u_1}'':=\delta_L(\kappa_L)$ and ${u_2}'':=u'' {{u_1}''}^{-1}$. The elements $u''$ and ${u_1}''$ both lie in the kernel of the map $\alpha_L$, so by Lemma \ref{85}, the element ${u_2}''$ also lies in the kernel of the map $\alpha_L$. Since $\kappa_L \in PGO^+(h)(L)$ is induced from an element $\kappa_K \in PGO^+(h)(K)$ under the restriction map, so the element ${u_1}''$ belongs to the image of the restriction map $H^1(K,Z) \to H^1(L,Z)$, hence $cor_{L/K}({u_1}'')$ is trivial in $H^1(K, Z)$ by the restriction–corestriction formula in Galois cohomology. We need to show that $cor_{L/K}({u_2}'') \in \Ker \ \alpha_K$.

%Note that $H^1(L,Z)= H^1(L,\mu_2) \times H^1(L,\mu_2)$ and the map $g_L$ is given explicitly by the product of the two components (see Subsection \ref{154}). 

Since $g_L(u'')=g_L({u_1}'')=[\pi]$, we have ${u_2}'' \in \Ker \ g_L$. By the exactness of the sequence (see Diagram \ref{139})

$$\dots  \longrightarrow H^1(L,\mu_2) \overset{f_L}\longrightarrow H^1(L,Z) \overset{g_L}\longrightarrow H^1(L,\mu_2) \longrightarrow \dots ,$$there exists an element $w\in L^*$ such that ${u_2}''=f_L([w])$. The element $[w]$ is then a spinor norm for $h$ over $L$. By Knebusch's norm principle, $cor_{L/K}[w]=[N_{L/K}(w)]$ is a spinor norm for $h$ over $K$. Since (see Diagram \ref{139})

$$\alpha_K(cor_{L/K}({u_2}''))=e_K(cor_{L/K}([w])),$$we have $cor_{L/K}({u_2}'') \in \Ker \ \alpha_K$.
    \end{proof}
\end{lem}

\textbf{Subcase 4.1,  when $h$ is unramified.}

By Remark \ref{66} we have $$h=<\alpha_1 x, \dots, \alpha_n x>,$$

where $\alpha_1 , \dots, \alpha_n \in {\mathcal{O}_K}^*$. Recall that by Lemma \ref{163}, $D_L$ is an $L$-division algebra and ramified over $L$. 

The restriction of the involution $\tau$ to the field $K(x)$ is a unitary involution, which we denote by $\tau$ (abusing notation).
We define the following unitary unit hermitian form over $(K(x), \tau)$  $$h':= <\alpha_1 x, \dots, \alpha_n x>   \ \   \textit{over} \  \ (K(x), \tau). $$

The hermitian form $h'$ will only be used in subcase 4.1. 

\begin{lem}\label{78}
    The unitary group $U(h')$ is a  subgroup of the special orthogonal group $O^+(h)$. 
    \begin{proof}

Let $g\in U(h')(K(x))$, so once we fix a linear representation $U(h')(K(x)) \longrightarrow GL_n(K(x))$, we can view $g$ as a matrix, which belongs to $GL_n(K(x))$. In fact, $g$ is also an element of $GL_n(D)$ as $GL_n(K(x))\subset GL_n(D)$. Let $H$ be the diagonal matrix $diag(\alpha_1 x, \dots, \alpha_n x)$ representing $h'$ (and $h$). Also let $\sigma_h$ be the involution on $M_n(D)$ adjoint to $h$, and $\sigma_{h'}$ be the involution on $M_n(K(x))$ adjoint to $h'$. Both of these involutions are given by the map $M \mapsto H \tau (M^t) H^{-1}$, where $M^t$ means transpose of $M$, and by $\tau(M^t)$ we mean the matrix obtained by the action of $\tau$ on all entries of $M^t$. The element $g$ satisfies the equation $g \sigma_{h'} (g)=Id$. By taking the determinant of both sides we have $N_{K(x)/K} (det \  g)=1$.
So the element $g$ clearly lies in $O^+(h)$ as well, which shows that $U(h')\subset O^+(h)$. 
    \end{proof}
\end{lem}

Note that the unitary group $U(h')$ has a subgroup $SU(h')$, which in fact can be viewed as
a subgroup of $O^+(h)$ by Lemma \ref{78}. Let $H$ be the preimage of the group $SU(h')$ under the morphism $Spin(h) \longrightarrow O^+(h)$.

\begin{lem}\label{70}

Let $H^0$ be the connected component of $H$. Then the restriction of the map $Spin(h) \longrightarrow O^+(h)$ to $H^0$ is an isomorphism of algebraic groups $H^0 \cong SU(h')$.
    \begin{proof}
        
The connected component $H^0$ of $H$ has the structure  $H^0=R_u(H^0)\rtimes {H^0}_{red}$ over $K^{sep}$, where $R_u(H^0)$ is its unipotent radical and ${H^0}_{red}$ is a reductive subgroup. Under the morphism $H^0 \rightarrow SU(h')$ the unipotent radical of $H^0$ is being mapped to the unipotent radical of $SU(h')$, which is trivial because $SU(h')$ is a semisimple group. Hence $R_u(H^0) \subset \Ker \ (H^0 \rightarrow SU(h'))$, and since the kernel is finite, $R_u(H^0)$ is trivial. So the connected component  $H^0$ is reductive, and it is an almost direct product of a central torus $S$ and its semisimple part, i.e. $H^0=S\cdot[H^0, H^0]$. But then since $S$ is in the kernel of the morphism $H\rightarrow SU(h')$, it is trivial because the kernel is finite. Therefore, $H^0$ is semisimple. So the map $H^0 \longrightarrow SU(h')$ is isomorphism since $SU(h')$ is simply connected and does not admit any nontrivial covering.
    \end{proof}
\end{lem}

Let $Z_1=\mu_2$ be the central subgroup in $Spin(h)$ which is the kernel of the map $Spin(h)\longrightarrow O^+(h)$. Note that $Z_1 \cap H^0= \{1\}$ by Lemma \ref{70}.  

We need to consider the cases below separately based on the parity of $n$ (the dimension of $h$ over $D$):\\

\textbf{Subcase 4.1.1 }: $n$ is even.

%Note that the unitary group $U(h')$ has a subgroup $SU(h')$, which in fact can be viewed as a subgroup of $O^+(h)$ by Lemma \ref{78}. 

The center of the group $O^+(h)$ is $\mu_2$, which is a central subgroup of $SU(h')$ as well. Let us denote by $Z_2$ its preimage in $H^0$. On the other hand we know that the center of $Spin(h)$, i.e. $Z$, has the structure $R_{F/K}(\mu_2)=R_{(K\times K)/K}(\mu_2)=\mu_2 \times \mu_2$ (since $disc(h)$ is trivial we have $F=K\times K$; see Lemma \ref{67}). The intersection of  $Z_1$ and $Z_2$ in $Spin(h)$ is trivial, and $Z_2 \leq H^0$. We have $Z=Z_1 \times Z_2$.\\

\textbf{Proof of the main theorem in subcase 4.1.1}: 

\begin{proof}
\textbf{Case A}. 
Suppose $\lambda$ is a unit in $L^*$ up to squares. Consider the following diagram:

\[\ 
\begin{tikzcd}
1 \arrow[r] & Z_1 \arrow[r] \arrow[d, "id"] & Z=Z_1 \times Z_2 \arrow[r] \arrow[d]  & Z_2 \arrow[r] \arrow[d]    & 1 \\
1 \arrow[r] & Z_1 \arrow[r] \arrow[d, "id"] & Z_1 \times SU(h') \arrow[r] \arrow[d] & SU(h') \arrow[r] \arrow[d] & 1 \\
1 \arrow[r] & Z_1 \arrow[r]                 & Spin(h) \arrow[r]                     & O^+(h) \arrow[r]           & 1
\end{tikzcd}
\]\

Now we apply Lemma \ref{88}: Put\\
$G_1:=Z_1$\\
$G_2:=Z$\\
$G_3:=Z_2$\\
$G_4:= Z_1 \times SU(h')$\\
$G_5:= SU(h')$\\
$G_6:= Spin(h)$\\
$G_7:= O^+(h)$\\

We need to show that conditions (a),(b),(c), and (d) are satisfied:\\

(a) The element $\chi_L(g_L(u))\in H^1(L, SU(h'))$ corresponds to the isometry class of the unitary hermitian form $\lambda h'$ over $(L(x), \tau_{L(x)})$, where $\tau_{L(x)}$ denotes the nontrivial $L$-automorphism of $L(x)$. In order to show that $u\in \Ker \ \chi_L \circ g_L$, we need to prove that $\lambda h' \cong h'$, as unitary hermitian forms over $(L(x), \tau_{L(x)})$.

Since $\lambda$ is a multiplier for the skew-hermitian form $h$ over $(D_L, \tau_L)$, we have the following isometry of skew-hermitian forms over $(D_L,\tau_L)$: $$\lambda h \cong h.$$

Recall that by definition of $h'$,  the diagonalization of $h$ and $h'$ are same. By \cite[Theorem~3.7]{15} (also see \cite[Section~5, Case~B.1.2]{27} for an explicit formula), $\overline{\lambda}\ \overline{h'} \cong \overline{h'}$ as unitary hermitian forms over $(l(\overline{x}), \overline{\tau})$ (Recall that $l(\overline{x})$ is the residue $l$-algebra of $D_L$ and the residue involution $\overline{\tau}$ is the nontrivial $l$-automorphism of $l(\overline{x})$; see \cite[Page~22]{22}). Applying  \cite[Theorem~3.7]{15} (or \cite[Section~5, Case~B.1.2]{27}) again, we have $\lambda h' \cong h'$ as unitary hermitian forms over $(L(x), \tau_{L(x)})$.\\

(b) This is equivalent to the $H^1$-variant of the norm principle for the pair $(Z_2, SU(h'))$ which holds by \cite[main~theorem]{1}, as $SU(h')$ is an algebraic group of type $A_n$.\\

(c) It follows from Theorem \ref{8}.\\

(d) The map $\eta: H^1(-,Z_1) \times H^1(-,Z_2) \to H^1(-,Z_1) \times H^1(-,SU(h'))$ can be written as the product of the identity map $id: H^1(-, Z_1) \to H^1(-, Z_1)$ and the map $\chi: H^1(-,Z_2) \to H^1(-,SU(h'))$. Then the same argument as in part (b) above shows that $cor_{L/K} (\Ker (\eta_L))\subset \Ker (\eta_K)$.\\

\textbf{Case B.} Now, assume that $\lambda$ is a uniformizer of $L$. 

Let $u'$ be the element in $Ker \  \alpha_L$ introduced in the proof of part (1) of Lemma \ref{956}. By the proof of case A, we have $cor_{L/K}(u')\in Ker \ \alpha_K$. Therefore, by part (2) of Lemma \ref{956}, we have $cor_{L/K}(u)\in Ker \ \alpha_K$.
\end{proof}

\textbf{Subcase 4.1.2}: $n$ is odd. \\

Let $\Tilde{U}$ be the preimage of the group $U(h')$ under the morphism $Spin(h)\to O^+(h)$. Recall that $H$ is the preimage of the group $SU(h')$ under the morphism $Spin(h) \longrightarrow O^+(h)$.

\begin{lem}
The connected component   ${\Tilde{U}}^0$ is an almost direct product of $H^0$ (the connected component of $H$) and a central torus $S$. In particular, ${\Tilde{U}}^0$ is reductive. The torus $S$ is mapped to the norm one torus ${R_{K(x)/K}}^{(1)}(\mathbb{G}_m)$ under the morphism $Spin(h)\to O^+(h)$ with kernel $Z_1=\mu_2$. Furthermore, $Z\subset S$. 
    \begin{proof}
        
The connected group ${\Tilde{U}}^0$ has an almost direct product structure  ${\Tilde{U}}^0=R_u({\Tilde{U}}^0) \cdot {{\Tilde{U}}^0}_{red}$ over $K^{sep}$ where $R_u({\Tilde{U}}^0)$ is its unipotent radical and ${{\Tilde{U}}^0}_{red}$ is a reductive subgroup. Under the morphism ${\Tilde{U}}^0 \rightarrow U(h')$ the unipotent radical of $\Tilde{U}^0$ is being mapped to the unipotent radical of $U(h')$, which is trivial because $U(h')$ is a reductive group. Hence $R_u(\Tilde{U}^0) \subset \Ker \ (\Tilde{U}^0 \rightarrow U(h'))$, and since the kernel is finite, $R_u(\Tilde{U}^0)$ is trivial. So the connected component  $\Tilde{U}^0$ is reductive, and it is an almost direct product of a central torus $S$ and its semisimple part, i.e. ${\Tilde{U}}^0=S\cdot[{\Tilde{U}}^0, {\Tilde{U}}^0]$. 

On the other hand, the unitary group $U(h')$ is an almost direct product of the central  torus ${R_{K(x)/K}}^{(1)}(\mathbb{G}_m)$  and the special unitary group $SU(h')$, i.e. $U(h')={R_{K(x)/K}}^{(1)}(\mathbb{G}_m)\cdot SU(h')$.

Under the map ${\Tilde{U}}^0 \rightarrow U(h')$, the central torus $S$ goes to ${R_{K(x)/K}}^{(1)}(\mathbb{G}_m)$, and the semisimple part of ${\Tilde{U}}^0$, i.e. $[{\Tilde{U}}^0,{\Tilde{U}}^0]$, goes to the semisimple part of $U(h')$, i.e. $SU(h')$. Since the group $SU(h')$ is simply connected, the covering $[{\Tilde{U}}^0, {\Tilde{U}}^0]\rightarrow SU(h')$ has to be trivial, so $${\Tilde{U}}^0=S \cdot H^0.$$ 

The center of the group $O^+(h)$, which is $\mu_2$, is in fact a subgroup of ${R_{K(x)/K}}^{(1)}(\mathbb{G}_m)$: once we fix a linear representation $O^+(h) \to GL_n(D)$, the nontrivial element of the center of $O^+(h)$ is in fact the diagonal matrix $diag(-1, \dots, -1)$, whose determinant (i.e. $-1$) has reduced norm $1$, so this element belongs to ${R_{K(x)/K}}^{(1)}(\mathbb{G}_m)$. This implies that the preimage of   ${R_{K(x)/K}}^{(1)}(\mathbb{G}_m)$, i.e. $S$, contains $Z$ (the center of $Spin(h)$), and we have an exact sequence $$1 \to Z_1 \to S \to {R_{K(x)/K}}^{(1)}(\mathbb{G}_m) \to 1.$$
    \end{proof}
\end{lem}

\begin{lem}\label{212}

The multiplier $\lambda$ belongs to $N_{L(x)/L} (L(x)^*)$.
\begin{proof}
Since $n$ is odd, the determinant of $h$ is of the form $\alpha x^n$, for a unit $\alpha \in \mathcal{O}_K$ (recall that $h\cong h_0$). So the discriminant extension of $L$ is $L(x)$ ($L(x)\cong L[t]/<t^2 - Nrd_{{D_L}/L}(det \ h)>$). Recall that $Z=\Ker (Norm: R_{L(x)/L}(\mu_4) \longrightarrow \mu_4)$ (\cite[Page 332]{19}), and $H^1(L,Z)={U(L)}/{U_0(L)}$, where $U \subseteq \mathbb{G}_m \times \mathbb{G}_m$ is the subgroup defined by $$U(L):=\{(q,p) \in L^* \times L(x)^* | q^4=Norm_{L(x)/L}(p)\}$$ and $U_0 \subseteq U$ is the subgroup defined by $$U_0(L):=  \{(N_{L(x)/L}(p),p^4) | p \in L(x)^*\}.$$

The map $g_L: H^1(L, Z) \rightarrow H^1(L,\mu_2)$ sends $[q,p]$ to  $[N_{L(x)/L}(p_0)]$, where $p_0 \in L(x)^*$ is such that $p_0 \psi_L {(p_0)}^{-1}=q^{-2}p$, where $\psi_L$ is the nontrivial $L-$automorphism of $L(x)$ (see 
\cite[Proposition~13.36]{13}). Since $\lambda$ is in the image of the map $g_L$, then it belongs to $N_{L(x)/L} (L(x)^*)$.
\end{proof}
\end{lem}

\textbf{Proof of the main theorem in subcase 4.1.2}: 

\begin{proof}
    Consider the following diagram:
\[\ 
\begin{tikzcd}
1 \arrow[r] & Z_1 \arrow[r] \arrow[d, "id"] & Z \arrow[r] \arrow[d] & Z_2 \arrow[r] \arrow[d]                              & 1 \\
1 \arrow[r] & Z_1 \arrow[r] \arrow[d, "id"] & S \arrow[r] \arrow[d] & {R_{K(x)/K}}^{(1)}(\mathbb{G}_m) \arrow[r] \arrow[d] & 1 \\
1 \arrow[r] & Z_1 \arrow[r]                 & Spin(h) \arrow[r]     & O^+(h) \arrow[r]                                     & 1
\end{tikzcd}
\]\

Now we apply Lemma \ref{88}: Put\\
$G_1:=Z_1$\\
$G_2:=Z$\\
$G_3:=Z_2$\\
$G_4:= S$\\
$G_5:= {R_{K(x)/K}}^{(1)}(\mathbb{G}_m)$\\
$G_6:= Spin(h)$\\
$G_7:= O^+(h)$\\

We need to show that conditions (a),(b),(c), and (d) are satisfied:

(a) By Lemma \ref{212}, we know that the multiplier  $\lambda$ belongs to  $N_{L(x)/L}(L(x)^*)$. By Hilbert theorem 90, we have $H^1(L, {R_{L(x)/L}}^{(1)}(\mathbb{G}_m) )= {L^*}/{N_{L(x)/L}(L(x)^*)}$, hence $\chi_L([\lambda])=1$. Therefore, $u \in \Ker (\chi_L \circ g_L)$.

(b) This is true due to commutativity of ${R_{K(x)/K}}^{(1)}(\mathbb{G}_m)$.

(c) It follows from Theorem \ref{8}.

(d) It holds because of commutativity of $S$.\\
\end{proof}

\textbf{Subcase 4.2, when $h=h_1$ is ramified.}  

By Remark \ref{66} we can assume that $$h=<(\beta_1 + \gamma_1 x)y, \dots, (\beta_n + \gamma_n x)y>,$$

where $\beta , \dots, \beta_n, \gamma_1, \dots, \gamma_n \in {\mathcal{O}_K}$ and the elements $\beta_i + \gamma_i x$ are all units.

Consider the following orthogonal involution on $D$:

$$\tau_{y}: D   \longrightarrow D$$

$$ \   \  \   \  \   \   \   \  \   \  \  \   \  \  \      s \mapsto  y^{-1} \tau(s) y.$$

Elements in the subfield $K(x)$ are invariant under $\tau_y$. Then
$$h':=<\beta_1 + \gamma_1 x, \dots, \beta_n + \gamma_n x>$$
is a unit quadratic form over $K(x)$. The following is a unit hermitian form over $(D, \tau_y)$ (since the elements of the maximal subfield $K(x)$ are invariant under $\tau_y$):

$$h''=<\beta_1 + \gamma_1 x, \dots, \beta_n + \gamma_n x>.$$

The quadratic form $h'$ is in fact the restriction of the hermitian form $h''$ to $K(x)^n$.

Let $H=diag(\beta_1 + \gamma_1 x, \dots, \beta_n + \gamma_n x)$ be the diagonal matrix representing $h'$ (or $h''$). Then the skew-symmetric matrix $Hy$ (with respect to the canonical involution $\tau$) represents $h$.

Let $\sigma_h$ and $\sigma_{h''}$ be the corresponding adjoint involutions on $M_n (D)$. Then for all $M\in M_n(D)$ we have

 \begin{equation*} 
\begin{split}
\sigma_{h''}(M) & =  H \tau_y(M^t)H^{-1} \\
 & = Hy \tau(M^t)y^{-1}H^{-1} \\
 & = Hy \tau(M^t) (Hy)^{-1} \\
 & = \sigma_h(M). 
\end{split}
\end{equation*}

Hence the algebraic groups $O^+(h)$ and $O^+(h'')$ (respectively, $Spin(h)$ and $Spin(h'')$) are identical. So it is enough to show the norm principle for the pair $(Z, Spin(h''))$ over $L/K$ (abusing notation, we denote the center of $Spin(h'')$ by $Z$ as well).\\

Note that the forms $h'$ and $h''$ will only be used in subcase 4.2.

Let $G:=R_{K(x)/K} (O^+ (h'))$.

\begin{lem}\label{92}
   $G$ is a  subgroup of $O^+(h'')$. 
    \begin{proof}

Let $g\in G(K)$. By the natural identification $G(K)=O^+ (h')(K(x))$ and abusing notation, we view $g$ as an element in $O^+(h')(K(x))$. Once we fix a linear representation $O^+(h') \longrightarrow GL_n(K(x))$, we can view $g$ as a matrix, which belongs to $GL_n(K(x))$. In fact, $g$ is also an element of $GL_n(D)$ as $GL_n(K(x))\subset GL_n(D)$. Recall that the matrix $H$ represents both $h'$ and $h''$. Furthermore $\sigma_{h''}$ is the involution on $M_n(D)$ adjoint to $h''$,  and  $\sigma_{h'}$ is the involution on $M_n(K(x))$ adjoint to $h'$. The involution $\sigma_{h''}$ is given by the map $M \mapsto H {\tau_y} (M^t) H^{-1}$, and the involution $\sigma_{h'}$ is the restriction of $\sigma_{h''}$ to $M_n(K(x))$ (so it is given by the same formula).  The element $g$ satisfies the equation $g \sigma_{h'} (g)=Id$ and $N_{K(x)/K} (det \  g)=1$. Clearly, $g$ lies in the group $O^+(h'')$ since

$$O^+(h''):=\{f \in GL_n(D)| f \sigma_{h''} (f)=Id, \  \textit{and} \ Nrd(f)=1\}.$$
    \end{proof}
\end{lem}

%Let $o_K$ be the induced map $o_K: H^1(K, G) \rightarrow H^1(K, O^+(h''))$. 

%\begin{lem}\label{89}
 %   Kernel of the map $o_K$ is trivial. Furthermore, if $D_L$ is an $L-$division algebra, then kernel of the map $o_L$ is also trivial.
  %  \begin{proof}
   %     Let $[\mathfrak{h}]\in \ \Ker (o_K)$. By the identification $H^1(K,G)=H^1(K(x), O^+ (h'))$, we can view  $\mathfrak{h}$ as a quadratic form over $K(x)$ which is isometric to the hermitian form $h''$ after extending the scalars from $K(x)$ to $D$. This in particular means that $\mathfrak{h}$ is a unit quadratic form by \cite[Theorem~3.7]{15}. We have $\overline{\mathfrak{h}} \cong \overline{h''}$ as unitary hermitian forms over $(k(\overline{x}), \overline{\tau_y})$. 
  %      Note that the field $k(\bar{x})$ is element-wise invariant under the involution $\overline{\tau_y}$, so we can view both $\overline{\mathfrak{h}}$ and $ \overline{h''}$ as quadratic forms over $k(\overline{x})$. 

     %   By Springer's decomposition theorem for quadratic forms, we have $\mathfrak{h} \cong h'$  over $K(x)$. Therefore $\Ker (o_K)$ is trivial. 
        
   %     A similar argument works for $o_L$, if $D_L$ is an $L-$division algebra. 
  %  \end{proof}
%\end{lem}

Since $G=R_{K(x)/K}(O^+ (h'))$, we have $C(G)=R_{K(x)/K}(C(O^+(h')))=R_{K(x)/K} (\mu_2)$ and 
$|C(G)|=4$. Also $G(K(x))=O^+(h')(K(x)) \times O^+(\widetilde{h'})(K(x))$ (direct product), where $\widetilde{h'}=<\beta_1 - \gamma_1 x, \dots, \beta_n - \gamma_n x>$. Over the algebraic closure $G$ is isomorphic to the direct product of two copies of ${O_n}^+$ and the center of each of these components is $\mu_2$.

Let $G'$ be the preimage of the group $G$ under the morphism $Spin(h'') \longrightarrow O^+(h'')$:

\[\begin{tikzcd}
            & 1 \arrow[d]                     & 1 \arrow[d]       \\
            & \mu_2 \arrow[d] \arrow[r, "id"] & \mu_2 \arrow[d]   \\
1 \arrow[r] & G' \arrow[r] \arrow[d]           & Spin(h'') \arrow[d] \\
1 \arrow[r] & G \arrow[r] \arrow[d]      & O^+(h'') \arrow[d]   \\
            & 1                               & 1                
\end{tikzcd}\]

Over the algebraic closure, we denote by ${O_m}^+$, $Spin_m$, and $C_{m}$ the special orthogonal group, the spinor group, and the even Clifford algebra of any $m-$ dimensional quadratic form, respectively (for any $m\in \mathbb{N}$). We have $G(\overline{K})\cong {O_n}^+ (\overline{K}) \times {O_n}^+ (\overline{K})$. The preimage of each component ${O_n}^+$ of $G$ under the map $Spin_{2n}\to {O_{2n}}^+$ is $Spin_n$. So there is a surjective map $Spin_n \times Spin_n \to G'$ which means that $G'$ is an almost direct product of two copies of $Spin_n$, i.e. $G'=Spin_n \cdot Spin_n$. The kernel of the map  $Spin_n \times Spin_n \to G'$ consist of all pairs $(w, w^{-1})$ where $w$ is an element in the intersection of the two components $Spin_n$. Let us view such an element $(w, w^{-1})$ in $C_{2n}$. The element $w$ then belongs to each component $C_n$ of the product $C_n \cdot C_n$ inside $C_{2n}$, but we know that the two copies of $C_n$ have only the scalars in common. Furthermore, based on the  explicit realization of elements in the even Clifford algebra (see 
\cite[Section~2.2]{3}), the only scalars which belong to each component $Spin_n$ are the elements in $\mu_2$. In conclusion, the intersection of the two components $Spin_n$ in the almost direct product $Spin_n \cdot Spin_n$  is $\mu_2$. The center of $G'$ is isomorphic to the almost direct product of $C(Spin_n)$ with itself, the intersection being $\mu_2$. Hence $C(G')$ is an abelian group of order $8$ (Note that $|C(Spin_n)|=4$) containing $Z=C(Spin(h''))$. Finally we get the following diagram, whose vertical arrows are inclusion and the rows are exact:

\[\
\begin{tikzcd}
1 \arrow[r] & \mu_2 \arrow[r] \arrow[d, "id"] & Z \arrow[r] \arrow[d]  & \mu_2 \arrow[r] \arrow[d] & 1 \\
1 \arrow[r] & \mu_2 \arrow[r] \arrow[d, "id"] & G' \arrow[r] \arrow[d] & G \arrow[r] \arrow[d]     & 1 \\
1 \arrow[r] & \mu_2 \arrow[r]                 & Spin(h'') \arrow[r]    & O^+(h'') \arrow[r]        & 1
\end{tikzcd}
\]\

Note that $G$ and $G'$ are both semisimple (in particular connected), because of their structure as direct or almost direct products of semisimple  groups over the algebraic closure. Also the inclusion $\mu_2 \to G$ is induced via the diagonal embedding of $\mu_2$ into the center of $G$ (over the algebraic closure), i.e. $$\mu_2 \to R_{K(x)/K} (\mu_2)=C(G) \to G.$$

\begin{lem}\label{165}
    The $H^1$-variant of the norm principle holds for the pair $(Z,G')$ over $L/K$.
    \begin{proof}
Let 
$$\overline{h'}:=<\overline{\beta_1 + \gamma_1 x}, \dots, \overline{\beta_n + \gamma_n x}>$$
be the residue quadratic form over $k(\overline{x})$. By our assumption, the $H^1$-variant of the norm principle holds for the pair $(C(Spin(\overline{h'})), Spin(\overline{h'}))$ over $k(\overline{x})$. By 
\cite[Theorem~5.1]{4}, the $H^1$-variant of the norm principle holds for the pair $(C(Spin(h')), Spin(h'))$ over $K(x)$. Then the proof of Theorem \ref{113} (2$\implies$1) shows  that the $H^0$-variant of the norm principle holds for the group $\Omega(h')$  (see \cite[Section~23.B]{13} for the definition of the extended Clifford group $\Omega$) over $K(x)$; this is obtained by assuming that the reductive and the semisimple groups in the proof of Theorem $\ref{113}$ are $\Omega(h')$ and $Spin(h')$, respectively. By 
\cite[Lemma~2.6]{1}, the $H^0$-variant of the norm principle holds for the group $R_{K(x)/K}(\Omega(h'))$ over $K$. The same argument as in the proof of Theorem \ref{113} (1$\implies$2) shows that the $H^1$-variant of the norm principle holds for the pair $(C(R_{K(x)/K}(Spin(h'))), R_{K(x)/K}( Spin(h')))$ over $K$. Note that the group $R_{K(x)/K}(Spin(h'))$ is the commutator subgroup of $R_{K(x)/K}(\Omega(h'))$:

$$R_{K(x)/K}(Spin(h'))= R_{K(x)/K}([\Omega(h'), \Omega(h')])=[R_{K(x)/K}(\Omega(h')), R_{K(x)/K}(\Omega(h'))].$$

By Lemma \ref{31}, the $H^1$-variant of the norm principle holds for the pair $(C(G'), G')$ over $K$, since $R_{K(x)/K}(Spin(h'))$ is the simply connected cover of $G'$ (this follows from the covering $Spin(h'') \to O^+(h'')$). Since $C(G')$ contains $Z$, by applying Lemma \ref{37} we obtain the $H^1$-variant of the norm principle for the pair $(Z,G')$ over $K$, in particular over $L/K$.
    \end{proof}
\end{lem}

\begin{lem}\label{167}
    The $H^1$-variant of the norm principle holds for the pair $(\mu_2,G)$ over $L/K$.
    \begin{proof}
        
 Scharlau's norm principle for quadratic forms states that the $H^0$-variant of the norm principle holds for the group $GO(h')$ over $K(x)$. By \cite[Lemma~2.6]{1}, the $H^0$-variant of the norm principle holds for the group $R_{K(x)/K}(GO(h'))$ over $K$. 

Repeating the argument in the proof of Theorem $\ref{113}$ (by assuming that the reductive group and the semisimple group in the proof of Theorem $\ref{113}$ are $R_{K(x)/K}(GO(h'))$ and $G=R_{K(x)/K}(O^+(h'))$ respectively), we obtain the $H^1$-variant of the norm principle for the pair $(C(G), G)$ over $K$; note that the group $G$ is the commutator subgroup of $R_{K(x)/K}(GO(h'))$:
$$G=R_{K(x)/K}(O^+(h'))= R_{K(x)/K}([GO(h'), GO(h')])=[R_{K(x)/K}(GO(h')), R_{K(x)/K}(GO(h'))].$$ Since $C(G)$  contains $\mu_2=C(O^+(h''))$, by applying Lemma \ref{37} we conclude that the $H^1$-variant of the norm principle holds for the pair $(\mu_2,G)$ over $K$, in particular over $L/K$.  \end{proof}
\end{lem}

We split the rest of the argument is case 4.2 into two subcases:\\

\textbf{Subcase 4.2.1, when $D_L$ is an $L-$division algebra.}\\

\textbf{Proof of the main theorem in subcase 4.2.1:}\\

\begin{proof}
    Recall that it is enough to show the norm principle for the pair $(Z, Spin(h''))$ over $L/K$.

    \textbf{Case A}: Assume that $\lambda$ is a unit in $L^*$ up to squares. 
    
    Consider the following diagram:

\[\
\begin{tikzcd}
1 \arrow[r] & \mu_2 \arrow[r] \arrow[d, "id"] & Z \arrow[r] \arrow[d]  & \mu_2 \arrow[r] \arrow[d] & 1 \\
1 \arrow[r] & \mu_2 \arrow[r] \arrow[d, "id"] & G' \arrow[r] \arrow[d] & G \arrow[r] \arrow[d]     & 1 \\
1 \arrow[r] & \mu_2 \arrow[r]                 & Spin(h'') \arrow[r]    & O^+(h'') \arrow[r]        & 1
\end{tikzcd}
\]\

Now we apply Lemma \ref{88}: Put\\
$G_1:=\mu_2$\\
$G_2:=Z$\\
$G_3:=\mu_2$\\
$G_4:= G'$\\
$G_5:= G$\\
$G_6:= Spin(h'')$\\
$G_7:= O^+(h'')$

We need to show that conditions (a), (b), (c), and (d) are satisfied:

(a) The element $\chi_L(g_L(u))\in H^1(L, G)= H^1(L(x), O^+(h'))$ corresponds to the isometry class of the quadratic  form $\lambda h'$ over $L(x)$. In order to show that $u\in \Ker \ \chi_L \circ g_L$, we need to prove that $\lambda h' \cong h'$, as quadratic forms over $L(x)$.

Since $\lambda$ is a multiplier for the skew-hermitian form $h$ over $(D_L, \tau_L)$, we have the following isometry of skew-hermitian forms over $(D_L,\tau_L)$: $$\lambda h \cong h.$$

Recall that according to the diagonalization of $h$ and  $h'$, and by applying  \cite[Theorem~3.7]{15} (see also \cite[Section~5, Case~B.1.2]{27} for an explicit formula), we obtain the following isometry of quadratic forms over $l(\overline{x})$:

$$\overline{\lambda} \overline{h'} \cong \overline{h'}.$$

By Springer's theorem for quadratic forms, we have $\lambda h' \cong h'$ as quadratic forms over $L(x)$.

(b) This was proved in Lemma \ref{167}.

(c) It is equivalent to Knebusch's norm principle for quadratic forms (see Theorem \ref{147}).

(d) This was proved in Lemma \ref{165}.\\

\textbf{Case B}: Assume that $\lambda$ is a uniformizer of $L^*$.

Let $u'$ be the element in $Ker \  \alpha_L$ introduced in the proof of part (1) of Lemma \ref{956}. By the proof of case A, we have $cor_{L/K}(u')\in Ker \ \alpha_K$. Therefore, by part (2) of Lemma \ref{956}, we have $cor_{L/K}(u)\in Ker \ \alpha_K$.

\end{proof}

\textbf{Subcase 4.2.2, when $D_L$ is split.}\\

When $D_L$ splits then $L\cong K(x)$, the unique unramified field extension of $K$ in $D_K$ (see Lemma \ref{94}). Once we fix an isomorphism $D_L\cong M_2(L)$, we have a Morita correspondence between $m$ dimensional hermitian forms over $D_L$ (with respect to the involution $\tau_y$) and $2m$ dimensional quadratic forms over $L$, up to scalars in $L$, for any $m\in \mathbb{N}$. Recall that $$h''=<\beta_1 + \gamma_1 x, \dots, \beta_n + \gamma_n x>.$$ 

Under Morita equivalence, the one dimensional hermitian form $<\beta_i + \gamma_i x>$ corresponds to the quadratic form $\theta<\beta_i + \gamma_i x, -(\beta_i - \gamma_i x)\pi>$ for a scalar $\theta \in L$. Since Morita equivalence respects direct sums, the form $h''$ corresponds to the quadratic form $\mathfrak{h}:=\theta h' \perp (-\theta \pi \widetilde{h'})$ over $L$. The element $\lambda$ is a multiplier for $h$ over $L$, so it is a multiplier for $h''$ over $L$ as well. The hermitian form $h''\perp (-\lambda h'')$ is hyperbolic over $L$, hence the quadratic form $\mathfrak{h}\perp (-\lambda \mathfrak{h})$ is hyperbolic over $L$ (again, because Morita equivalence respects scalar products, orthogonal sums, and hyperbolicity). So $\lambda$ is a multiplier for $\mathfrak{h}$ over $L$, and therefore it is also a multiplier for the quadratic form $\theta^{-1}\mathfrak{h}=h' \perp (-\pi \widetilde{h'})$. We have the $L$-isometry of quadratic forms $${h_L}' \perp (-\pi \widetilde{{h_L}'})\cong \lambda {h_L}' \perp (-\lambda \pi \widetilde{{h_L}'}).$$

We split the rest of the proof in subcase 4.2.2 into two parts: when $\lambda$ is a uniformizer of $L$, and then the case that $\lambda$ is a unit in $L$. \\

\textbf{Subcase 4.2.2.1, when $\lambda$ is a uniformizer of $L$.}\\

\textbf{Proof of the main theorem in subcase 4.2.2.1:}

\begin{proof}
    By Springer's decomposition theorem we have $-\pi \widetilde{{h_L}'}\cong \lambda {h_L}'$, hence $\mathfrak{h}_L\cong <1,\lambda>\otimes (\theta {h_L}')$ and $PGO^+(\mathfrak{h})(L)$ is $R$-trivial by 
    \cite[Lemma~4.2]{4}. Therefore $PGO^+(h'')(L)$ and $PGO^+(h)(L)$ are also $R-$trivial and we are done by Theorem \ref{7}. \\
\end{proof}

\textbf{Subcase 4.2.2.2, when $\lambda$ is a unit in $L$.}\\

\textbf{Proof of the main theorem in subcase 4.2.2.2:}

\begin{proof}
    By Springer's decomposition theorem we have the following $L$ isometries of quadratic forms:
$${h_L}'\cong \lambda  {h_L}'\ \ \ \ \  \textit{and} \ \ \ \ \  
 \widetilde{{h_L}'}\cong \lambda \widetilde{{h_L}'}.\ \ \ \ \ \ $$

Recall that it is enough to show the norm principle for the pair $(Z, Spin(h''))$ over $L/K$. Consider the following diagram:

\[\
\begin{tikzcd}
1 \arrow[r] & \mu_2 \arrow[r] \arrow[d, "id"] & Z \arrow[r] \arrow[d]  & \mu_2 \arrow[r] \arrow[d] & 1 \\
1 \arrow[r] & \mu_2 \arrow[r] \arrow[d, "id"] & G' \arrow[r] \arrow[d] & G \arrow[r] \arrow[d]     & 1 \\
1 \arrow[r] & \mu_2 \arrow[r]                 & Spin(h'') \arrow[r]    & O^+(h'') \arrow[r]        & 1
\end{tikzcd}
\]\

Now we apply Lemma \ref{88}: Put\\
$G_1:=\mu_2$\\
$G_2:=Z$\\
$G_3:=\mu_2$\\
$G_4:= G'$\\
$G_5:= G$\\
$G_6:= Spin(h'')$\\
$G_7:= O^+(h'')$

We need to show that conditions (a),(b),(c), and (d) are satisfied:

(a) The image of $[\lambda]$ under the following map

$$\chi_L: H^1(L,\mu_2) \rightarrow H^1(L,G)=H^1(L,O^+(h')) \times H^1(L,O^+(\widetilde{h'}))$$is the pair $([\lambda h'], [\lambda \widetilde{h'}])=([h'], [\widetilde{h'}])=1\in H^1(L,G)$. So $u \in \Ker (\chi_L \circ g_L)$.

(b) This was proved in Lemma \ref{167}.

(c) It is equivalent to Knebusch's norm principle for quadratic forms (see Theorem \ref{147}).

(d) This was proved in Lemma \ref{165}.\\
\end{proof}


\begin{thebibliography}{90}

\bibitem{1} P. Barquero and A. Merkurjev, Norm principle for reductive algebraic groups, Algebra,
arithmetic and geometry, Part I, II (Mumbai, 2000), Tata Inst. Fund. Res. Stud.
Math., vol. 16, Tata Inst. Fund. Res., Bombay, 2002, pp. 123–137. MR 1940665
\bibitem{28} Nivedita Bhaskhar, On Serre’s injectivity question and norm principle, Comment.
Math. Helv. 91 (2016), no. 1, 145–161, DOI 10.4171/CMH/381. MR3471940
\bibitem{3} Nivedita Bhaskhar, R-equivalence and norm principles in algebraic groups, ProQuest
LLC, Ann Arbor, MI, 2016, Thesis (Ph.D.)–Emory University. MR 3593322
\bibitem{4} Nivedita Bhaskhar, Vladimir Chernousov, and Alexander Merkurjev, The norm principle for type $D_n$ groups over complete discretely valued fields, Trans. Amer. Math.
Soc. 372 (2019), no. 1, 97–117. MR 3968764
\bibitem{5} Armand Borel, Linear algebraic groups, second ed., Graduate Texts in Mathematics,
vol. 126, Springer-Verlag, New York, 1991. MR 1102012


\bibitem{1001} S. Bosch, W. Lutkebohmert, M. Raynaud, Neron Models, Ergeb. Math.
Grenz., 21, Springer, New York–Heidelberg–Berlin, 1990.



\bibitem{6} V. Chernousov and A. Merkurjev, R-equivalence and special unitary groups, J. Algebra
209 (1998), no. 1, 175–198. MR 1652122

\bibitem{235} M. Demazure, Sous-groupes paraboliques des groupes réductifs, in Schémas en groupes (SGA 3), Séminaire de Géométrie Algébrique du Bois Marie 1962–1964, Lecture Notes in Mathematics, vol. 153, Springer-Verlag, 1970, Exposé XXVI.




\bibitem{7} Antonio J. Engler and Alexander Prestel, Valued fields, Springer Monographs in Mathematics, Springer-Verlag, Berlin, 2005. MR 2183496
\bibitem{29} Philippe Gille, R-´equivalence et principe de norme en cohomologie galoisienne
(French, with English and French summaries), C. R. Acad. Sci. Paris S´er. I Math.
316 (1993), no. 4, 315–320. MR1204296
\bibitem{8} Philippe Gille, La R-équivalence sur les groupes algébriques réductifs définis sur un
corps global, Inst. Hautes Études Sci. Publ. Math. (1997), no. 86, 199–235 (1998). MR
1608570
\bibitem{9} Philippe Gille and Tamás Szamuely, Central simple algebras and Galois cohomology,
Cambridge Studies in Advanced Mathematics, vol. 101, Cambridge University Press,
Cambridge, 2006. MR 2266528

\bibitem{25} G. Harder, Uber die Galoiskohomologie halbeinfacher algebraischer Grup- ¨
pen. III (German), J. Reine Angew. Math. 274/275 (1975), 125–138, DOI
10.1515/crll.1975.274-275.125. Collection of articles dedicated to Helmut Hasse on
his seventy-fifth birthday, III. MR0382469
\bibitem{10} James E. Humphreys, Linear algebraic groups, Graduate Texts in Mathematics, No.
21, Springer-Verlag, New York-Heidelberg, 1975. MR 396773
\bibitem{11} Nikita A. Karpenko, Hyperbolicity of orthogonal involutions, Doc. Math. (2010),
no. Extra vol.: Andrei A. Suslin sixtieth birthday, 371–392, With an appendix by
Jean-Pierre Tignol. MR 2804259
\bibitem{12} Max-Albert Knus, Quadratic and Hermitian forms over rings, Grundlehren der mathematischen Wissenschaften [Fundamental Principles of Mathematical Sciences], vol.
294, Springer-Verlag, Berlin, 1991, With a foreword by I. Bertuccioni. MR 1096299


\bibitem{26} Martin Kneser, Galois-Kohomologie halbeinfacher algebraischer Gruppen ¨uber
p-adischen K¨orpern. II (German), Math. Z. 89 (1965), 250–272, DOI
10.1007/BF02116869. MR0188219 

\bibitem{13} Max-Albert Knus, Alexander Merkurjev, Markus Rost, and Jean-Pierre Tignol, The
book of involutions, American Mathematical Society Colloquium Publications, vol. 44,
American Mathematical Society, Providence, RI, 1998, With a preface in French by
J. Tits. MR 1632779


\bibitem{1000} Amit Kulshrestha and R. Parimala, R-equivalence in adjoint classical groups over
fields of virtual cohomological dimension 2, Trans. Amer. Math. Soc. 360 (2008),
no. 3, 1193–1221, DOI 10.1090/S0002-9947-07-04300-0. MR2357694



\bibitem{14} T. Y. Lam, Introduction to quadratic forms over fields, Graduate Studies in Mathematics, vol. 67, American Mathematical Society, Providence, RI, 2005. MR 2104929
\bibitem{15} Douglas W. Larmour, A Springer theorem for Hermitian forms, Math. Z. 252 (2006),
no. 3, 459–472. MR 2207754



\bibitem{30}A. S. Merkurjev, The norm principle for algebraic groups (Russian, with Russian
summary), Algebra i Analiz 7 (1995), no. 2, 77–105; English transl., St. Petersburg
Math. J. 7 (1996), no. 2, 243–264. MR1347513
\bibitem{16} A. S. Merkurjev, R-equivalence and rationality problem for semisimple adjoint classical
algebraic groups, Inst. Hautes Études Sci. Publ. Math. (1996), no. 84, 189–213 (1997).
MR 1441008
\bibitem{17} Yevsey A. Nisnevich, Espaces homogènes principaux rationnellement triviaux et arithmétique des schémas en groupes réductifs sur les anneaux de Dedekind, C. R. Acad.
Sci. Paris Sér. I Math. 299 (1984), no. 1, 5–8. MR 756297
\bibitem{18} Alex Ondrus, Minimal anisotropic groups of higher real rank, Michigan Math. J. 60
(2011), no. 2, 355–397, With an appendix by V. Chernousov and A. Merkurjev. MR
2825267
\bibitem{19} Vladimir Platonov and Andrei Rapinchuk, Algebraic groups and number theory, Pure
and Applied Mathematics, vol. 139, Academic Press, Inc., Boston, MA, (1994). Translated from the 1991 Russian original by Rachel Rowen. MR1278263
\bibitem{20} Winfried Scharlau, Quadratic and Hermitian forms, Grundlehren der mathematischen
Wissenschaften [Fundamental Principles of Mathematical Sciences], vol. 270, SpringerVerlag, Berlin, 1985. MR 770063
\bibitem{21} Jean-Pierre Serre, Galois cohomology, english ed., Springer Monographs in Mathematics, Springer-Verlag, Berlin, 2002, Translated from the French by Patrick Ion and
revised by the author. MR 1867431
\bibitem{27} Amin Soofiani, An explicit formula for Larmour's decomposition of hermitian forms, preprint available online at https://arxiv.org/pdf/2412.11110
\bibitem{270} Amin Soofiani, Hensel’s lemma for the norm principle for type $D_n$ groups, PhD thesis, University of British Columbia, 2025. https://dx.doi.org/10.14288/1.0448236
\bibitem{22} Jean-Pierre Tignol and Adrian R. Wadsworth, Value functions on simple algebras,
and associated graded rings, Springer Monographs in Mathematics, Springer, Cham,
2015. MR 3328410
\bibitem{23} Adrian R. Wadsworth, Extending valuations to finite-dimensional division algebras,
Proc. Amer. Math. Soc. 98 (1986), no. 1, 20–22. MR 848866



\bibitem{24} William C. Waterhouse, Introduction to affine group schemes, Graduate Texts in
Mathematics, vol. 66, Springer-Verlag, New York-Berlin, 1979. MR 547117






\end{thebibliography}
\end{document}